%% file: main.tex
\documentclass[twoside,11pt]{article}

\usepackage{blindtext}
\usepackage{microtype}
\usepackage{graphicx}
\usepackage{subcaption}
\usepackage{booktabs} 
\usepackage{amssymb} 
\usepackage{soul}
\usepackage{tcolorbox}
\usepackage[dvipsnames]{xcolor}

\usepackage{xcite}       
\usepackage{url}
\usepackage{booktabs}
\usepackage{array}
\usepackage{amsfonts}
\usepackage{nicefrac}
\usepackage{microtype}
\usepackage{amsmath}
\usepackage{enumitem}
\usepackage{multirow}
\usepackage{tablefootnote}
\usepackage{bm} 
\usepackage{bbm}
\usepackage{amssymb}
\usepackage{xspace}
\usepackage{tikz}
\usepackage[colorinlistoftodos]{todonotes}
\usepackage{pgfplots}
\usepgfplotslibrary{groupplots}
\pgfplotsset{compat=1.17}
\usepackage{graphicx}
\usepackage{subcaption}
\usepackage{float}
\usepackage{wrapfig}
\usepackage[normalem]{ulem}
\usepackage{pdfpages}
\usepackage{mathtools}
\usepackage{algorithm}
\usepackage{algorithmic}

\usepackage[preprint]{jmlr2e}
\usepackage{cleveref}

\input{macros}

\usepackage{lastpage}
\jmlrheading{23}{2026}{1-\pageref{LastPage}}{1/21; Revised 5/22}{9/22}{21-0000}{Yixuan Florence Wu, Yilun Zhu, Naichen Shi}

\ShortHeadings{JK-EGW for Multimodal Alignment}{Wu, Zhu, and Shi}
\firstpageno{1}

\begin{document}

\title{Multimodal Alignment Through Joint Kernel Entropic Gromov--Wasserstein Optimal Transport}

\author{\name Yixuan Florence Wu \email yixuanwu@u.northwestern.edu \\
       \addr Department of Engineering Sciences and Applied Mathematics \\
       Northwestern University \\
       Evanston, IL 60208, USA
       \AND
       \name Yilun Zhu \email allanzhu@umich.edu \\
       \addr Deparment of Electrical Engineering and Computer Science \\
       University of Michigan\\
       Ann Arbor, MI 48109, USA
       \AND
       \name Naichen Shi \email naichen.shi@northwestern.edu \\
       \addr Department of Industrial Engineering and Management Sciences \\
       Department of Mechanical Engineering \\
       Northwestern University \\
       Evanston, IL 60208, USA
       }

\editor{My editor}

\maketitle

\begin{abstract}
We study the problem of aligning data from multiple modalities into a shared representation space, focusing on settings where strong pretrained unimodal encoders are available but cross-modal paired data are scarce. We propose a structure-preserving alignment framework, joint kernel entropic Gromov--Wasserstein Optimal Transport (JK-EGW), which maps multiple modalities into a common latent space by minimizing a quadratic optimal transport objective. JK-EGW leverages fine-grained similarity relationships within and across modalities to construct a global affinity kernel instead of relying on raw feature-space distances. Our framework naturally provides explicit control over the geometry and distribution of the latent embedding. On the theory side, we establish parametric sample complexity rate of $n^{-1/2}$, matching the corresponding rates for standard, entropic and Gromov--Wasserstein optimal transport. On the algorithmic side, we derive a scalable alternating procedure to solve JK-EGW with entropic optimal transport (EOT) updates through a low-rank kernel approximation and a variational lifting. This lifting scheme effectively relieves the burden of a quadratic objective, and allowing us to take the advantage of existing EOT solvers. Empirically, we focus on post-hoc alignment of embeddings from pretrained encoders in data-scarce regimes, and show that our proposed method achieves improved multimodal retrieval performance compared to existing alignment baselines.
\end{abstract}



\section{Introduction}
Multimodal learning aims to integrate information from heterogeneous data sources, such as images, text, audio, and sensor measurements. Since different modalities are often represented in distinct feature spaces, a fundamental challenge is to construct a shared representation in which semantically related observations from different modalities become comparable~\citep{wang2022git,liu2023visual,zhao2025qlip}. Such aligned representations provide a common interface for downstream tasks~\citep{liang2024foundations}, including cross-modal retrieval~\citep{karpathy2015deep,faghri2017vse++}, zero-shot classification~\citep{frome2013devise,clip}, multimodal prediction~\citep{liang2021multibench}, and representation transfer~\citep{moschella2022relative}, etc. Developing principled methods for multimodal alignment is therefore a basic step toward learning systems that can use complementary information across modalities.

A large number of techniques focus on pretraining joint representation encoders via auto-encoders~\citep{suzuki2016joint,sutter2024unity}, contrastive learning~\citep{clip,elizalde2023clap}, and auto-regressive learning~\citep{team2024gemini,aiello2024jointly}. To circumvent the need for extensive training data and computational resources, 
numerous methods~\citep{huang2025towards,licsa,groger2025limited} also propose to leverage pretrained single-modal encoders, such as DINOv3~\citep{simeoni2025dinov3},  ModernBERT~\citep{warner2025smarter}, or BEATs~\citep{chen2022beats}, to construct semantically aligned representations across modalities. We focus on this setting of representation alignment.

One limitation of these approaches, however, is that the distribution of the joint representations remains uncontrolled. The geometry of the shared space is largely determined from the source encoders and shaped by the alignment objective, without being explicitly matched to a prescribed reference distribution. Consequently, the learned representations may exhibit undesirable geometric properties, such as anisotropy, distortions, or even mode collapse. This motivates our central research question: 
can we construct joint latent representations that preserve semantic structure from multiple modalities while conforming to a prescribed marginal distribution? 

We address this question from the perspective of Gromov--Wasserstein optimal transport (GWOT)~\citep{memoli2011gromov}. GWOT constructs a probabilistic coupling between distributions supported on distinct spaces by minimizing the discrepancy between their within-space pairwise relations. The resulting coupling therefore aligns the two distributions by preserving, as closely as possible, the relational structure of the source representations.

It is therefore natural to extend the standard GWOT formulation~\citep{peyre2016gromov} to multimodal representation learning. Such an extension should integrate two types of information. The first is the semantic structure of the observed features within each modality, which should be preserved when modality-specific knowledge is transferred to the joint representation space. The second is the set of cross-modal correspondences, which characterizes the relationships among observations from different modalities and guides their alignment in the joint space.

\begin{figure}[ht]
    \centering
    \includegraphics[width=0.7\linewidth]{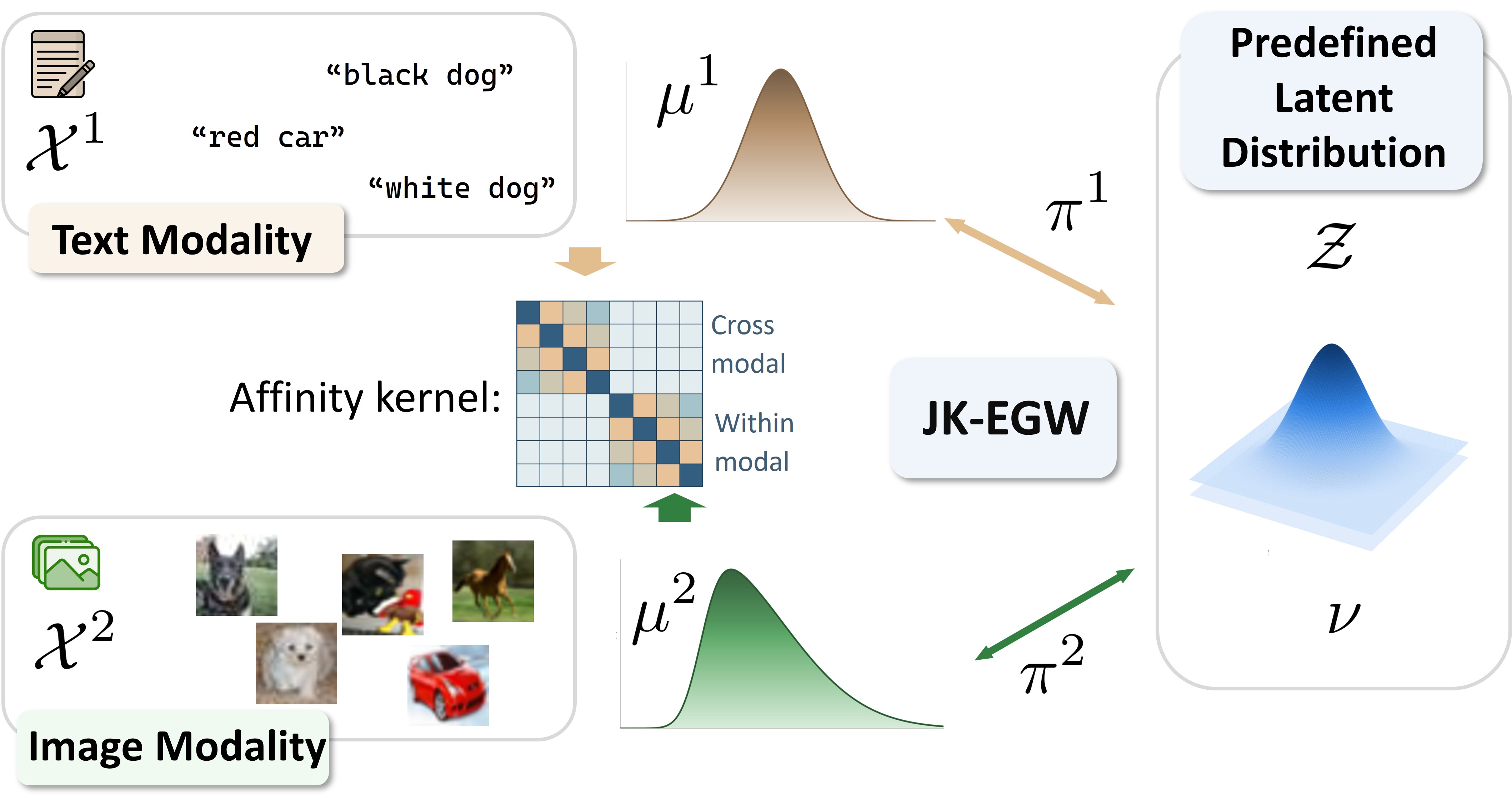}
    \caption{A global affinity kernel $K$ scores observation pairs within and across modalities. Each modality is coupled to a prescribed latent reference $\nu$. 
    }
    \label{fig:intro_illustration}
\end{figure}

To capture both sources of information, we construct a global affinity kernel $K$ on the disjoint union of all modality-specific representation spaces. Consider, for example, an image--text alignment problem. Within the image modality, an image of a black dog may have high affinity with another image of a black dog. Across modalities, the same image may also have high affinity with the text ``\texttt{black dog}'' (see Figure~\ref{fig:intro_illustration}). More generally, we propose to construct the within-modality component of $K$ by feature similarity metrics, and the cross-modality component of $K$ by domain knowledge, including class labels~\citep{cristianini2001kernel}, expert knowledge~\citep{kulis2005semi}, or similarities derived from pretrained representations~\citep{alvarez2018gromov}.

Building on this affinity kernel, we introduce the joint kernel entropic Gromov--Wasserstein (JK-EGW) objective, which formulates joint representation learning as an optimal transport problem. JK-EGW aligns the distribution of each modality with a common reference distribution in the latent space through  probabilistic couplings. The objective encourages observations with high affinity to link to nearby latent representations, thereby unifying within-modality structure preservation and cross-modal alignment.

In practical applications, instance-level cross-modal pairs are costly to collect and are therefore often scarce. It is thus important to characterize the finite-sample statistical error of JK-EGW. Establishing such guarantees is challenging, however, because the general GWOT objective takes the form of a nonconvex quadratic program~\citep{pmlr-v162-scetbon22b}, whose global optimization can be NP-hard~\citep{burkard1984quadratic,kravtsova2024np}. This nonconvexity complicates the statistical analysis and the optimization of empirical solutions.

We address this challenge by studying a specific class of kernel GWOT problems whose affinity kernel $K$ approximately admits a separable decomposition in terms of products of basis functions. For a finite-rank truncation of $K$, we show that the resulting quadratic GWOT problem can be lifted to a linear optimal transport (LOT) problem whose cost function depends on an additional set of auxiliary variables. This reformulation makes the statistical error amenable to analysis using tools developed for LOT. In particular, we establish finite-sample guarantees showing that the empirical optimal value converges to its population counterpart at the parametric rate $O(n^{-1/2})$. For the full-kernel objective, the bound includes an additional truncation term that quantifies the discrepancy between the original kernel and its finite-rank approximation. Our results therefore distinguish the statistical error arising from finite sampling from the approximation error introduced by kernel truncation.

Computationally, the basis-product structure naturally motivates an efficient alternating optimization algorithm for JK-EGW. We first construct a low-rank approximation of the empirical kernel matrix and use it to derive a variational lifting of the original quadratic objective. The lifted optimization problem should optimize over both the transport plan and the auxiliary variables. With fixed transport plan, the auxiliary variables admit a closed-form update. With fixed auxiliary variables, the transport-plan update reduces to a standard entropic optimal transport (EOT) problem over discrete distributions. The latter can therefore be solved using existing scalable OT solvers, such as Sinkhorn-type algorithms~\citep{cuturi2013sinkhorn}. Therefore, we can construct an alternating minimization algorithm by iteratively update the two parameter blocks. 

We evaluate JK-EGW on image--text alignment tasks derived from MS--COCO~\citep{lin2014microsoft}. The results show that JK-EGW learns geometrically coherent shared representations and achieves decent performance in bidirectional image--text retrieval.

In summary, our main contributions are as follows:
\begin{enumerate}
\item We formulate joint kernel entropic Gromov--Wasserstein optimal transport (JK-EGW) for aligning observations from multiple modalities into a shared representation space whose marginal distribution follows a prescribed reference distribution.

\item For kernels admitting a separable basis decomposition, we establish a finite-sample error bound of JK-EGW consisting of a statistical term of order $O(n^{-1/2})$ and a kernel-truncation term. Our analysis introduces a new variational formulation and develops techniques for establishing the regularity of optimal transport plans under costs induced by the kernel basis.

\item Building on the variational formulation of the empirical JK-EGW problem, we develop a scalable alternating optimization algorithm. One parameter block admits a closed-form update, while the other reduces to an entropic optimal transport (EOT) problem that can be solved using standard EOT solvers.

\item We empirically evaluate JK-EGW on digit multi-feature alignment, image--text retrieval tasks to demonstrate its applicability.

\end{enumerate}

\subsection{Related Works}


\paragraph{Alignment of pretrained unimodal encoders.}
Multiple works leverage frozen pretrained unimodal encoders and learn only a lightweight post-hoc alignment mechanism from scarce cross-modal correspondence data. Early works, such as APE~\citep{rosenfeld2022ape} and \textsc{ImageBind}~\citep{girdhar2023imagebind}, align several modalities onto a fixed, pretrained image–text embedding space by minimizing InfoNCE loss. Recent vision--language alignment methods further develop this idea by assessing or selecting compatible frozen encoder pairs and training simple projectors or alignment layers between their embedding spaces~\citep{maniparambil2025harnessing,zhang2025assessing}. More directly related to our empirical baselines, CSA~\citep{licsa} maps frozen unimodal features into a multimodal feature space and uses a weighted similarity score. STRUCTURE~\citep{groger2025limited} improves limited-data alignment by regularizing the learned representation to preserve the multiscale neighborhood geometry of pretrained unimodal latent spaces. Together, these methods show that post-hoc alignment of frozen encoders is a practical and data-efficient alternative to training multimodal models from scratch. Our work follows a similar rationale but studies a different formulation of this aligning problem, which allows for an explicitly chosen marginal distribution on the shared latent space.  Prior work suggests that latent geometry influences compositional and cross-modal generation~\citep{girdhar2023imagebind}, robustness to missing or noisy modalities~\citep{tsai2018learning}, downstream task performance~\citep{jiang2023understanding}, and the emergence and impact of modality gaps~\citep{liang2022mind}. 

\paragraph{Structure-aware alignment via optimal transport.}
Optimal transport (OT) provides a distributional framework for aligning heterogeneous datasets through learned couplings. Several works have used OT to improve multimodal or cross-domain alignment by exploiting structure beyond independent instance matching. For example, hierarchical OT leverages clustered structure to improve alignment in noisy or ambiguous settings~\citep{lee2019hierarchical}, while MuLOT uses self-attention for intra-modal correspondence and OT for cross-modal correspondence in resource-constrained multimodal sarcasm and humor detection~\citep{pramanick2022multimodal}. More recently, \citet{xu2025optimal} formulates speech--text alignment in spoken language models as a linear OT regularization between speech and transcript embeddings, using the resulting transport plan to construct a regularization loss during training. These works demonstrate the usefulness of OT for cross-domain and cross-modal correspondence learning. However, the linear OT formulations are often used to align features with the same dimensions. Gromov--Wasserstein-type OT (GWOT)~\citep{memoli2011gromov} can handle distributions with different dimensions as it compares distributions through pairwise relational structure rather than direct pointwise costs. GWOT has also been popular in aligning tasks, such as for words across two different languages~\citep{alvarez2018gromov}, and for single-cell genomics~\citep{klein2024genot}. Our work develops a joint kernel entropic Gromov--Wasserstein formulation that allows for the alignment across multiple distributions.

\paragraph{Statistical analysis for multimodal alignment and optimal transport.}
Recent work has begun to develop statistical theories for both contrastive multimodal representation learning and optimal-transport-based alignment. On the contrastive learning side, \citet{oko2025statistical} studies why contrastive pretraining can yield transferable multimodal representations and derives sample complexity guarantees for downstream multimodal learning. On the optimal transport side, parametric sample complexity bounds have been established for Sinkhorn divergences and EOT under suitable conditions. \citet{genevay2019sample} established the first sample complexity bounds for Sinkhorn divergences, showing that the empirical Sinkhorn divergence converges to its population counterpart at the parametric rate $n^{-1/2}$. \citet{mena2019statistical} strengthened these guarantees for EOT with the squared-Euclidean cost between sub-Gaussian measures in arbitrary dimension, obtaining the same parametric $n^{-1/2}$ rate but improving the $\varepsilon$-dependence of the constants exponentially over \citet{genevay2019sample} and extending the bound to unbounded measures. \citet{groppe2024lower} further showed that empirical EOT obeys a ``lower complexity adaptation'' principle, under which the convergence rate is governed by the lower-dimensional of the two measures. For Gromov--Wasserstein (GW) distances, \citet{zhang2024gromov} established duality and sample complexity results for a specific form of quadratic GW and entropic GW (EGW) in the Euclidean spaces, which allows for a parametric empirical convergence rate. Related work by \citet{rioux2024entropic} further develops stability and algorithmic guarantees for EGW. Our analysis is closest in spirit to the EGW theory, but aims at mapping several feature distributions to a common shared latent space under a kernelized distance, instead of in the Euclidean space. The resulting analysis requires new tools from finite-rank approximation. 

\subsection{Notation}
We denote the standard Euclidean norm and inner product by $\norm{\cdot}$ and $\langle\cdot,\cdot\rangle$, respectively. For matrices $A$ and $A'$ of the same size, $\left<A,A'\right>_F \coloneqq \operatorname{tr}{A^\top A'}$ denotes the Frobenius inner product and $\norm{A}_F \coloneqq \left<A,A\right>_F^{1/2}$ the Frobenius norm. The closed ball in $\mathbb{R}^d$ centered at the origin with radius $r$ is denoted by $B_d(r)$.

Let $\mathcal{X}^i$ represent the feature space for the $i$-th modality, and let $\mathcal{Z}$ be the shared latent space across different modalities. We denote their respective probability measures by $\mu^i \in \mathcal{P}(\mathcal{X}^i)$ and $\nu \in \mathcal{P}(\mathcal{Z})$, , where $\mathcal{P}(\cdot)$ denotes the set of Borel probability measures on the corresponding space. A probability measure $\rho \in \mathcal{P}(\mathbb{R}^d)$ is defined as $\beta$-sub-Weibull with parameter $\sigma^2$ ($\sigma \geq 0$) if $\int \exp(\norm{x}^{\beta}/2\sigma^2)d\rho(x)\leq 2$. Notably, the condition that $X\sim\rho$ is 4-sub-Weibull is equivalent to $\norm{X}^2$ being sub-Gaussian, satisfying the condition $\int \exp(t\norm{x}^2)d\rho(x)\leq 2\exp(t^2\sigma^2/2)$.

For any $p\in[1,\infty)$ and measure $\rho\in\mathcal{P}(\mathbb{R}^d)$, $L^p(\rho)$ refers to the space of measurable functions $f$ with finite norm $\norm{f}_{L^p(\rho)} \coloneqq(\int_{\mathbb{R}^d} |f|^p d\rho)^{1/p}$. 

Given a function class $\mathcal{F}$ equipped with metric $m$, we denote the $\epsilon$-covering number by $N(\mathcal{F},m,\epsilon)$ and the bracketing number by $N_{[\,]}(\mathcal{F},m,\epsilon)$. The relation $\lesssim_x$ indicates inequality up to a constant factor depending solely on $x$. Finally, for $a,b \in \BR$, we use $a \vee b$ to denote $\max\{a,b\}$.

\section{Joint Kernel Gromov--Wasserstein Optimal Transport}\label{sect:2}

\subsection{Problem formulation}

We consider a multimodal alignment problem with $M$ modalities, indexed by $[M]\coloneqq\{1,\ldots,M\}$. For each modality $i\in[M]$, observations take values in a feature space $\mathcal{X}^i$ and are distributed according to a probability measure $\mu^i\in\mathcal{P}(\mathcal{X}^i)$. The feature spaces may differ across modalities. Without loss of generality, we assume our dataset consists of $n$ observations from each modality, $\{x^1_k,\dots, x^M_k \}_{k=1}^n$. To use a unified notation, we define their disjoint union $\mathcal{X}\coloneqq\bigsqcup_{i=1}^M \mathcal{X}^i$ so that each $x\in\mathcal{X}$ belongs to a unique modality-specific component. When the modality is clear from context, we suppress the modality index $i$. 

We aim to align observations from different modalities into a shared latent space $\mathcal{Z}\subseteq\mathbb{R}^{d_z}$ that represents the shared information across modalities. Importantly, we equip $\mathcal{Z}$ with a fixed marginal distribution $\nu\in\mathcal{P}(\mathcal{Z})$ that can be defined by the user. $\nu$ specifies the distributional landscape of $\mathcal{Z}$. The alignment between a modality $i \in [M]$ and the latent space is represented by a probabilistic coupling $\pi^i\in\Gamma(\mu^i,\nu)$, where $\Gamma(\mu^i,\nu)$ denotes the set of probability measures on $\mathcal{X}^i\times\mathcal{Z}$ with first marginal $\mu^i$ and second marginal $\nu$. Naturally, the conditional distribution $\pi^i(\cdot\mid x)$ describes how $x$ is associated with locations in the latent space. 

Under the constraint of the marginal distributions of $\nu$ and $\mu^i$'s, the task of multi-modal alignment transforms to finding the probabilistic couplings $\pi^i$'s. Ideally, these couplings should preserve the relational structure of the observations. Since different modalities could have different dimensions, Gromov--Wasserstein optimal transport (GWOT) provides a natural framework for finding couplings that minimize an objective defined to induce alignment. 

Specifically, the standard quadratic GWOT aligns distributions by comparing pairwise relational structures in different spaces. One specifies relational functions $c_{\mathcal{X}}(x,x')$ and $c_{\mathcal{Z}}(z,z')$ on the source and target spaces, and seeks a coupling whose paired samples minimize the relational difference in aggregate, 
\begin{equation*}
    \inf_\pi \iint
    \bigl|c_{\mathcal{X}}(x,x') - c_{\mathcal{Z}}(z,z')\bigr|^2
    \,d\pi(x,z)\,d\pi(x',z').
\end{equation*}
Despite its simplicity, quadratic GWOT may be insufficient to align heterogeneous modalities as it calculates the distances of the raw features. When observations are collected from different spaces, using distances based on raw features may fail to reflect modality level nuances. This motivates us to replace the source-space geometric relation $c_{\mathcal{X}}$ and $c_{\mathcal{Z}}$ with a modality-aware and task-specific relation. 

To specify such relational structure, we assume access to a global affinity kernel
\begin{equation*}
K:\mathcal{X}\times\mathcal{X}\to\mathbb{R}_{\ge 0} \, .
\end{equation*}
This kernel is defined on the disjoint union $\mathcal{X}$, so that a single quantity $K(x,x')$ can describe both within-modality similarity when $x,x'\in\mathcal{X}^i$, and cross-modality similarity when $x\in\mathcal{X}^i$ and $x'\in\mathcal{X}^j$ with $i\neq j$. We interpret larger values of $K(x,x')$ as stronger task-relevant semantic similarity. The kernel may be specified using labels, expert knowledge, pretrained representations, or other domain-specific similarity scores. 

The affinity score encoded by $K$ determines how strongly pairs of observations should be encouraged to remain close after alignment. For coupled samples $(x,z)$ and $(x',z')$, we assign the pairwise latent penalty $K(x,x')\|z-z'\|^2$. The penalty is large when observations with high affinity are assigned distant latent representations, and is small when either the observations have weak affinity or their latent representations are close. Therefore, minimizing the aggregate penalty encourages observations that are related under $K$ to concentrate near one another in the shared latent space. 

Motivated by the kernel-weighted pairwise penalty shown above, we define our JK-EGW objective as
\begin{equation}\label{eq:main_objective}
\begin{aligned}
    S_{\epsilon}(\{\mu^i\}_{i=1}^M,\nu)
    \coloneqq
   \inf_{\substack{\pi^i \in \Gamma(\mu^i,\nu) \\ i \in [M]}} 
    \Bigg\{
    &\sum_{i,j=1}^M
    \int_{\mathcal{X}^i\times\mathcal{Z}}
    \int_{\mathcal{X}^j\times\mathcal{Z}}
    K(x,x')\|z-z'\|^2
    \,d\pi^i(x,z)\,d\pi^j(x',z')
    \\
    &\quad
    +\epsilon\sum_{i=1}^M
    D_{\mathrm{KL}}\!\left(
        \pi^i\,\middle\|\,\mu^i\otimes\nu
    \right)
    \Bigg\}.
\end{aligned}
\end{equation}
Here, $\epsilon>0$ is the entropic regularization parameter, and
$D_{\mathrm{KL}}$ denotes the Kullback--Leibler (KL) divergence 
\[
    D_{\mathrm{KL}}(\alpha\|\beta)
    =
    \int \log\!\left(\frac{d\alpha}{d\beta}\right)\,d\alpha \, ,
\]
where $\alpha$ and $\beta$ are probability measures. In the first term of~\Cref{eq:main_objective}, the variables $x$ and $x'$ are integrated over $\mathcal{X}^i$ and $\mathcal{X}^j$, respectively, through the couplings $\pi^i$ and $\pi^j$. The ordered sum over $(i,j)$ aggregates relations within and across modalities. The terms with $i=j$ encourage observations within the same modality to remain close in the latent space, while the terms with $i\neq j$ encourage cross-modality alignment. The KL term is the entropic regularization in~\Cref{eq:main_objective} that adds a strictly convex regularizing component to the GW-type objective. The problem of multi-modal alignment becomes choosing couplings $\{\pi^i\}_{i=1}^M$ so that pairs of observations with large $K(x,x')$ are assigned by nearby latent representations $z$ and $z'$.

\subsection{Finite-rank kernel approximation}
The objective in~\Cref{eq:main_objective} is quadratic in the unknown couplings because its first term averages a pairwise interaction over products of the form $d\pi^i(x,z)d\pi^j(x',z')$. The problem is a constrained nonconvex problem and potentially difficult to analyze. To make the problem amenable for analysis, we assume a special class of affinity kernel, which are decomposable kernels that can be approximated by finite-rank expansions~\citep{aronszajn1950theory,mercer1909xvi,weisse2006kernel}. The spectral representation of $K$ separates the dependence on $x$ and $x'$, so that the pairwise interaction can be expressed through feature-latent moments of the couplings. 
This subsection discusses the key assumptions and the induced error of kernel truncation, preparing for the variational lifting of~\Cref{eq:main_objective} in the next subsection. 

To state the spectral condition, let $\rho \coloneqq \frac{1}{M}\sum_{i=1}^M \mu^i$ denote the mixture distribution on the disjoint union $\mathcal{X}$. We use $\rho$ as the reference measure for the spectral analysis of $K$, since it aggregates the modality-specific marginals into a single probability measure on $\mathcal{X}$ where the kernel is defined. Throughout this paper, we assume that $K$ is symmetric positive semidefinite, and thus admits
a Mercer decomposition~\citep{mercer1909xvi} with respect to $\rho$,
\[
    K(x,x')
    =
    \sum_{r=1}^{\infty}
    \lambda_r \phi'_r(x)\phi'_r(x'),
\]
where $\lambda_r\ge 0$ and $\{\phi'_r\}_{r\ge 1}$ are orthonormal eigenfunctions
in $L^2(\rho)$ of the integral operator $T_Kf=\int_{\mathcal X}K(\cdot,x')\,f(x')\,d\rho(x')$. Such a decomposition exists under mild regularity conditions on $K$ and $\rho$. For instance, it follows from Mercer's theorem when $\mathcal X$ is compact and $K$ is continuous. A discussion on more general domains and kernels is given by~\citet{steinwart2012mercer}. We use $\phi_r(x) \coloneqq \sqrt{\lambda_r}\phi'_r(x)$ to denote the scaled eigenfunction. For a truncation rank $R\in\mathbb N$, define $\Phi_R(x) \coloneqq
\bigl(\phi_1(x),\ldots,\phi_R(x)\bigr)^\top\in\mathbb R^R$ as the vector comprised of the first $R$ scaled eigenfunction, 
and
$$
    K_R(x,x')
    \coloneqq
    \Phi_R(x)^\top\Phi_R(x')
    =
    \sum_{r=1}^R \phi_r(x)\phi_r(x')
$$
as the rank-$R$ approximation of the kernel $K$.

The following assumption collects the regularity conditions used in this paper to construct the finite-rank approximation and to control its truncation error.
\begin{assumption}[Kernel regularity and Mercer truncation]
\label{assumption2.1}

We impose the following conditions.
\begin{enumerate}
    \item \emph{Uniform diagonal bound:} there exists $\kappa>0$ such that $
        K(x,x)\le \kappa^2
        $, for all $x\in\mathcal X .
    $
    Consequently, $
    \|\Phi_R(x)\|^2=K_R(x,x)\le K(x,x)\le \kappa^2$, for all $R\ge 1,\ x\in\mathcal X .$

    \item \emph{Uniform Mercer-tail control:} we define the kernel approximate residual
    \[
        \mathcal E_R(x,x')
        \coloneqq
        K(x,x')-K_R(x,x')
        =
        \sum_{r>R}\phi_r(x)\phi_r(x'),
    \]
    and its upper bound
    \[
        \eta_R
        \coloneqq
        \sup_{x\in\mathcal X}\mathcal E_R(x,x)
        =
        \sup_{x\in\mathcal X}
        \sum_{r>R}\phi_r(x)^2 .
    \]
    We assume $
        \eta_R\to 0
        $ as $R\to\infty$.
\end{enumerate}
\end{assumption}

From the definition of $\mathcal E_R$, we can apply the Cauchy--Schwarz inequality
to obtain an upper bound on $|\mathcal E_R(x,x')|$,
\[
    |\mathcal E_R(x,x')|
    \le
    \sqrt{\mathcal E_R(x,x)\mathcal E_R(x',x')}
    \le
    \eta_R .
\]

The subsequent theoretical analysis in this paper depends heavily on the kernel tail properties via $
    \eta_R
    =
    \sup_{x\in\mathcal X}\sum_{r>R}\lambda_r\,(\phi'_r(x))^2$. Thus, we shall first verify the two conditions in Assumption~\ref{assumption2.1} for popular kernels. The next proposition shows that both conditions hold for \emph{any} continuous positive semidefinite kernel on a compact domain. It also estimates the decay rates for the radial basis function (RBF) kernel. 

\begin{proposition}[Satisfiability and Mercer-tail rates]
\label{prop:kernel}
\leavevmode
    Let $\mathcal X$ be compact, let $\rho$ be a finite Borel measure on $\mathcal X$ with $\text{supp}(\rho)=\mathcal X$, and let $K:\mathcal X\times\mathcal X\to\mathbb R$ be a continuous symmetric positive semidefinite kernel. Then $K$ admits the Mercer decomposition of Assumption~\ref{assumption2.1}, where the uniform diagonal bound holds with $       \kappa^2=\max_{x\in\mathcal X}K(x,x)<\infty$,
    and $\eta_R\to0$ as $R\to\infty$. 
    
    In particular, the radial basis function (RBF) kernel
    \[
        K(x,x')=\exp\!\Bigl(-\frac{\|x-x'\|^2}{2\ell^2}\Bigr), \qquad \ell>0,
    \]
    satisfies Assumption~\ref{assumption2.1} on every compact domain. For the RBF kernel, $K(x,x)=1$ for every $x$, so the diagonal bound holds with $\kappa^2=1$. More specifically, suppose $\mathcal X=[0,1]^{d_x}$ with $\rho$ the uniform measure, and let $K$ be stationary, i.e. $K(x,x')=k(x-x')$ for some function $k:\mathcal X\to\mathbb R$. 
Then the following hold.
    
    \begin{enumerate}
        \item[(a)] There exists a Mercer eigensystem satisfying Assumption~\ref{assumption2.1}, under which, for every $R\ge1$,    $\eta_R\le 2\sum_{r>R}\lambda_r$, where $\{\lambda_r\}_{r\ge1}$ is the non-increasing rearrangement of $\{\lambda_m\}_{m\in\mathbb Z^{d_x}}$.
        \item[(b)] The eigenvalues of the periodized RBF kernel $
        K(x,x')=\sum_{w\in\mathbb Z^{d_x}}\exp\!\Bigl(-\frac{\|x-x'+w\|^2}{2\ell^2}\Bigr)$
satisfy
\[
    \lambda_m\lesssim\exp(-2\pi^2\ell^2\|m\|^2),
    \qquad m\in\mathbb Z^{d_x},
\]
and consequently $\eta_R\lesssim e^{-c'R^{2/d_x}}$ for some constant $c'>0$ depending only on $\ell$ and $d_x$.
    \end{enumerate}
\end{proposition}

Proposition~\ref{prop:kernel} shows that the first condition of Assumption~\ref{assumption2.1} holds for every continuous positive semidefinite kernel on a compact domain with full support. As a special case, Proposition~\ref{prop:kernel} also provides an upper bound on the $\eta_R$ for periodized RBF kernels, which vanishes exponentially when $R\rightarrow\infty$. The proof of Proposition~\ref{prop:kernel} is deferred to Section~\ref{proof:prop_kernel}. We consider kernels that satisfy Assumption~\ref{assumption2.1} hereon.


Fix a truncation rank $R\in\mathbb{N}$. Under Assumption~\ref{assumption2.1}, the finite-rank term $\Phi_R(x)^\top\Phi_R(x')$ separates the dependence of the kernel on $x$ and $x'$. We can thus introduce a rank-$R$ JK-EGW objective by replacing $K$ by its rank-$R$ approximation,
\begin{equation}\label{eq:main_kernel_objective}
\begin{aligned}
    S^{R}_{\epsilon}(\{\mu^i\}_{i=1}^M,\nu)
    \coloneqq
    \inf_{\substack{\pi^i \in \Gamma(\mu^i,\nu) \\ i \in [M]}}
    \Bigg\{
    &\sum_{i,j=1}^M
    \int_{\mathcal{X}^i\times\mathcal{Z}}
    \int_{\mathcal{X}^j\times\mathcal{Z}}
    \Phi_R(x)^\top \Phi_R(x') \|z-z'\|^2
    \,d\pi^i(x,z)\,d\pi^j(x',z')
    \\
    &\quad
    +\epsilon \sum_{i=1}^M
    D_{\mathrm{KL}}\!\left(
        \pi^i\,\middle\|\,\mu^i\otimes\nu
    \right)
    \Bigg\}.
\end{aligned}
\end{equation}

Replacing $K$ by its rank-$R$ truncation should not introduce a large error, and the error it does introduce is explicitly controllable through the choice of $R$. Since the problems in~\Cref{eq:main_objective} and~\Cref{eq:main_kernel_objective} share the same feasible set and the same entropic regularization term, their values differ only through the quadratic term, where the kernel is perturbed by the residual $\mathcal{E}_R$, which Assumption~\ref{assumption2.1} bounds uniformly by $\eta_R$. The truncation therefore introduces an error controlled by the Mercer tail of the kernel, and the following proposition quantifies this error.
\begin{proposition}[Kernel truncation error]
\label{prop:kernel_truncation_error}
Suppose Assumption~\ref{assumption2.1} holds and $\nu$ is $4$-sub-Weibull with parameter $\sigma^2$. Then, for every $\epsilon>0$ and every $R\ge 1$, 
\[
\begin{aligned}
    \left|S_\epsilon(\{\mu^i\}_{i=1}^M,\nu)-S_\epsilon^R(\{\mu^i\}_{i=1}^M,\nu)\right|\lesssim M^2\sigma\,\eta_R.
\end{aligned}
\]
\end{proposition}

Proposition~\ref{prop:kernel_truncation_error} suggests that the error induced by kernel truncation error bound is dependent only on the number of modalities $M$, the sub-Weibull parameter $\sigma$ of the marginal distribution of the latent, and the Mercer-tail quantity $\eta_R$ defined in Assumption~\ref{assumption2.1}. Hence, the error vanishes with $\eta_R\rightarrow0$. The proof of Proposition~\ref{prop:kernel_truncation_error} is deferred to Section~\ref{proof:kernel_error}. 

\subsection{Variational lifting to entropic optimal transport}
\label{sec:population_lifting}
The finite-rank objective in~\Cref{eq:main_kernel_objective} presents new structures that are not apparent in the original formulation~\Cref{eq:main_objective}. One key observation is that, after replacing $K(x,x')$ by $\Phi_R(x)^\top\Phi_R(x')$, the quadratic interaction among couplings depends on each coupling only through $R$-dimensional feature eigenfunctions $\Phi$. This allows the quadratic problem to be lifted into a linear optimization problem, which an EOT solver could tackle. The resulting variational representation is the main bridge between the rank-$R$ JK-EGW objective and standard EOT tools.

To introduce the lifting, we define $
    m_\Phi
    \coloneqq
    \sum_{j=1}^M
    \int_{\mathcal{X}^j} \Phi_R(x')\,d\mu^j(x')
    \in \mathbb{R}^R$ as the averaged eigenfunction. 
For an arbitrary matrix $A\in\mathbb{R}^{R\times d_z}$, we introduce the modality-specific cost
\[
    c_{A,m_\Phi}^i(x,z)
    \coloneqq
    2\,\Phi_R(x)^\top m_\Phi\,\|z\|^2
    -
    4\left\langle A,\Phi_R(x)z^\top\right\rangle_F ,
    \qquad
    x\in\mathcal{X}^i,\ z\in\mathcal{Z}.
\]
We further define a linear EOT problem with the cost defined by $c_{A,m_\Phi}^i(x,z)$,
\[
    \OT^i_{A,m_\Phi,\epsilon}(\mu^i,\nu)
    \coloneqq
    \inf_{\pi^i\in\Gamma(\mu^i,\nu)}
    \left\{
        \int c_{A,m_\Phi}^i(x,z)\,d\pi^i(x,z)
        +
        \epsilon
        D_{\mathrm{KL}}\!\left(
            \pi^i\,\middle\|\,\mu^i\otimes\nu
        \right)
    \right\}.
\]
The following theorem presents the transformation from the quadratic regularized optimization problem to a linear regularized optimization problem with the additional optimization variable $A$.

\begin{theorem}[Variational lifting of the rank-$R$ EGW objective]
\label{thrm:pop_duality}
Fix $\epsilon>0$. Suppose Assumption~\ref{assumption2.1} holds and
$\nu$ is a 4-sub-Weibull distribution with parameter $\sigma^2>0$. Then
\begin{equation}\label{eq:threm2.1}
    S^R_{\epsilon}(\{\mu^i\}_{i=1}^M,\nu)
    =
    \inf_{A\in\mathbb{R}^{R\times d_z}}
    \left\{
        2\|A\|_F^2
        +
        \sum_{i=1}^M
        \OT^i_{A,m_\Phi,\epsilon}(\mu^i,\nu)
    \right\}.
\end{equation}
Moreover, the infimum is achieved at $A^*\in\mathcal{A}_M \coloneqq\left\{A\in\mathbb{R}^{R\times d_z}:
        \|A\|_F\le M\kappa(2\sigma^2)^{1/4}\right\}$. 
\end{theorem}

The proof of Theorem~\ref{thrm:pop_duality} is deferred to
Section~\ref{proof:pop_duality}. The lifted formulation introduces an auxiliary matrix $A$ to capture the aggregate feature-latent moment, while each inner term $\OT^i_{A,m_\Phi,\epsilon}(\mu^i,\nu)$ is a standard EOT problem with cost $c_{A,m_\Phi}^i$. Since only $A$ has the quadratic terms $\|A\|_F^2$, the lifted formulation separates the nonconvex quadratic interaction from the modality-wise transport problems. This lifted formulation, enabled by the decomposition of the affinity kernel, is the foundation for the statistical analysis in Section~\ref{sec:theory}. In literature, the similar duality of quadratic EGW is analyzed in~\citet{zhang2024gromov}, whose auxiliary matrix linearizes the cross-covariance of raw coordinates. In contrast, our lifting linearizes the moment $\int \Phi_R(x)z^\top\,d\pi^i$ of the kernel feature map against the latent coordinates. The lifted formulation is thus set by the  rank $R$ rather than by the raw feature dimension. 

\section{Theoretical Analysis}\label{sec:theory}
In almost any application of multi-modality alignment, the raw feature distribution $\mu^i$'s are rarely observed directly. We have access only to a finite set of paired samples, from which we form an empirical version of the JK-EGW objective. In the data-scarce regime, we therefore are interested in how fast the empirical objective converges to its population counterpart as the number of samples $n$ grows. This is the central question of the proposed alignment that will be studied in this section. 

To obtain the sample complexity guarantee of the objective $S_\epsilon$, we decompose the analysis in two stages. We first establish sample complexity bounds for the rank-truncated objective $S_\epsilon^R$, whose finite-dimensional representation allows the problem to be reduced to a family of EOT problems. We then combine this statistical bound with kernel truncation estimates to obtain a corresponding sample complexity guarantee for the full objective $S_\epsilon$. The resulting parametric bounds quantify how the empirical objective approaches its population counterpart as the number of samples increases, with constants governed by the latent dimension, the truncation level, the entropic regularization parameter, and the spectral decay of the kernel.

\subsection{Sample Complexity}
Let $x^1,\ldots,x^M \overset{\mathrm{i.i.d.}}{\sim} \rho$ gives a dataset of $\{x_k^1,\dots,x_k^M \}_{k=1}^n$, and let $z_1,\ldots,z_n \overset{\mathrm{i.i.d.}}{\sim} \nu$, independently of the modality samples. We define the empirical measures as
\[
    \hat{\mu}_n^i = \frac{1}{n}\sum_{k=1}^n \delta_{x_k^i},
    \qquad
    \hat{\nu}_n = \frac{1}{n}\sum_{\ell=1}^n \delta_{z_\ell}.
\]
Throughout this subsection, we work under the Assumption~\ref{assumption2.1}. Recall that we assume the latent reference distribution $\nu$ is $4$-sub-Weibull with parameter $\sigma^2$, while no tail condition is imposed on the modality distributions $\mu^i$. The following theorem studies both the two-sample deviation
\[
    S_\epsilon(\{\hat{\mu}_n^i\},\hat{\nu}_n)
    -
    S_\epsilon(\{\mu^i\},\nu),
\]
where both the modality distributions and the latent reference distribution are sampled, and the one-sample deviation
\[
    S_\epsilon(\{\hat{\mu}_n^i\},\nu)
    -
    S_\epsilon(\{\mu^i\},\nu),
\]
which is relevant when the latent reference distribution is specified or controlled by design.

\begin{theorem}[Full-kernel EGW sample complexity]
\label{cor:full}
Fixed $\epsilon>0$. Suppose Assumption~\ref{assumption2.1} holds, $\nu$ is $4$-sub-Weibull with parameter $\sigma^2$, and $\hat{\mu}_n^i$ and $\hat{\nu}_n$ are the empirical measures from $n$ samples, then the two-sample deviation is upper bounded by
\begin{align*}
&\mathbb E\left[
\left|
S_\epsilon(\{\mu^i\}_{i=1}^M,\nu)
-
S_\epsilon(\{\hat{\mu}_n^i\}_{i=1}^M,\hat{\nu}_n)
\right|
\right] \\
&\qquad\lesssim_{M,R,d_z,\kappa}
\epsilon
\left[
    1+\left(\frac{\sigma}{\epsilon}\right)+\left(\frac{\sigma}{\epsilon}\right)^{2q_{R,z}+2}
    +\left(\frac{\sigma}{\epsilon}\right)^{4q_{R,z}+4}
\right]
\frac{1}{\sqrt n}
+
\sigma\eta_R,
\end{align*}
and the one-sample deviation is upper bounded by
\begin{align*}
&\mathbb E\left[
\left|
S_\epsilon(\{\mu^i\}_{i=1}^M,\nu)
-
S_\epsilon(\{\hat{\mu}_n^i\}_{i=1}^M,\nu)
\right|
\right] \\
&\qquad\lesssim_{M,R,d_z,\kappa}
\epsilon
\left[
    1+\left(\frac{\sigma}{\epsilon}\right)+\left(\frac{\sigma}{\epsilon}\right)^{4q_{R,z}+4}
\right]
\frac{1}{\sqrt n}+\sigma\eta_R,
\end{align*}
where $q_{R,z}=\left\lceil \frac{R\vee d_z}{2}\right\rceil$, and the expectation is taken over the samples.
\end{theorem}

Theorem~\ref{cor:full} shows that, for fixed $\epsilon$, the empirical full-kernel objective converges to its population counterpart up to a kernel truncation error of order $\eta_R$ and a statistical error that decays at the parametric rate $n^{-1/2}$, with constants depending on $M$, $R$, $d_z$, and $\kappa$. The kernel truncation error is proportional to $\eta_R$, which will diminish with large rank $R$ since $\eta_R\rightarrow0$ as $R\rightarrow\infty$. The full error bound separates the statistical error induced by empirical sampling from the approximation error induced by the rank-$R$ kernel truncation. 

Interestingly, the error bound has no explicit dependence on the ambient input dimension of the raw feature spaces $d_x$. The statistical error term depends on the latent dimension $d_z$ and the truncation level $R$, while the kernel approximation error is governed by the spectral decay of $K$. Consequently, when the kernel has a rapidly decaying spectrum, the effective statistical complexity is controlled by the latent dimension and the effective kernel rank, instead of the raw dimension of the input spaces. This property makes~\Cref{eq:main_objective} amenable to align multi-modal data whose raw features are high-dimensional but structured.

The $n^{-1/2}$ rate is consistent with known convergence rates for EGW distances with Euclidean costs~\citep{zhang2024gromov}, as well as related sample complexity results for EOT~\citep{mena2019statistical}. The main difference is that our analysis applies to a kernelized feature-space cost, which is natural in heterogeneous multimodal settings where raw Euclidean distances across modalities are not directly comparable. 

\subsection{Proof Sketch}
To prove Theorem~\ref{cor:full}, which is a sample complexity guarantee for the full-kernel objective $S_\epsilon$, we will need to decompose the bound into three parts. There are two terms accounting for the truncation error of replacing kernel $K$ with its rank-$R$ counterpart for the population distributions and empirical distributions, respectively. There is another term accounting for the statistical error of the rank-truncated objective. More specifically, we have the following triangle inequality that decompose the error into three terms:
\begin{align}
&\left|
S_\epsilon(\{\mu^i\}_{i=1}^M,\nu)
-
S_\epsilon(\{\hat{\mu}_n^i\}_{i=1}^M,\hat{\nu}_n)
\right| \notag \\
&\le
\left|
S_\epsilon(\{\mu^i\}_{i=1}^M,\nu)
-
S_\epsilon^R(\{\mu^i\}_{i=1}^M,\nu)
\right| \tag{Truncation error on population distributions}\\
&
+
\left|
S_\epsilon^R(\{\mu^i\}_{i=1}^M,\nu)
-
S_\epsilon^R(\{\hat{\mu}_n^i\}_{i=1}^M,\hat{\nu}_n)
\right| \tag{Statistical error on rank-$R$ JK-EGW}\\
&
+
\left|
S_\epsilon^R(\{\hat{\mu}_n^i\}_{i=1}^M,\hat{\nu}_n)
-
S_\epsilon(\{\hat{\mu}_n^i\}_{i=1}^M,\hat{\nu}_n)
\right| \tag{Truncation error on empirical distributions}.
\end{align}
The first term on the RHS is controlled by the population truncation bound in Proposition~\ref{prop:kernel_truncation_error}. It remains to control the second term, which is the rank-$R$ JK-EGW statistical error bound, and the third term, which is a kernel truncation error evaluated at the empirical measures. The following Proposition analyzes the second term.

\begin{proposition}[Rank-truncated JK-EGW sample complexity]
\label{thrm:sample_complexity}
Suppose that Assumption~\ref{assumption2.1} holds. Fix $\epsilon>0$. Assume that the latent reference distribution $\nu$ is $4$-sub-Weibull with parameter $\sigma^2$. Let $\hat{\mu}_n^i$ and $\hat{\nu}_n$ be the empirical measures, with samples drawn independently across modalities and from the latent reference distribution. Set $    q_{R,z}=\left\lceil \frac{R\vee d_z}{2}\right\rceil$. 
Then,
\begin{equation*}
\begin{aligned}
    &\mathbb E\left[
    \left|
    S_\epsilon^R(\{\mu^i\}_{i=1}^M,\nu)
    -
    S_\epsilon^R(\{\hat{\mu}_n^i\}_{i=1}^M,\hat{\nu}_n)
    \right|
    \right] \lesssim_{M,R,d_z,\kappa}
    \epsilon
    \left[
        1+\left(\frac{\sigma}{\epsilon}\right)+\left(\frac{\sigma}{\epsilon}\right)^{2q_{R,z}+2}+
        \left(\frac{\sigma}{\epsilon}\right)^{4q_{R,z}+4}
    \right]
    \frac{1}{\sqrt n},
    \\
    &\mathbb E\left[
    \left|
    S_\epsilon^R(\{\mu^i\}_{i=1}^M,\nu)
    -
    S_\epsilon^R(\{\hat{\mu}_n^i\}_{i=1}^M,\nu)
    \right|
    \right] \lesssim_{M,R,d_z,\kappa}
    \epsilon
    \left[
        1+\left(\frac{\sigma}{\epsilon}\right)+\left(\frac{\sigma}{\epsilon}\right)^{4q_{R,z}+4}
    \right]
    \frac{1}{\sqrt n},
\end{aligned}
\end{equation*}
\end{proposition}

The proof of Proposition~\ref{thrm:sample_complexity} is deferred to Section~\ref{proof:thrm:sample_complexity}. We highlight a few technical challenges. First, since the transport cost of the lifted EOT problem is defined as $c_{A,m_\Phi}^i=2\Phi_R(x)^\top m_\Phi\norm{z}^2-4\left\langle A,\Phi_R(x)z^\top\right\rangle_F$, this term is dependent on the marginal distribution $\mu$ through $m_\Phi$.  The empirical version of $m_\Phi$ is defined as $\hat m_\Phi=\sum_{j=1}^M\int_{\mathcal X^j}\Phi_R(x')\,d\hat\mu_n^j(x')=\frac1n\sum_{j=1}^M\sum_{k=1}^n\Phi_R(x_k^j)$, which results in a different EOT transport cost. Such difference is not present in the previous statistical analysis of GWOT with Euclidean distance cost \citep{zhang2024gromov}.  To tackle the challenge, we develop an error bound on $\mathbb E\norm{m_\Phi-\hat m_\Phi}$ through the regularity conditions on $\Phi_R$. 


Second, our analysis works with kernel features without adding strong assumptions on the raw feature spaces $\mathcal X^i$. Therefore, new regularity conditions of the transport plan should be established. However, Assumption~\ref{assumption2.1} only assumes a uniform bound $\norm{\Phi_R(x)}\le\kappa$ on the feature map, which does not directly imply the smoothness of the Mercer eigenfunctions and or the conditions on the $\mu^i$. To establish the regularity, we leverage the fact that the cost $c^i_{A,m_\Phi}$ depends on $x$ only through $\Phi_R(x)$, which enables rebuilding the regularity theory behind the empirical-process argument in the feature space. Full characterizations of the regularity are in Lemma~\ref{lemma:1}, which is a core step in supporting the proof of Proposition~\ref{thrm:sample_complexity}. 
As a result of the new analysis, the regularity does not require tail or moment condition on the $\mu^i$. This is the mechanism behind the effective dimension $R\vee d_z$ in Theorem~\ref{cor:full} and Proposition~\ref{thrm:sample_complexity}. In contrast, the error analysis in~\citet{zhang2024gromov} and~\citet{mena2019statistical} involve constants depend on the ambient dimensions of the sample spaces, whose bounds require tail conditions on both marginals.

With the statistical term controlled by Proposition~\ref{thrm:sample_complexity}, the remaining ingredient toward Theorem~\ref{cor:full} is the kernel truncation error evaluated at the empirical measures. The following corollary records the empirical analogue of the population truncation estimate in Proposition~\ref{prop:kernel_truncation_error}.

\begin{corollary}[Empirical kernel truncation error]
\label{cor:empirical_truncation}
Suppose that Assumption~\ref{assumption2.1} holds, and assume that $\nu$ is $4$-sub-Weibull with parameter $\sigma^2$. Let $\hat{\mu}_n^i$ and $\hat{\nu}_n$ be the empirical measures from $n$ samples defined above. Then, for every $\epsilon>0$ and $R\ge 1$, 
\[
\begin{aligned}
    &\mathbb E \left[
    \left|
    S_\epsilon(\{\hat{\mu}_n^i\}_{i=1}^M,\hat{\nu}_n)
    -
    S_\epsilon^R(\{\hat{\mu}_n^i\}_{i=1}^M,\hat{\nu}_n)
    \right|
    \right]
    \lesssim
    M^2\sigma\,\eta_R,
\end{aligned}
\]
where the expectation is taken over the samples. Consequently, since $\eta_R\to0$ by Assumption~\ref{assumption2.1}, the
empirical truncation error vanishes in expectation as $R\to\infty$.
\end{corollary}

Indeed, applying Proposition~\ref{prop:kernel_truncation_error} conditionally on the samples gives
\[
\begin{aligned}
    &\left|
    S_\epsilon(\{\hat{\mu}_n^i\}_{i=1}^M,\hat{\nu}_n)
    -
    S_\epsilon^R(\{\hat{\mu}_n^i\}_{i=1}^M,\hat{\nu}_n)
    \right| \le
    M^2\eta_R
    \iint_{\mathcal Z\times\mathcal Z}
    \|z-z'\|^2\,d\hat{\nu}_n(z)d\hat{\nu}_n(z').
\end{aligned}
\]
With $\hat{\nu}_n=\frac{1}{n}\sum_{\ell=1}^n\delta_{z_\ell}$, we have
\[
    \mathbb E
    \iint
    \|z-z'\|^2\,d\hat{\nu}_n(z)d\hat{\nu}_n(z')
    =
    \left(1-\frac1n\right)\mathbb E\|Z-Z'\|^2
    \le
    2\mathbb E\|Z\|^2,
\]
where $Z,Z'\sim\nu$ are independent. Since $\nu$ is $4$-sub-Weibull with parameter $\sigma^2$,
\[
    \mathbb E\|Z\|^2
    \le
    \left(\mathbb E\|Z\|^4\right)^{1/2}
    \lesssim
    \sigma .
\]
Therefore,
\[
    \mathbb E \left[
    \left|
    S_\epsilon(\{\hat{\mu}_n^i\}_{i=1}^M,\hat{\nu}_n)
    -
    S_\epsilon^R(\{\hat{\mu}_n^i\}_{i=1}^M,\hat{\nu}_n)
    \right|
    \right]
    \lesssim
    M^2\sigma\,\eta_R .
\]

Recall that the kernel truncation error of the population distributions is bounded by Proposition~\ref{prop:kernel_truncation_error}. The rank-truncated statistical error is bounded by Proposition~\ref{thrm:sample_complexity}, and the kernel truncation error of the empirical distributions is bounded by Corollary~\ref{cor:empirical_truncation}. The decay rate of $\eta_R$ can be determined by substituting specific kernels, as in Proposition~\ref{prop:kernel}. The aforementioned propositions and corollary together build up to Theorem~\ref{cor:full}.




\section{Optimization for Discrete JK-EGW}\label{sec:optimization}

Having established the closeness between $S_\epsilon(\{\mu^i\}_{i=1}^M,\nu)$ and $S_\epsilon(\{\hat{\mu}_n^i\}_{i=1}^M,\hat{\nu}_n)$, we now focus on $S_\epsilon(\{\hat{\mu}_n^i\}_{i=1}^M,\hat{\nu}_n)$ and develop an algorithm to solve it. Recall that for each modality $i\in[M]$ we have $n$ i.i.d.\ samples $\{x_k^i\}_{k=1}^n \subset \mathcal X^i \subseteq \mathbb R^{d_i}$, yielding the empirical measure $\hat{\mu}_n^i = \frac{1}{n}\sum_{k=1}^n \delta_{x_k^i}$. We further generate $n$ latent support points $\{z_\ell\}_{\ell=1}^n \subset \mathcal Z$ by randomly identically sampling from the latent reference measure $\nu$. The shared empirical latent measure is constructed by $\hat{\nu}_n = \frac{1}{n}\sum_{\ell=1}^n \delta_{z_\ell}$. Our goal is to jointly compute the modality-wise empirical transport plans $\{\hat{\pi}^i\}_{i=1}^M$, where each $\hat{\pi}^i \in \Gamma(\hat{\mu}_n^i,\hat{\nu}_n)$ couples $\hat{\mu}_n^i$ to the shared $\hat{\nu}_n$ under uniform marginals. These empirical transport plans are parametrized by $n\times n$ (scaled) doubly stochastic matrices,
\begin{align*}
    \hat{\pi}^i \in \mathbb R^{n\times n}_+,
    \qquad \text{subject to}\qquad
    \hat{\pi}^i \vecone = \frac{1}{n}\vecone,
    \qquad
    \hat{\pi}^{i\top}\vecone = \frac{1}{n}\vecone .
\end{align*}
We define $\hat{\pi}_{k\ell}^i$ as the mass transported from sample $x_k^i$ in modality $i$ to the latent support point $z_\ell$. The discrete version of the exact-kernel JK-EGW objective is thus:
\begin{equation} \label{eq:sum_obj}
\begin{aligned}
\hat{S}_{\epsilon}\big(\{\hat{\mu}_n^i\}_{i=1}^M,\hat{\nu}_n\big)
=\min_{\{\hat{\pi}^i\}_{i=1}^M}\;
& \sum_{i,j} \sum_{k,k',\ell,\ell'}
W^{i,j}_{k,k'}\,\hat{\pi}^i_{k,\ell}\,\hat{\pi}^j_{k',\ell'}\,\norm{z_{\ell}-z_{\ell'}}^2
\;+\;\epsilon\sum_{i=1}^M D_{\mathrm{KL}}\!\big(\hat{\pi}^i \,\big\|\, \hat{\mu}_n^i\otimes\hat{\nu}_n\big)\\
\text{subject to}\quad
& \hat{\pi}^i\vecone = \tfrac{1}{n}\vecone,
\qquad
\hat{\pi}^{i\top}\vecone = \tfrac{1}{n}\vecone,
\qquad i\in[M].
\end{aligned}
\end{equation}

In~\Cref{eq:sum_obj}, $W^{i,j}\in\BR^{n\times n}$ is the empirical similarity matrix with its $k,k'$-th entry defined by $W^{i,j}_{k,k'}=K(x^i_k,x^j_{k'})$. We further concatenate  all $W^{i,j}$ matrices into $W\in\mathbb{R}^{N\times N}$, where $N=Mn$. More specifically, 
\[
    W
    =
    \begin{bmatrix}
        W^{1,1} & W^{1,2} & \cdots & W^{1,M}\\
        W^{2,1} & W^{2,2} & \cdots & W^{2,M}\\
        \vdots  & \vdots  & \ddots & \vdots \\
        W^{M,1} & W^{M,2} & \cdots & W^{M,M}
    \end{bmatrix}
    \in\BR^{Mn\times Mn},
    \qquad
    W^{i,j}_{k,k'}=K(x^i_k,x^j_{k'}).
\]

At the empirical level, the interaction term in~\Cref{eq:sum_obj} is also quadratic in the transport plans and therefore expensive to optimize directly. We similarly employ a low-rank approximation of the block kernel matrix $W$, followed by similar variational lifting procedures discussed in Section~\ref{sec:population_lifting}.

\subsection{Finite-rank Decomposition and Lifting}
We first approximate the empirical kernel matrix through a Nystr\"om-style feature factorization
\[
    W \approx LL^\top,
    \qquad
    L\in\BR^{N\times \hat R},
    \qquad
    \hat R\ll N,
\]
that can be computed by Pivoted Cholesky using kernel evaluations of the entries of $W$~\citep{williams2000using,harbrecht2012low}. Under the positive-semidefinite kernel assumption, the Pivoted Cholesky factorization constructs $L$ greedily from selected $\hat R$ pivot columns, terminating once the residual trace $\operatorname{tr}(W-LL^\top)$ falls below a prescribed tolerance $\tau_{\hat R}>0$, and yields a positive-semidefinite residual $W-LL^\top\succeq 0$ at each truncation level. We define $\rho_{\hat R}\coloneqq\max_{a,b\in[N]}\left|W_{ab}-(LL^\top)_{ab}\right|$ as the $\ell_\infty$ norm of the residual. Positive-semidefiniteness of the residual gives $\left|W_{ab}-(LL^\top)_{ab}\right|\le\sqrt{(W-LL^\top)_{aa}(W-LL^\top)_{bb}}\le\operatorname{tr}(W-LL^\top)$ for all $a,b\in[N]$, so the returned factor $L$ satisfies $\rho_{\hat R}\le\tau_{\hat R}$. Because each pivot step involves only the residual diagonal and one additional column of $W$, computing $L$ requires $O(N\hat R)$ kernel evaluations, $O(N\hat R^{2})$ arithmetic operations, and $O(N\hat R)$ memory~\citep{harbrecht2012low}. This linear scaling in $N$ makes Pivoted Cholesky well suited to large $N$. The calculation of $L$ is performed only once before our optimization algorithm.

We therefore denote by $\hat S_{\epsilon}^{\hat R}$ the empirical objective obtained from \Cref{eq:sum_obj} after replacing the global block matrix $W$ by its rank-$\hat R$ approximation $LL^\top$. 
\begin{equation} \label{eq:emp_sum_obj}
\begin{aligned}
\hat{S}_{\epsilon}^{\hat R}\big(\{\hat{\mu}_n^i\}_{i=1}^M,\hat{\nu}_n\big)
=\min_{\{\hat{\pi}^i\}_{i=1}^M}\;
& \sum_{i,j} \sum_{k,k',\ell,\ell'}
(LL^\top)_{k,k'}^{i,j}\,\hat{\pi}^i_{k,\ell}\,\hat{\pi}^j_{k',\ell'}\,\norm{z_{\ell}-z_{\ell'}}^2
\;+\;\epsilon\sum_{i=1}^M D_{\mathrm{KL}}\!\big(\hat{\pi}^i \,\big\|\, \hat{\mu}_n^i\otimes\hat{\nu}_n\big)\\
\text{subject to}\quad
& \hat{\pi}^i\vecone = \tfrac{1}{n}\vecone,
\qquad
\hat{\pi}^{i\top}\vecone = \tfrac{1}{n}\vecone,
\qquad i\in[M].
\end{aligned}
\end{equation}

Notice that the computational rank $\hat R$ is conceptually distinct from the analytic Mercer truncation rank $R$ in Section~\ref{sec:theory}, where the former is used for approximating the empirical kernel in practice while the later is used in the theoretical analysis. The additional error caused by the Nystr\"om-style kernel approximation is analyzed separately in later Proposition~\ref{prop:pivoted_cholesky_empirical_error}. 

We are now ready to derive the finite-dimensional lifting that would catalyze an optimization algorithm. For transport plans $\{\hat{\pi}^i\}_{i=1}^M$, define the vertically stacked
transport matrix 
\[
    \hat{\Pi}
    =
    \begin{bmatrix}
        \hat{\pi}^1\\
        \vdots\\
        \hat{\pi}^M
    \end{bmatrix}
    \in\BR^{N\times n},
    \qquad N=Mn.
\]
We use $\hat Z\in\BR^{n\times d_z}$ to denote the latent embedding matrix whose $\ell$-th row is
$z_\ell^\top$, and define $
    \hat Z_{\mathrm{norm}}
    =
    (\|z_1\|^2,\ldots,\|z_n\|^2)^\top$ as the row-norm vector of $\hat Z$. Define
\[
    \Omega
    =
    \frac{1}{n}\operatorname{Rowsum}(LL^\top)
    =
    \frac{1}{n}L(L^\top\vecone_N)
    \in\BR^N,
\]
where $\vecone_N$ is the all-one vector in $\BR^N$. The following result rewrites the quadratic low-rank empirical objective $\hat S_{\epsilon}^{\hat R}$ as a minimization over an auxiliary matrix and $M$ EOT subproblems.

\begin{theorem}[Low-rank sample kernel EGW lifting]
\label{thrm:sample_duality}
Fix $\epsilon>0$. Suppose Assumption~\ref{assumption2.1} holds, then
\[
\hat S_{\epsilon}^{\hat R}
\big(\{\hat{\mu}_n^i\}_{i=1}^M,\hat{\nu}_n\big)
=
\min_{\hat A\in\BR^{\hat R\times d_z}}
\left\{
    2\|\hat A\|_F^2
    +
    \sum_{i=1}^M
    \OT^i_{\hat C^i(\hat A),\epsilon}
    (\hat{\mu}_n^i,\hat{\nu}_n)
\right\},
\]
where $\hat C^i(\hat A)\in\BR^{n\times n}$ is the $i$-th row block of
\[
    \hat C(\hat A)
    =
    2\Omega \hat Z_{\mathrm{norm}}^\top
    -
    4L\hat A \hat Z^\top
    \in\BR^{N\times n}.
\]
Equivalently, $\hat C^i(\hat A)$ consists of rows $((i-1)n+1):in$ of $\hat C(\hat A)$. For each modality $i$, the EOT subproblem is
\[
    \OT^i_{\hat C^i(\hat A),\epsilon}
    (\hat{\mu}_n^i,\hat{\nu}_n)
    =
    \min_{\hat{\pi}^i\in\Gamma(\hat{\mu}_n^i,\hat{\nu}_n)}
    \left\{
        \langle \hat C^i(\hat A),\hat{\pi}^i\rangle
        +
        \epsilon
        D_{\mathrm{KL}}\!\left(
            \hat{\pi}^i
            \,\middle\|\,
            \hat{\mu}_n^i\otimes\hat{\nu}_n
        \right)
    \right\}.
\]
Moreover, for fixed transport plans $\{\hat{\pi}^i\}_{i=1}^M$, the minimizing
auxiliary variable satisfies
\[
    \hat A
    =
    L^\top \hat{\Pi} \hat Z .
\]
\end{theorem}

The proof of Theorem~\ref{thrm:sample_duality} is deferred to Section~\ref{proof:sample_duality}, which follows by applying the similar variational lifting mechanisms as in Theorem~\ref{thrm:pop_duality}. 
Importantly, Theorem~\ref{thrm:sample_duality} suggests an alternating minimization scheme to optimize~\Cref{eq:emp_sum_obj}, which proceeds in iterations. For each iteration, we first fix $\hat A$ and optimize the transport plans with linear cost calculated by $\hat C(\hat A)$. This subroutine decouples into $M$ EOT problems that can be solved by off-the-shelf EOT solvers, such as the Sinkhorn algorithm~\citep{cuturi2013sinkhorn}. We then fix transport plans and optimize the auxiliary variable $\hat A$, which has a closed-form solution $\hat A=L^\top\hat\Pi \hat Z$. The algorithm proceeds by iterations until convergence. We summarize this procedure in Algorithm~\ref{alg:alt_opt}, where the {\color{RoyalPurple} purple} lines execute the sub-optimization for the transport plans $\hat\Pi$ and the {\color{Maroon} maroon} lines optimize $\hat A$.





\begin{algorithm}[ht]
\caption{Alternating minimization algorithm for solving $\hat S_{\epsilon}^{\hat R}$}
\label{alg:alt_opt}
\begin{algorithmic}[1]
\REQUIRE Finite-rank kernel factor $L \in \mathbb{R}^{N \times \hat R}$, latent matrix $\hat Z \in \mathbb{R}^{n \times d_z}$, initial auxiliary matrix $\hat A^{(0)}\in\mathbb R^{\hat R\times d_z}$, regularization parameter $\epsilon>0$, tolerance $\tau>0$, maximum iterations $T_{\max}$
\ENSURE Transport plan $\hat{\Pi}$, auxiliary matrix $\hat A$, approximate objective value $\mathcal J$

\STATE Compute $\hat Z_{\mathrm{norm}} \leftarrow (\|z_1\|^2,\ldots,\|z_n\|^2)^\top$
\STATE Compute $\Omega \leftarrow \frac{1}{n}L(L^\top\vecone_N)$
\STATE Initialize $\mathcal J^{(0)} \leftarrow +\infty$

\FOR{$t=1,\ldots,T_{\max}$}
    {\color{RoyalPurple}\STATE $\hat C^{(t)} \leftarrow 2\Omega \hat Z_{\mathrm{norm}}^\top - 4L\hat A^{(t-1)}\hat Z^\top$
    \FOR{$i=1,\ldots,M$}
        \STATE Let $\hat C^{i,(t)}$ be the $i$-th row block of $\hat C^{(t)}$
        \STATE $\hat\pi^{i,(t)} \leftarrow \mathrm{Sinkhorn}\!\left(\hat C^{i,(t)},\epsilon;\frac{1}{n}\vecone,\frac{1}{n}\vecone\right)$
    \ENDFOR
    \STATE $\hat\Pi^{(t)} \leftarrow
    \begin{bmatrix}
        \hat\pi^{1,(t)}\\
        \vdots\\
        \hat\pi^{M,(t)}
    \end{bmatrix}$}
    {\color{Maroon} \STATE $\hat A^{(t)} \leftarrow L^\top \hat\Pi^{(t)}\hat Z$

    \STATE $\mathcal J^{(t)} \leftarrow
    2\|\hat A^{(t)}\|_F^2
    +
    \left\langle
        2\Omega \hat Z_{\mathrm{norm}}^\top - 4L\hat A^{(t)}\hat Z^\top,
        \hat\Pi^{(t)}
    \right\rangle
    +
    \epsilon\sum_{i=1}^M
    D_{\mathrm{KL}}\!\left(
        \hat\pi^{i,(t)}
        \,\middle\|\,
        \hat\mu_n^i\otimes\hat\nu_n
    \right)$}

    \IF{$t>1$ and $\frac{|\mathcal J^{(t-1)}-\mathcal J^{(t)}|}{1+|\mathcal J^{(t-1)}|}<\tau$}
        \STATE \textbf{break}
    \ENDIF
\ENDFOR

\STATE \textbf{return} $(\hat\Pi^{(t)},\hat A^{(t)},\mathcal J^{(t)})$
\end{algorithmic}
\end{algorithm}

Algorithm~\ref{alg:alt_opt} implements the Sinkhorn updates~\citep{cuturi2013sinkhorn} to solve the modality-wise EOT subproblems whose costs are determined by current auxiliary matrix $\hat A$. Each Sinkhorn update costs $O(n^2)$ per transport plan. The total per-iterate cost of Sinkhorn update is thus $O(Mn^2)$. 

\subsection{Error Analysis}
We can now relate the computable low-rank empirical objective in~\Cref{eq:emp_sum_obj} back to the population objective in~\Cref{eq:main_objective}, whose error is the ultimate quantity we care about. Similar to Theorem~\ref{cor:full}, we can analyze the error bound of running Algorithm~\ref{alg:alt_opt} to solve the surrogate objective of Eq~\eqref{eq:main_objective}. More specifically, we introduce the Theorem~\ref{cor:end_to_end_computable_error} to bound the expectation between $S_\epsilon$ and $\hat S_\epsilon^{\hat R}$. 

\begin{theorem}[End-to-end error for the computable low-rank objective]
\label{cor:end_to_end_computable_error}
Suppose Assumption~\ref{assumption2.1} holds and assume that $\nu$ is $4$-sub-Weibull with parameter $\sigma^2$. Assume also that the latent support points $z_1,\ldots,z_n$ are sampled i.i.d.\ from $\nu$. Let $
    q_{R,z}
    =
    \left\lceil\frac{R\vee d_z}{2}\right\rceil$. Then, for every $\epsilon>0$,
\[
\begin{aligned}
    &
    \mathbb E
    \left[
    \left|
    S_\epsilon(\{\mu^i\}_{i=1}^M,\nu)
    -
    \hat S_\epsilon^{\hat R}(\{\hat\mu_n^i\}_{i=1}^M,\hat\nu_n)
    \right|
    \right] \\
    &\qquad\lesssim_{M,R,d_z,\kappa}
    \epsilon
    \left[
        1+\left(\frac{\sigma}{\epsilon}\right)+\left(\frac{\sigma}{\epsilon}\right)^{2q_{R,z}+2}+
        \left(\frac{\sigma}{\epsilon}\right)^{4q_{R,z}+4}
    \right]
    \frac1{\sqrt n} \\
    &\qquad\quad
    +
    \sigma
    \left(
        2\eta_R
        +\mathbb E[\rho_{\hat R}]
    \right),
\end{aligned}
\]
where the expectation over $\rho_{\hat R}$ is taken over modality samples. In particular, if the Pivoted Cholesky/Nystr\"om approximation satisfies $\rho_{\hat R}\le\tau_{\hat R}$ almost surely, then
\[
\begin{aligned}
    &
    \mathbb E
    \left[
    \left|
    S_\epsilon(\{\mu^i\}_{i=1}^M,\nu)
    -
    \hat S_\epsilon^{\hat R}(\{\hat\mu_n^i\}_{i=1}^M,\hat\nu_n)
    \right|
    \right] \\
    &\qquad\lesssim_{M,R,d_z,\kappa}
    \epsilon
    \left[
        1+\left(\frac{\sigma}{\epsilon}\right)+\left(\frac{\sigma}{\epsilon}\right)^{2q_{R,z}+2}+\left(\frac{\sigma}{\epsilon}\right)^{4q_{R,z}+4}
    \right]
    \frac1{\sqrt n} \\
    &\qquad\quad+\sigma
    \left(
        2\eta_R+\tau_{\hat R}
    \right).
\end{aligned}
\]
\end{theorem}

\subsection{Proof sketch}
The bound in Theorem~\ref{cor:end_to_end_computable_error} can be decomposed into several components via triangle inequalities: 
\begin{align}
&
\left|
S_\epsilon(\{\mu^i\}_{i=1}^M,\nu)
-
\hat S_\epsilon^{\hat R}(\{\hat\mu_n^i\}_{i=1}^M,\hat\nu_n)
\right| \notag \\
&\le
\left|
S_\epsilon(\{\mu^i\}_{i=1}^M,\nu)
-
S_\epsilon^R(\{\mu^i\}_{i=1}^M,\nu)
\right| \tag{Truncation error on population distributions} \\
&+
\left|
S_\epsilon^R(\{\mu^i\}_{i=1}^M,\nu)
-
S_\epsilon^R(\{\hat\mu_n^i\}_{i=1}^M,\hat\nu_n)
\right| \tag{Statistical error on rank-$R$ JK-EGW}  \\
&+
\left|
S_\epsilon^R(\{\hat\mu_n^i\}_{i=1}^M,\hat\nu_n)
-
\hat S_\epsilon^{\hat R}(\{\hat\mu_n^i\}_{i=1}^M,\hat\nu_n)
\right|  \tag{Empirical analytic-to-computable error} .
\end{align}

Here, the first term is the population kernel truncation error bounded by Proposition~\ref{prop:kernel_truncation_error}, the second term is the rank-$R$ statistical error bounded by Proposition~\ref{thrm:sample_complexity}. It remains to bound the third term, which is the empirical analytic-to-computable error. The following proposition gives a direct comparison between the analytic rank-R empirical objective and the computable rank-$\hat R$ empirical objective (``analytic'' and ``computable'' refer to the kernel). 

\begin{proposition}[Empirical analytic-to-computable kernel approximation error]
\label{cor:empirical_analytic_to_computable}
Suppose Assumption~\ref{assumption2.1} holds and $\nu$ is $4$-sub-Weibull with parameter $\sigma^2$. Then, for every $\epsilon>0$,
\[
\begin{aligned}
    &
    \mathbb E \left[ \left|
    S_\epsilon^R(\{\hat\mu_n^i\}_{i=1}^M,\hat\nu_n)
    -
    \hat S_\epsilon^{\hat R}(\{\hat\mu_n^i\}_{i=1}^M,\hat\nu_n)
    \right| \right] \\
    &\qquad\le
    M^2\sigma(\eta_R+\mathbb E[\rho_{\hat R}]),
\end{aligned}
\]
where the expectation over $\rho_{\hat R}$ is taken over modality samples. In particular, if $\rho_{\hat R}\le\tau_{\hat R}$ almost surely, then
\[
    \mathbb E
    \left[
    \left|
    S_\epsilon^R(\{\hat\mu_n^i\}_{i=1}^M,\hat\nu_n)
    -
    \hat S_\epsilon^{\hat R}(\{\hat\mu_n^i\}_{i=1}^M,\hat\nu_n)
    \right|
    \right]
    \lesssim
    M^2\sigma\left(\eta_R+\tau_{\hat R}\right).
\]
\end{proposition}

To prove Proposition~\ref{cor:empirical_analytic_to_computable}, we rely on an additional layer of triangle inequality:
\[
\begin{aligned}
&
\left|
S_\epsilon^R(\{\hat\mu_n^i\}_{i=1}^M,\hat\nu_n)
-
\hat S_\epsilon^{\hat R}(\{\hat\mu_n^i\}_{i=1}^M,\hat\nu_n)
\right| \\
&\qquad\le
\left|
S_\epsilon^R(\{\hat\mu_n^i\}_{i=1}^M,\hat\nu_n)
-
\hat S_\epsilon(\{\hat\mu_n^i\}_{i=1}^M,\hat\nu_n)
\right|
+
\left|
\hat S_\epsilon(\{\hat\mu_n^i\}_{i=1}^M,\hat\nu_n)
-
\hat S_\epsilon^{\hat R}(\{\hat\mu_n^i\}_{i=1}^M,\hat\nu_n)
\right|,
\end{aligned}
\]
where the first term is the discrepancy due to replacing $\Phi_R(x)^\top\Phi_R(x)$ by $W$ on the empirical sample pairs. More specifically, the following proposition quantifies this error in terms of the Mercer-tail quantity $\eta_R$ and the empirical second moment of the latent support. 

\begin{proposition}[Empirical Mercer realization error]
\label{prop:empirical_mercer_realization}
Suppose Assumption~\ref{assumption2.1} holds and $\nu$ is $4$-sub-Weibull with parameter $\sigma^2$. Then, for every $\epsilon>0$,
\[
\begin{aligned}
    \mathbb E \left[ \left|
    S_\epsilon^R(\{\hat\mu_n^i\}_{i=1}^M,\hat\nu_n)-
    \hat S_\epsilon(\{\hat\mu_n^i\}_{i=1}^M,\hat\nu_n)
    \right| \right]\le
    M^2\sigma\,\eta_R.
\end{aligned}
\]
\end{proposition}

Recall that there are two distinct kernel approximation mechanisms. The first is replacing the population kernel $K$ by the rank-$R$ Mercer truncation $\Phi_R^\top\Phi_R$, which is bounded in Proposition~\ref{prop:kernel_truncation_error}. The other error comes from replacing the empirical realization $W$ by the Pivoted Cholesky/Nystr\"om factorization $LL^\top$, which contributes to the second term of the triangle inequality above. More specifically, the second term is bounded by the following Proposition:

\begin{proposition}[Pivoted Cholesky empirical approximation error]
\label{prop:pivoted_cholesky_empirical_error}
Suppose that $\nu$ is $4$-sub-Weibull with parameter $\sigma^2$. Then, for every $\epsilon>0$,
\[
\begin{aligned}
    \mathbb E
    \left[
    \left|
    \hat S_\epsilon(\{\hat\mu_n^i\}_{i=1}^M,\hat\nu_n)
    -\hat S_\epsilon^{\hat R}(\{\hat\mu_n^i\}_{i=1}^M,\hat\nu_n)
    \right|
    \right] \lesssim
    M^2\sigma\,\mathbb E[\rho_{\hat R}],
\end{aligned}
\]
where the expectation over $\rho_{\hat R}$ is taken over modality samples. In particular, if the Pivoted Cholesky routine returns a factor $L$ satisfying
$\rho_{\hat R,n}\le\tau_{\hat R}$ almost surely, then
\[
    \mathbb E
    \left[
    \left|
    \hat S_\epsilon(\{\hat\mu_n^i\}_{i=1}^M,\hat\nu_n)
    -
    \hat S_\epsilon^{\hat R}(\{\hat\mu_n^i\}_{i=1}^M,\hat\nu_n)
    \right|
    \right]
    \lesssim
    M^2\sigma\,\tau_{\hat R}.
\]
\end{proposition}

The proofs of Proposition~\ref{prop:empirical_mercer_realization} and Proposition~\ref{prop:pivoted_cholesky_empirical_error} are deferred to Section~\ref{pf:prop:empirical_mercer_realization} and Section~\ref{pf:pivoted_cholesky_empirical_error}, respectively. 

\section{Numerical Results}
We present the numerical analysis of the proposed JK-EGW in this section. We will first use multi-modality MNIST as a proof-of-concept experiment to illustrate the concept of multimodal alignment and the performance of JK-EGW. Then, we will validate the statistical error on a synthetic example. Finally, we implement Algorithm~\ref{alg:alt_opt} on an image-text alignment task to investigate its performance on real multimodal datasets.  

\subsection{Multi-modality Numerals }

 The purpose of this example is to visualize how the transport plans learned by Algorithm~\ref{alg:alt_opt} map heterogeneous feature representations into a shared latent reference. The multi-modality numerals is one UCI dataset for the multi-feature representation of handwritten numerals from 0 to 9 ~\citep{multiple_features_72}.  For numerical digits, multiple features are constructed through several feature sets. We take two of them as the $M=2$ modalities. The modality first consists of 76 Fourier coefficients of the character shapes, which describe the geometry of each digit in the frequency domain. The second consists of 64 Karhunen--Lo\`eve coefficients, which represent each digit image through its projections onto principal components computed from the images. The two feature sets take values in Euclidean spaces of different dimensions, and their coordinates often does not have common interpretation. Therefore, pointwise comparisons of features are not well defined across modalities. JK-GWOT provides an alternative tool to align features from the two modalities.

We construct the within-modality similarity blocks $W^{i,i}$ using the RBF kernels. The cross-modality similarity block $W^{i,j}$ is defined by label agreement. More specifically, we set $W^{i,j}_{k,k'}=1$ when $x_k^i$ and $x_{k'}^j$ have the same digit label and $W^{i,j}_{k,k'}=0$ otherwise. The shared latent support is sampled from a two-dimensional Gaussian reference distribution for visualization. The displayed
point for sample $x_k^i$ is its barycentric latent embedding $\hat y_k^i = n\sum_{\ell=1}^n \hat\pi_{k\ell}^i z_\ell$. We use the POT library~\citep{flamary2021pot} that implements~\citet{cuturi2013sinkhorn} to compute Sinkhorn updates in Algorithm~\ref{alg:alt_opt}.

\begin{figure}[h]
    \centering
    \includegraphics[width=1\linewidth]{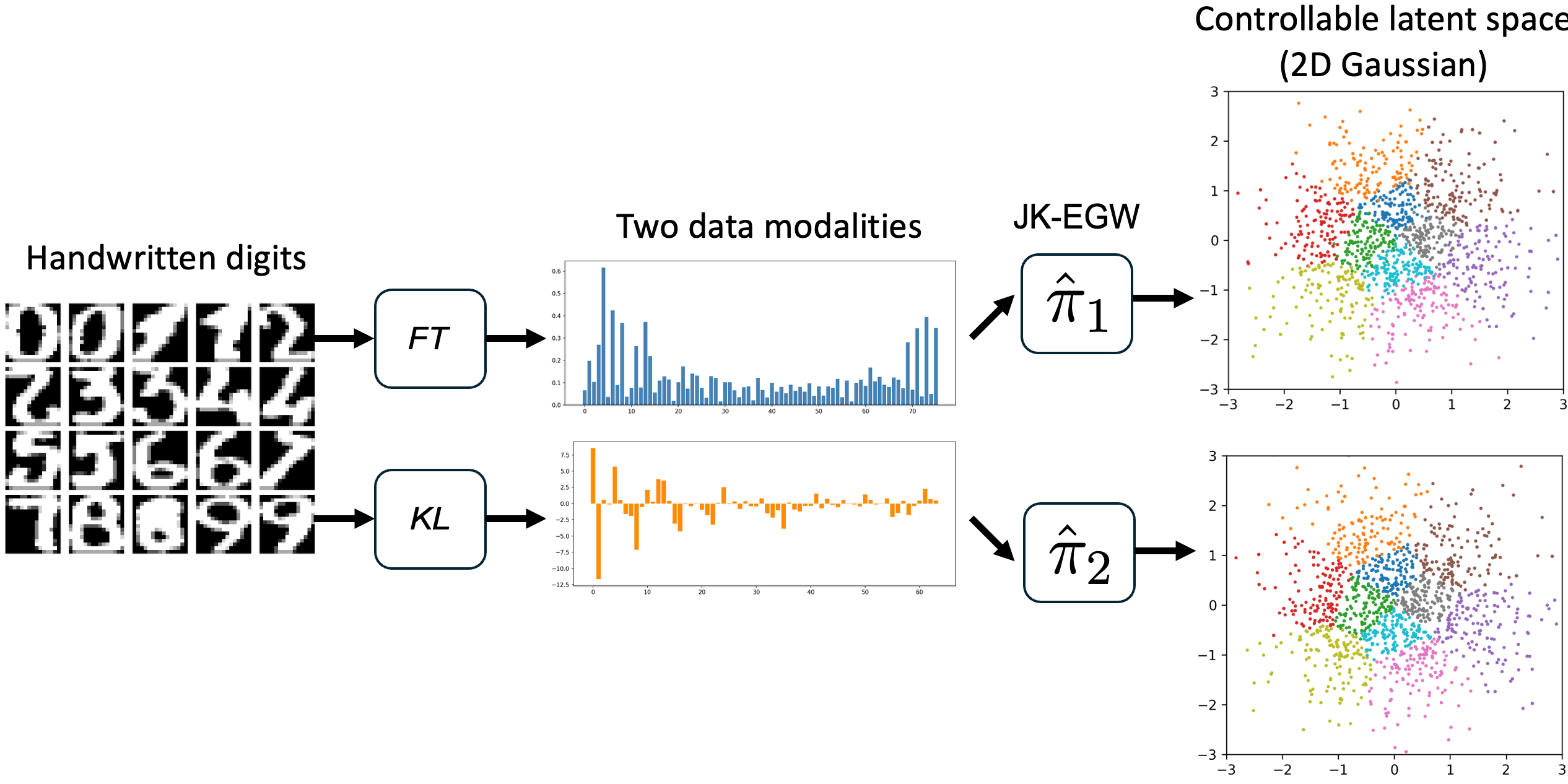}
    \caption{Worked illustration of Algorithm~\ref{alg:alt_opt} on a multi-view handwritten digits dataset. Two feature views, Fourier coefficients (FT) and Karhunen--Lo\`eve coefficients (KL), are treated as two modalities and coupled to a shared two-dimensional Gaussian latent reference. Points show the barycentric latent embeddings induced by the learned transport plans and are colored by digit class.}
    \label{fig:handwritten-digits}
\end{figure}

Figure~\ref{fig:handwritten-digits} illustrates the output of Algorithm~\ref{alg:alt_opt}. The resulting embeddings place samples with the same digit label in comparable regions of the latent space across the two views, illustrating how the optimized couplings encode both within-modality geometry and cross-modality affinity.

\subsection{Synthetic Example}
\begin{figure}[h]
    \centering
    \includegraphics[width=0.45\linewidth]{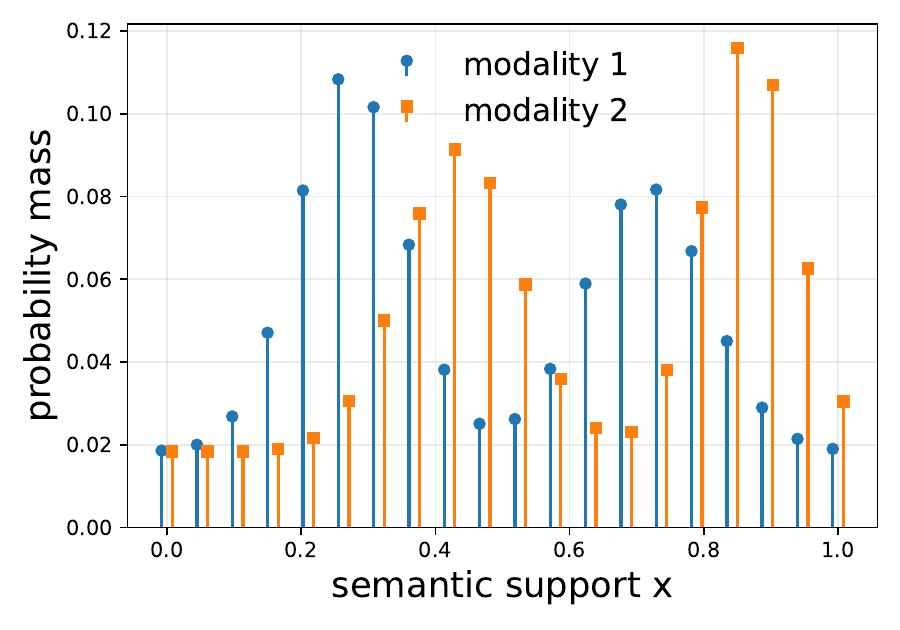}
    \includegraphics[width=0.45\linewidth]{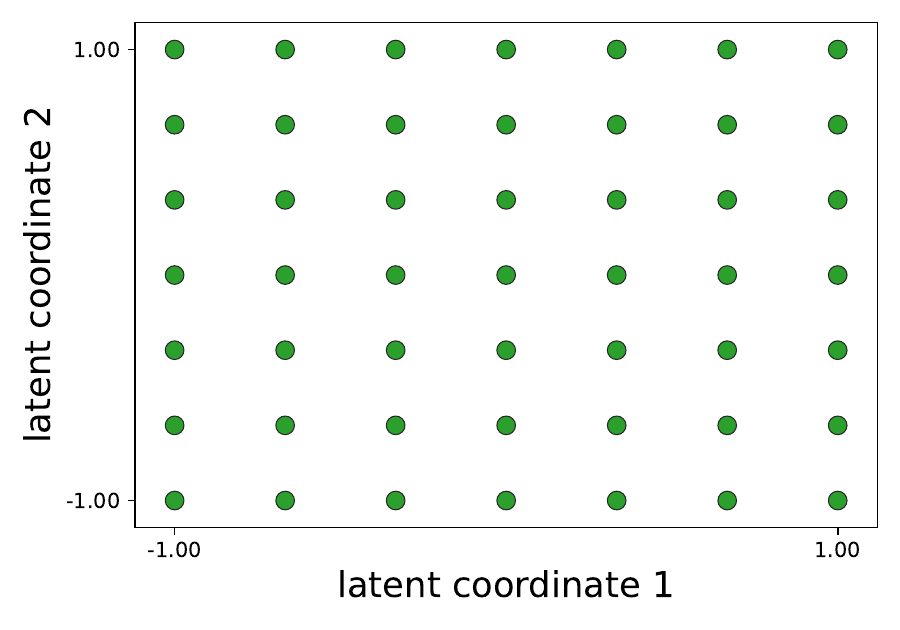}

    \caption{ The feature space follow two discrete mixtures of Gaussian distributions. The latent distribution is uniform on a 2D grid. }    \label{fig:theories_distribution_plots}
\end{figure}

To validate the theoretical findings in Theorem~\ref{cor:full}, we construct a synthetic problem in which all distributions have finite support, so that the population objective $S_\epsilon(\{\mu^i\}_{i=1}^M,\nu)$ is itself finite-dimensional. In this setting, the exact Mercer feature map of $K$ is available through eigendecomposition of its Gram matrix, and running Algorithm~\ref{alg:alt_opt} on the population weights returns the ground-truth value up to solver tolerances. There are $M=2$ modalities. The marginal distribution of each is supported on a discrete 1D set of a uniform grid of $20$ points in $[0,1]$. The distributions are created from discrete mixture of Gaussian and visualized in Figure~\ref{fig:theories_distribution_plots}. The latent reference $\nu$ is the uniform distribution on a $7\times 7$ grid in $[-1,1]^2$. Such discrete distribution is $4$-sub-Weibull. We use the RBF kernel $K(x,x')=\gamma\exp\bigl(-\norm{x-x'}^2/2\ell^2\bigr)$ with bandwidth $\ell=0.1$ as the similarity kernels, where $\gamma=1$ for within-modality affinity and $\gamma=0.6$ for cross-modality affinity. $K$ is positive semidefinite for $\gamma\in[0,1]$, and $K(x,x)=1$ gives $\kappa=1$. We implement Algorithm~\ref{alg:alt_opt} where the optimization over transport plans is solved by log-domain Sinkhorn~\citep{cuturi2013sinkhorn, flamary2021pot}.

\begin{figure}[h]
    \centering
    \includegraphics[width=0.46\linewidth]{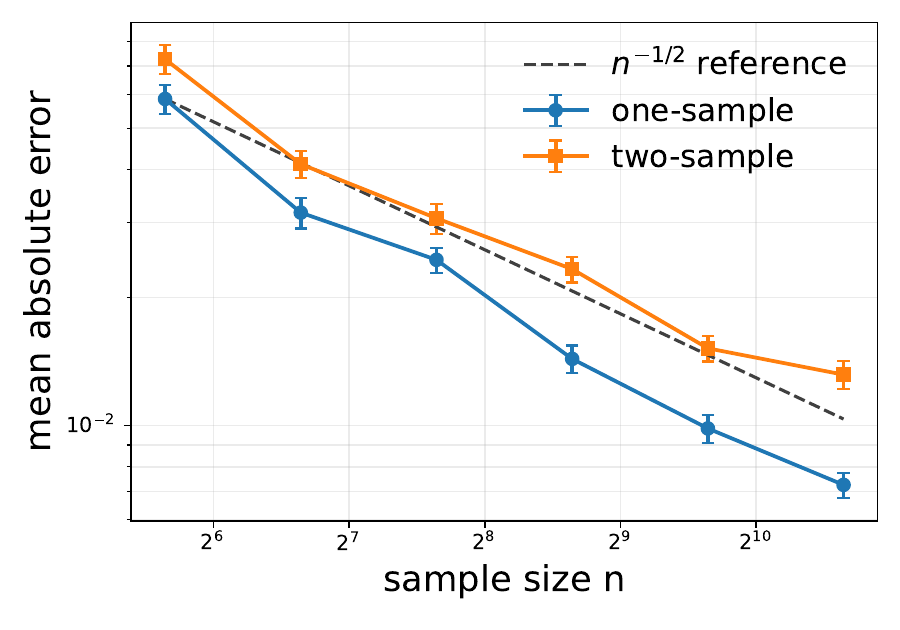}
    \includegraphics[width=0.46\linewidth]{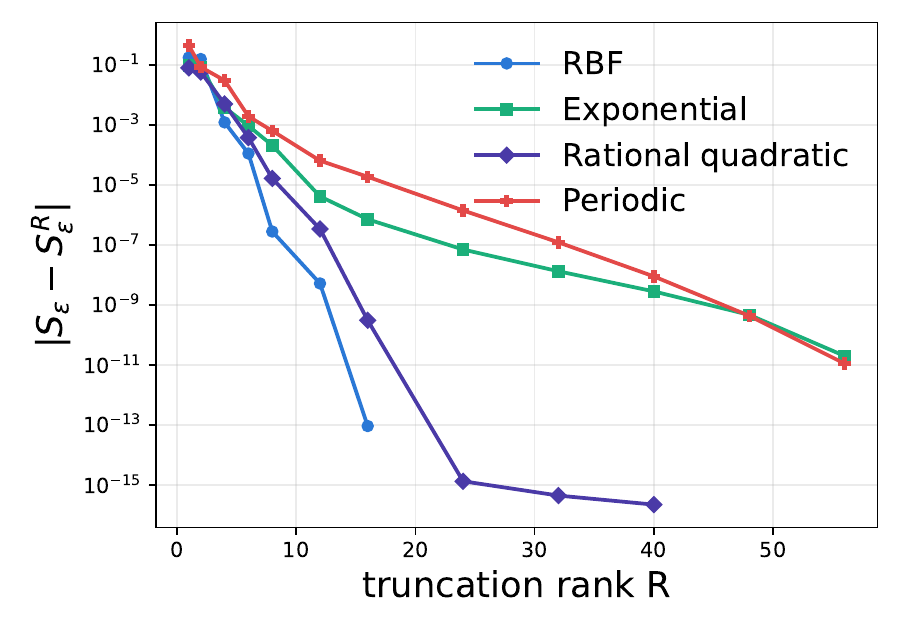}
    \caption{ \textbf{Left}: Sample complexity of the full-kernel entropic objective. Mean absolute error of the plug-in estimator versus $n$ (log--log), in the one-sample ($\nu$ known) and two-sample ($\nu$ estimated) regimes. \textbf{Right}: Objective truncation gap $|S_\epsilon - S_\epsilon^R|$ v.s. Mercer rank $R$.  Objective truncation gaps $|S_\epsilon - S_\epsilon^R|$ for four kernels. 
    Both curves follow the dashed $n^{-1/2}$ reference. }
    \label{fig:theories_plots}
\end{figure}

\noindent \underline{\textit{Parametric rate}}. For each $n\in\{50,100,200,400,800,1600\}$, we sample and construct empirical measures $\hat\mu^1_n,\hat\mu^2_n$, and $\hat\nu_n$, independently across modalities and the latent reference, and calculate $S_\epsilon(\{\hat\mu^i_n\},\hat\nu_n)$ and $S_\epsilon(\{\mu^i\},\nu)$ via Algorithm~\ref{alg:alt_opt}. Since the population distributions $\mu^1$, $\mu^2$, and $\nu$ are also discrete, Algorithm~\ref{alg:alt_opt} directly yields the population value $S_\epsilon({\mu^i},\nu)$. We record the two-sample error $|S_\epsilon(\{\hat\mu^i_n\},\hat\nu_n)-S_\epsilon(\{\mu^i\},\nu)|$. Similarly, in the one-sample regime, we calculate $S_\epsilon(\{\hat\mu^i_n\},\nu)$ using Algorithm~\ref{alg:alt_opt} and record the statistical error $|S_\epsilon(\{\hat\mu^i_n\},\nu)-S_\epsilon(\{\mu^i\},\nu)|$.  We run both experiments over $100$ Monte Carlo repetitions and plot the averaged error curves on Figure~\ref{fig:theories_plots} (\textbf{Left}) on log--log axes with the $n^{-1/2}$ reference. Both curves show trends consistent with the $O(n^{-1/2})$ upper bound provided in Theorem~\ref{cor:full}.

\noindent \underline{\textit{Kernel truncation error}}. In addition, we explore how the choice of $R$ affect kernel truncation error $\eta_R$ thus $|S_\epsilon-S_\epsilon^R|$. We experiment with 4 different types of common kernels, including the RBF kernel $K(x,x')=\exp\!\big(-\norm{x-x'}^2/2\ell^2\big)$, the exponential kernel $K(x,x')=\exp\!\big(-\norm{x-x'}/\ell\big)$, the rational quadratic kernel $K(x,x')=\big(1+\norm{x-x'}^2/2\ell^2\big)^{-1}$, and the periodic kernel $K(x,x')=\exp\!\big(-2\sin^2(\pi\norm{x-x'})/\ell^2\big)$. We set a common length scale $\ell=0.15$ to all the kernels that approximately matches the geometry in the feature space. We still use Algorithm~\ref{alg:alt_opt} to solve for the transport plan with the entropy regularization parameter $\epsilon=0.2$. Since the $\mu^1$ and $\mu^2$ have finite support, we can obtain the \textit{exact} Mercer feature map $\Phi$ by eigendecomposition of the kernel Gram matrix. Similarly, the rank-$R$ features $\Phi_R$ are formed by by retaining the top $R$ leading eigenpairs of the kernel Gram matrix. Figure~\ref{fig:theories_plots} (\textbf{Right}) reports the gap $|S_\epsilon - S_\epsilon^R|$ against different $R$. Since $M$ and $\sigma$ are the same for all kernel experiments, the errors are dominated by the Mercer-tail quantity $\eta_R$. We can observe that for all kernels we experimented, the gap in transport value decreases with the increase of $R$, which is consistent with the theoretical prediction in Theorem~\ref{cor:full} that smaller $\eta_R$ results in smaller error between $S_\epsilon$ and $S_\epsilon^R$. Among the 4 kernels, the smooth RBF and rational quadratic kernels reach solver precision faster, whereas the exponential and periodic kernels require larger $R$.

It is worth noting that in both synthetic experiments, the reported objective values are subject to optimization error from Algorithm~\ref{alg:alt_opt}, which cannot be avoided in numerical experiments. Nevertheless, we assume the optimization error is negligible compared to the statistical error in practice and use the algorithm output to validate Proposition~\ref{prop:kernel_truncation_error} and Theorem~\ref{cor:full}.

\subsection{Image-text Alignment}\label{sec:experiments}
We evaluate JK-EGW alignment on image--text alignment tasks constructed from MS--COCO~\citep{lin2014microsoft}. The experiments are designed to assess whether the proposed JK-EGW formulation learns a geometrically coherent shared latent space. 

\paragraph{Dataset.}
We construct image--text subsets from the MS--COCO 2017 training split~\citep{lin2014microsoft}, pairing each image with the concatenation of its five reference captions. To obtain controlled alignment benchmarks with class-level structure, we use a class-targeted sampling strategy. We first select a set of MS--COCO object categories and then sample a balanced number of images per category. We instantiate this strategy in two configurations with different sizes of the training dataset: a 6K subset with 6000 training and 3000 test images drawn from 20 categories (300 training and 150 test images per category), and a larger 9K subset with 9000 training and 3000 test images drawn from 30 categories (300 training and 100 test images per category). Within each subset, each image is assigned to at most one target category, and the training and test splits are disjoint.


\paragraph{Procedure.}
We use MPNet~\citep{song2020mpnet} as the frozen text encoder and DINOv3~\citep{simeoni2025dinov3} as the frozen image encoder. All features are extracted once and $\ell_2$-normalized before alignment and evaluation. The distribution $\nu$ is instantiates as a Gaussian mixture model with 20 components. Each component follows an isotropic Gaussian, whose means are placed along mutually orthogonal directions at a fixed distance from the origin. This is an informative choice of marginal distribution which does not require the prior knowledge of of the common representation space. We set the dimension of the common representations space $d=32$.

Algorithm~\ref{alg:alt_opt} returns one transport plan per modality, $\hat\pi^{\mathrm{txt}}$ and $\hat\pi^{\mathrm{img}}$, coupling the training samples of each modality to the latent support matrix $\hat Z$. We define the per-sample latent embeddings of the training set as $\hat{Z}^{\mathrm{txt}}_{\mathrm{train}} = \hat\pi^{\mathrm{txt}} \hat Z$ and $\hat{Z}^{\mathrm{img}}_{\mathrm{train}} = \hat\pi^{\mathrm{img}} \hat Z$, whose rows are, up to a uniform scaling, the barycentric latent embeddings introduced in Section~\ref{sec:optimization}. To extend this map out of the training set, we train a simple MLP for each modality that regresses the raw encoder features $X_{\mathrm{train}}$ onto the ($\ell_2$-normalized) training $\hat{Z}_{\mathrm{train}}$. Model selection uses a small validation split held out from the training set, based on validation retrieval recall at 1 (R@1). Test samples are used only at inference. The trained MLPs are then frozen and applied to the $\ell_2$-normalized test features, so each test point is embedded in the shared latent space by a learned approximation of the training transport map. 

\paragraph{Baselines.}
We compare \name with three most recent baseline methods for aligning features generated by pretrained single-modality encoders. CSA~\citep{licsa} is a closed-form canonical similarity analysis method that maps frozen unimodal features into a shared multimodal space via canonical correlation analysis. CSA+$\textit{STRUCTURE}$ augments CSA with the similarity-based layer-selection rule and STRUCTURE regularization in~\citet{groger2025limited}. MLP+$\textit{STRUCTURE}$ replaces the linear map in CSA with a trainable MLP alignment head optimized using a contrastive objective~\citep{clip} with similarity-based layer selection and STRUCTURE regularization~\citep{groger2025limited}.

\paragraph{Alignment.}
For MS--COCO alignment, we pool the aligned image and caption embeddings from the test split and report purity and silhouette scores using the target-category labels. Purity measures the agreement between a $k$-means clustering of the pooled embeddings and the target categories, while the silhouette score measures whether each pooled embedding lies closer to its own category than to the nearest other category under the cosine distance. They characterize how well the learned embeddings conform with the semantic structure. To evaluate cross-modal retrieval, we report recall at 1 (R@1) and recall at 5 (R@5). More specifically, retrieval is calculated based on the cosine similarity between the aligned test embeddings of the two modalities. In the image-to-text direction, each test image ranks all test caption embeddings by decreasing cosine similarity to its own embedding, and $r_i^{\mathrm{img}\to\mathrm{txt}}$ denotes the resulting rank of the caption paired with image $i$. The rank $r_j^{\mathrm{txt}\to\mathrm{img}}$ is defined symmetrically for each test caption $j$. For $k\in\{1,5\}$, the image-to-text recall is
\[
    \mathrm{R@}k^{\mathrm{img}\to\mathrm{txt}}
    =
    \frac{1}{n_{\mathrm{ev}}}
    \sum_{i=1}^{n_{\mathrm{ev}}}
    \vecone_{\{r_i^{\mathrm{img}\to\mathrm{txt}}\le k\}},
\]
where $n_{\mathrm{ev}}$ denotes the number of test pairs. The text-to-image recall $\mathrm{R@}k^{\mathrm{txt}\to\mathrm{img}}$ is defined analogously, and the reported value is
\[
    \mathrm{R@}k
    =
    \tfrac{1}{2}
    \left(
    \mathrm{R@}k^{\mathrm{img}\to\mathrm{txt}}
    +
    \mathrm{R@}k^{\mathrm{txt}\to\mathrm{img}}
    \right).
\]

\begin{table*}[h]
\centering
\caption{Alignment and retrieval performance of \name and baselines.}
\begin{tabular}{clcccc}
\toprule
   \multirow{2}{*}{Dataset} & \multirow{2}{*}{Model}  & \multicolumn{2}{c}{Alignment} & \multicolumn{2}{c}{Retrieval} \\
     &  & \textbf{Purity} & \textbf{Silhouette} & \textbf{R@1} & \textbf{R@5} \\
\hline
    \multirow{4}{*}{MS--COCO 6k} & CSA & 0.638 & -0.210 & 0.257 & 0.542 \\
    & CSA+$\textit{STRUCTURE}$ & 0.654 & -0.196 & 0.290 & 0.595 \\
    & MLP+$\textit{STRUCTURE}$ & 0.638 & -0.156 & 0.154 & 0.390 \\
    & \name(Ours) & \textbf{0.719} & \textbf{0.117} & \textbf{0.380} & \textbf{0.674} \\
\hline
    \multirow{4}{*}{MS--COCO 9k} & CSA & 0.673 & -0.150 & 0.263 & 0.549 \\
    & CSA+$\textit{STRUCTURE}$ & \textbf{0.689} & -0.143 & 0.304 & 0.613 \\
    & MLP+$\textit{STRUCTURE}$ & 0.638 & -0.133 & 0.084 & 0.239 \\
    & \name(Ours) & 0.633 & \textbf{0.104} & \textbf{0.413} & \textbf{0.706} \\
\hline
\end{tabular}
\label{tab:coco6k}
\end{table*}

Table~\ref{tab:coco6k} reports results for \name and three baseline alignment methods on the test splits of the 6K and 9K MS--COCO subsets. \name produces the strongest alignment and retrieval results across all dimensions. It achieves the highest purity in both datasets and is the only method with positive silhouette scores, while all baselines yield negative silhouette values. These results suggest that \name learns a shared representation that produces improved retrieval performance. We further conduct an ablation study on the latent dimension $d$ in Section~\ref{sec:additiona_exp} of the Appendix.


\begin{figure}[ht]
    \centering
    \begin{subfigure}[b]{0.46\linewidth}
        \centering
        \includegraphics[width=\linewidth]{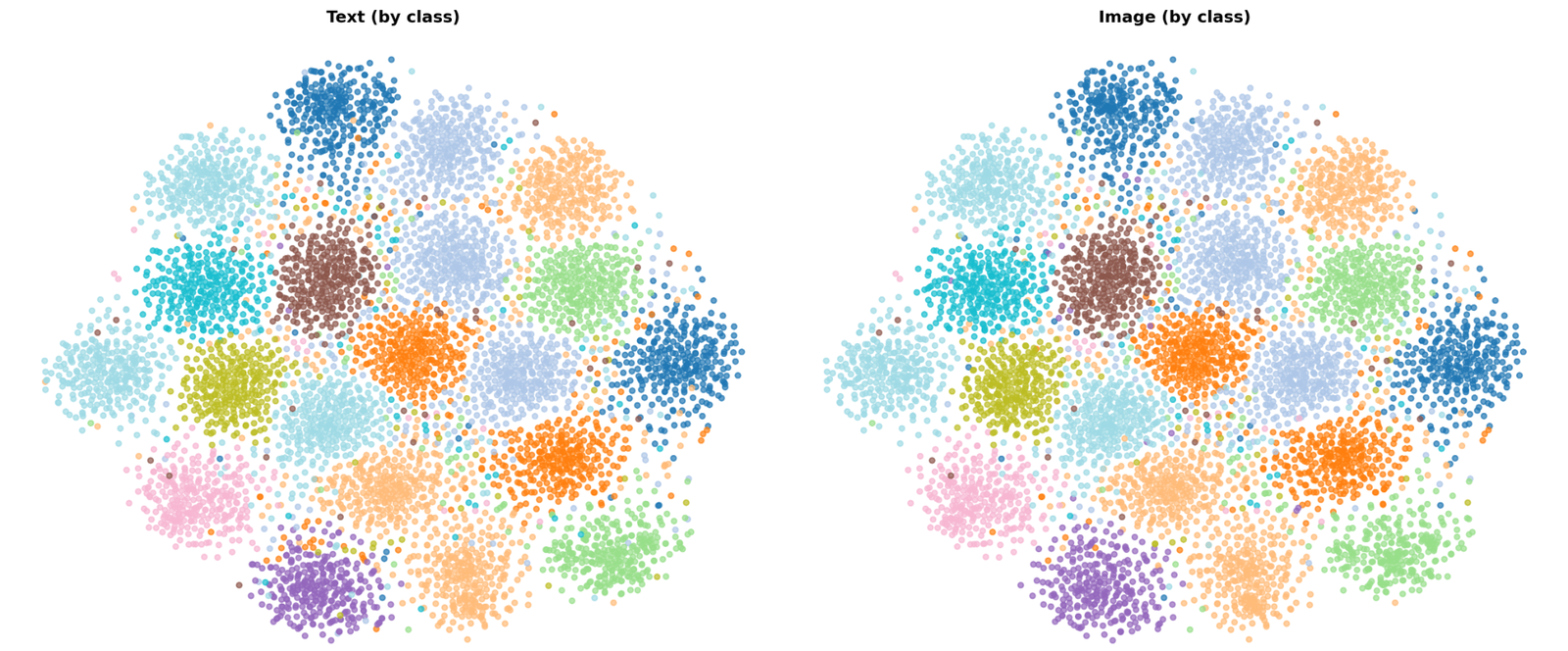}\\[2pt]
        \caption{JK-EGW(Ours)}
        \label{fig:tsne-ours}
    \end{subfigure}
    \hfill
    \begin{subfigure}[b]{0.46\linewidth}
        \centering
        \includegraphics[width=\linewidth]{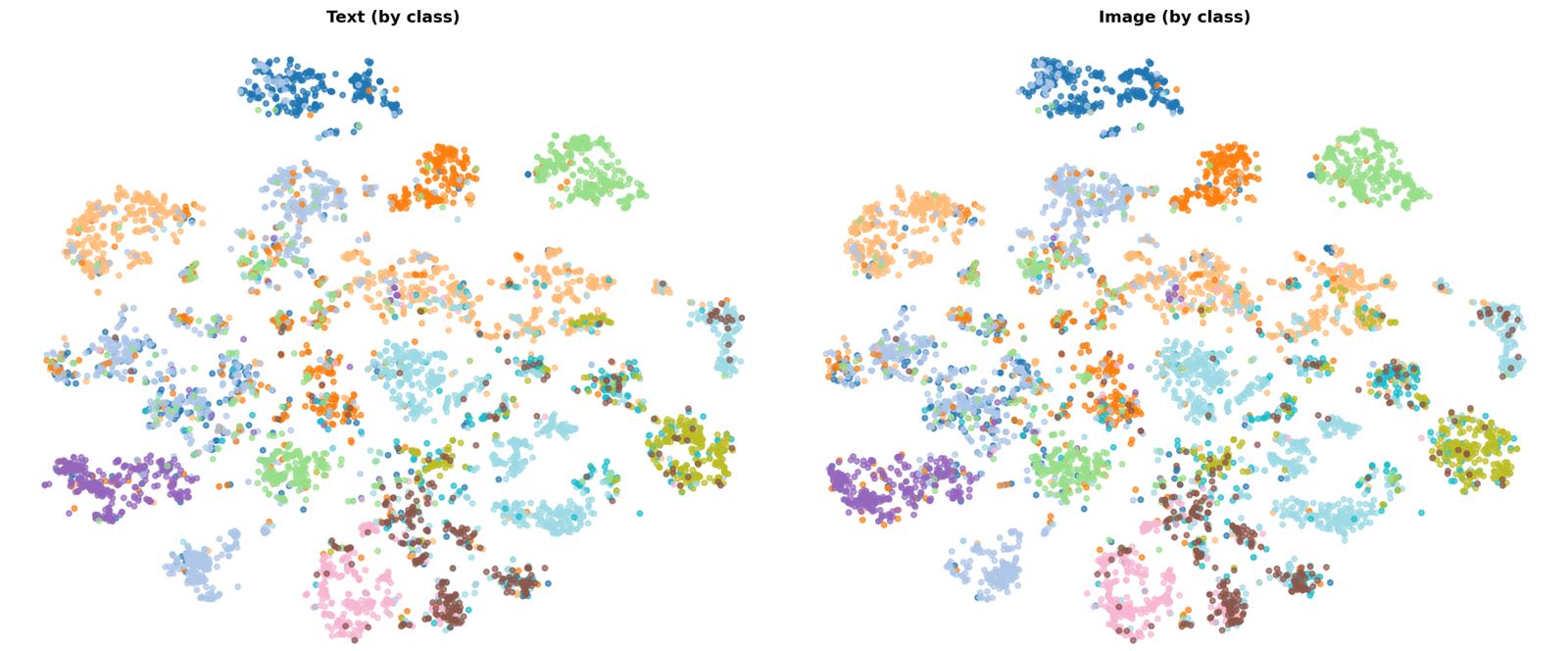}\\[2pt]
        \caption{CSA+$\textit{STRUCTURE}$}
        \label{fig:tsne-csars}
    \end{subfigure}

    \vspace{2mm}

    \begin{subfigure}[b]{0.46\linewidth}
        \centering
        \includegraphics[width=\linewidth]{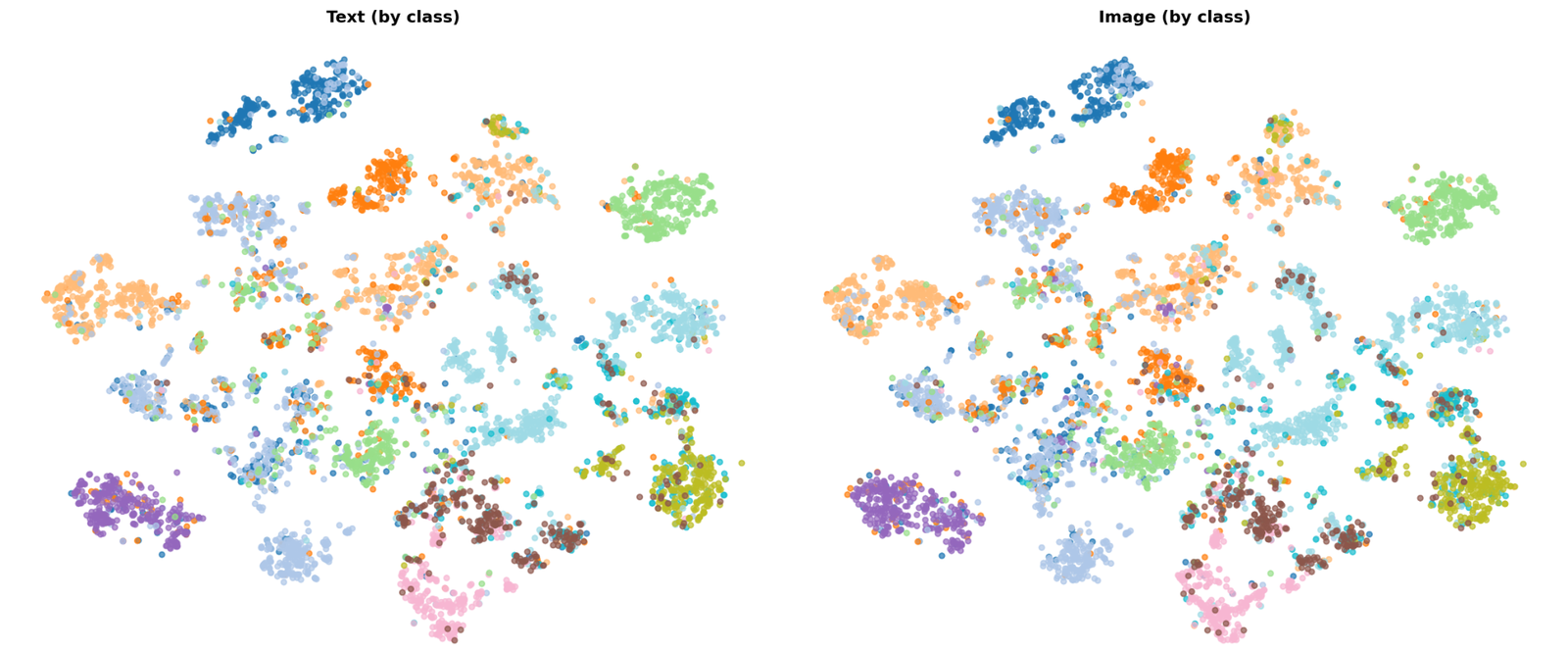}\\[2pt]
        \caption{CSA}
        \label{fig:tsne-csa}
    \end{subfigure}
    \hfill
    \begin{subfigure}[b]{0.46\linewidth}
        \centering
        \includegraphics[width=\linewidth]{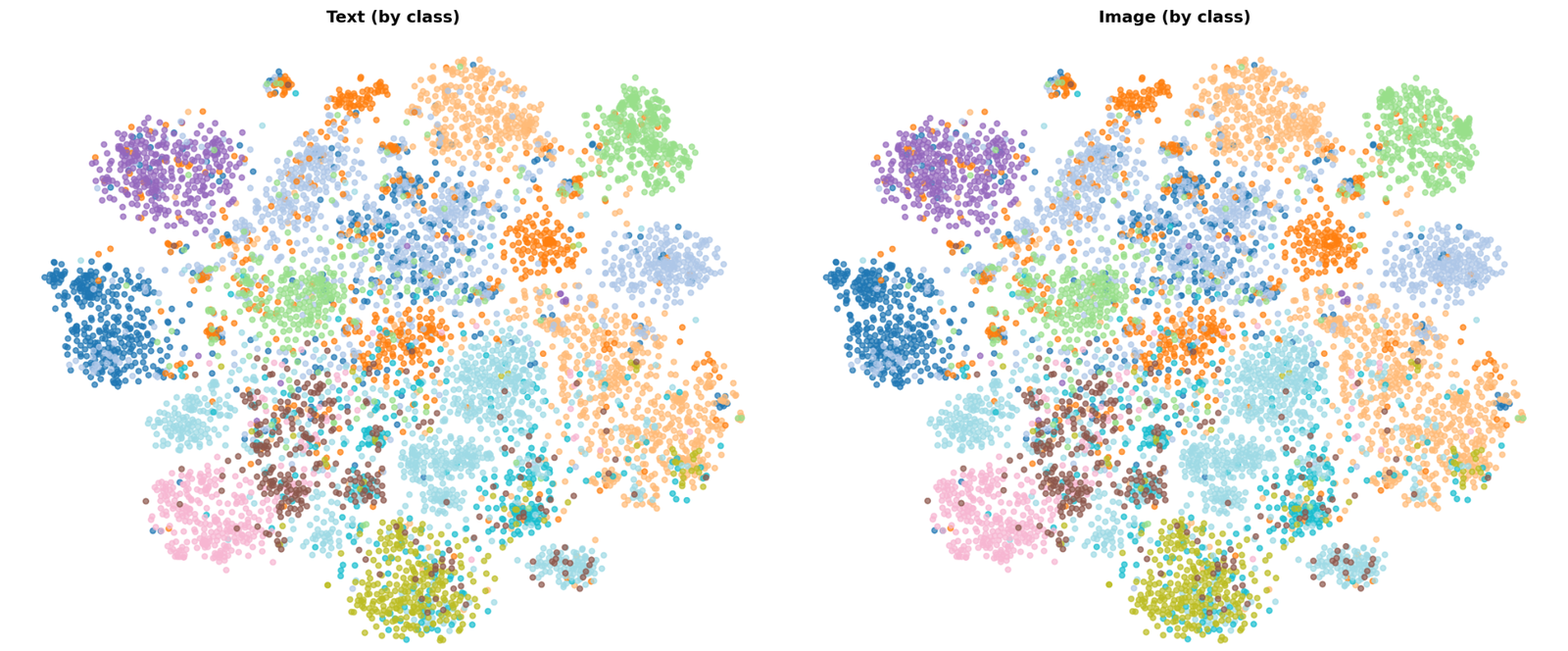}\\[2pt]
        \caption{MLP+$\textit{STRUCTURE}$}
        \label{fig:tsne-mlprs}
    \end{subfigure}

    \caption{Two-dimensional t-SNE embeddings of the shared representations learned by \name and the baselines on the MS--COCO 6K subset. 
    }
    \label{fig:ours-mscoco-6k}
\end{figure}

Figure~\ref{fig:ours-mscoco-6k} visualizes the shared text and image representations learned by \name and the baselines on the MS--COCO 6K subset. More specifically, we use the t-SNE~\citep{van2008visualizing} to generate 2D embedding coordinated of the shared representations. These plots visually support the claim that the distributions of the latent embeddings from different modalities are similar, while embeddings from different categories remain well separated. For \name, the embeddings of each class concentrate in a distinct contiguous region, with the text and image modalities visually showing the same arrangement of these regions. In contrast, the categories in the baseline representations are less compact.

\section{Conclusion}
We propose JK-EGW problem that maps several frozen unimodal encoders into a shared latent space with a prespecified  marginal distribution. Theoretically, when a finite-rank Mercer approximation of the affinity kernel exists, we can certify the statistical error of the JK-EGW problem, despite its nonconvexity. Algorithmically, we develop an alternating optimization method to solve JK-EGW. Empirically, \name learns geometrically coherent shared representations and strong bidirectional retrieval on MS--COCO.

Our study also points to several natural directions for future work. First, principled strategies for selecting or adapting the reference distribution $\nu$ to the downstream tasks would be valuable. Second, our setting focuses on finite number of modalities. Certifying the behavior of JK-EGW when $M$ is extremely large would be helpful in large-scale multi-modal alignment tasks.

\section*{Impact Statement}
This paper presents work whose goal is to advance the field of Machine Learning. There are many potential societal consequences of our work, none which we feel must be specifically highlighted here.




\bibliography{sample}

\newpage

\appendix
\section{Proofs for Section~\ref{sect:2}}
\subsection{Proof of Proposition~\ref{prop:kernel}}\label{proof:prop_kernel}
We first prove the claims for general kernels and specialize to the RBF kernel, then prove parts (a) and (b) by computing the Mercer eigensystem of a stationary kernel on $[0,1]^{d_x}$ explicitly.

\paragraph{General kernels on a compact domain.}
Under Proposition~\ref{prop:kernel}'s hypotheses, Mercer's theorem provides an orthonormal basis $\{\phi'_r\}_{r\ge1}$ of $L^2(\rho)$ consisting of eigenfunctions of the integral operator $T_Kf\coloneqq\int_{\mathcal X}K(\cdot,x')f(x')\,d\rho(x')$, with eigenvalues $\lambda_r\ge0$, such that
\begin{equation}\label{eq:mercer_expansion}
    K(x,x')=\sum_{r\ge1}\lambda_r\,\phi'_r(x)\,\phi'_r(x')
\end{equation}
converges absolutely and uniformly on $\mathcal X\times\mathcal X$, and $\phi'_r$ is continuous whenever $\lambda_r>0$~\citep{steinwart2012mercer}. Since $x\mapsto K(x,x)$ is continuous on the compact set $\mathcal X$, it attains its maximum, which gives the diagonal bound with $\kappa^2=\max_{x}K(x,x)<\infty$. The partial sums $K_R(x,x)=\sum_{r\le R}\lambda_r\phi'_r(x)^2$ are continuous and increase pointwise to the continuous limit $K(x,x)$, so by Dini's theorem the convergence is uniform on $\mathcal X$, that is,
\[
    \eta_R=\sup_{x\in\mathcal X}\bigl(K(x,x)-K_R(x,x)\bigr)\longrightarrow 0 .
\]
For the RBF kernel, and it is continuous, and positive semidefinite on every domain, because
\[
    \exp\!\Bigl(-\frac{\|x-x'\|^2}{2\ell^2}\Bigr)
    =e^{-\|x\|^2/(2\ell^2)}\,e^{-\|x'\|^2/(2\ell^2)}\sum_{j\ge0}\frac{1}{j!}\Bigl(\frac{x^\top x'}{\ell^2}\Bigr)^{j},
\]
where the linear kernel $x^\top x'$ is positive semidefinite, each power $(x^\top x')^j$ is then positive semidefinite by the Schur product theorem. Hence the RBF kernel satisfies Assumption~\ref{assumption2.1} on every compact domain, and $K(x,x)=1$ gives $\kappa^2=1$.

Now, let $\mathcal X=[0,1]^{d_x}$ with $\rho$ the uniform measure and $K(x,x')=k(x-x')$, with differences taken modulo $1$. Since $K$ is stationary, $K(x,x')$ depends only on $(x-x')\mod 1$, so $k$ extends to a $1$-periodic function on $\mathbb R^{d_x}$ through $k(t+w)=k(t)$ for $w\in\mathbb Z^{d_x}$. $k$ is even since $k(t)=K(x,x')=K(x',x)=k(-t)$ for such a pair. Two elementary facts about arithmetic modulo $1$ are used below. For $m\in\mathbb Z^{d_x}$ define the complex exponentials and the Fourier coefficients
\[
    \varphi_m(x)\coloneqq e^{2\pi i m^\top x},
    \qquad
    \lambda_m\coloneqq\hat k(m)\coloneqq\int_{[0,1]^{d_x}}k(t)\,e^{-2\pi i m^\top t}\,dt ,
\]
so the family $\{\lambda_m\}_{m\in\mathbb Z^{d_x}}$ in parts (a) and (b) is the family of Fourier coefficients of $k$. 

\paragraph{Proof of (a).}
We first want to verify that the predefined complex exponentials $\varphi_m$ are eigenfunctions of $T_K$, and the corresponding eigenvalues are the Fourier coefficients. Substituting $t=(x-x')\mod 1$ at fixed $x$, which is a negation and translation by $x$ for $x'$. Here, the negation and translation are bijections preserve the uniform measure, and
\[
    T_K\varphi_m(x)
    =\int_{[0,1]^{d_x}}k(x-x')\,e^{2\pi i m^\top x'}\,dx'
    =\Bigl(\int_{[0,1]^{d_x}}k(t)\,e^{-2\pi i m^\top t}\,dt\Bigr)e^{2\pi i m^\top x}
    =\hat k(m)\,\varphi_m(x),
\]
so every $\varphi_m$ is an eigenfunction of $T_K$ with eigenvalue $\hat k(m)$. These eigenvalues are nonnegative. Indeed, substituting the uniformly convergent expansion~\eqref{eq:mercer_expansion} and integrating termwise gives $\langle T_Kf,f\rangle=\sum_r\lambda_r|\langle f,\phi'_r\rangle|^2\ge0$ for every $f\in L^2(\rho)$, and in particular $\hat k(m)=\langle T_K\varphi_m,\varphi_m\rangle\ge0$. They are also even in $m$, since the measure preserving substitution $t\mapsto-t$ modulo $1$ and evenness of $k$ give $\hat k(-m)=\hat k(m)$.

Next, we want to verify the real trigonometric functions form a valid Mercer eigensystem. For $m\ne0$ the identities
\[
    \sqrt2\cos(2\pi m^\top x)=\tfrac{1}{\sqrt2}\bigl(\varphi_m(x)+\varphi_{-m}(x)\bigr),
    \qquad
    \sqrt2\sin(2\pi m^\top x)=\tfrac{1}{i\sqrt2}\bigl(\varphi_m(x)-\varphi_{-m}(x)\bigr)
\]
describe a unitary change of basis of $\mathrm{span}\{\varphi_m,\varphi_{-m}\}$. We therefore work in these normalized cosine and sine functions, together with the constant $1$. By linearity of the operator $T_K$ and $\hat k(m)=\hat k(-m)$, every element is an eigenfunction of $T_K$, and the constant $1$ has eigenvalue $\lambda_0$. The cosine and sine pair both have eigenvalue $\hat k(m)=\hat k(-m)$, so the eigenvalue multiset of this system is exactly $\{\lambda_m\}_{m\in\mathbb Z^{d_x}}$. This trigonometric system and the Mercer basis are two orthonormal eigenbases of the compact self adjoint operator $T_K$, so they carry the same multiset of nonzero eigenvalues, and
\[
    \sum_{m\in\mathbb Z^{d_x}}\lambda_m
    =\sum_{r\ge1}\lambda_r
    =\int_{[0,1]^{d_x}}K(x,x)\,dx
    =k(0)<\infty .
\]
In particular the non-increasing rearrangement $\{\lambda_r\}_{r\ge1}$ of $\{\lambda_m\}$ in the statement of part (a) is well defined. Summability also makes the series $\sum_m\lambda_m\,e^{2\pi i m^\top(x-x')}$ converge absolutely and uniformly, whose sum is a continuous $1$-periodic function with the same Fourier coefficients as $k$, hence $\sum_m\lambda_m\,e^{2\pi i m^\top(x-x')}=k(x-x')$. Splitting each exponential into real and imaginary parts, the imaginary contributions of $m$ and $-m$ cancel because $\lambda_m=\lambda_{-m}$, and the identity $\cos(a-b)=\cos a\cos b+\sin a\sin b$ with $a=2\pi m^\top x$ and $b=2\pi m^\top x'$ gives
\[
\begin{aligned}
    K(x,x')=\sum_{m\in\mathbb Z^{d_x}}\lambda_m\bigl[
    &\cos(2\pi m^\top x)\cos(2\pi m^\top x')\\
    &+\sin(2\pi m^\top x)\sin(2\pi m^\top x')\bigr].
\end{aligned}
\]
The summand is invariant under $m\mapsto-m$, so each pair $\{m,-m\}$ contributes twice its cosine product and twice its sine product, and each doubled product equals the product of the corresponding $\sqrt2$ normalized functions. With the constant term $\lambda_0$ included, the right hand side is therefore $\sum_{r\ge1}\lambda_r\,\phi'_r(x)\,\phi'_r(x')$ for the real system above.

Finally, since every element of the real system satisfies $\phi'_r(x)^2\le2$ pointwise, 
\[
    \eta_R
    =\sup_{x\in[0,1]^{d_x}}\sum_{r>R}\lambda_r\,\phi'_r(x)^2
    \le2\sum_{r>R}\lambda_r ,
\]
which is exactly claim (a).

\paragraph{Proof of (b).}
We first want to calculate the exact eigenvalues of the RBF kernel. Since differences are taken modulo $1$, the profile of the RBF kernel in part (b) is understood as the $1$-periodic function 
\[
    k(t)=\sum_{w\in\mathbb Z^{d_x}}\exp\Bigl(-\frac{\|t+w\|^2}{2\ell^2}\Bigr).
\]
The series converges uniformly by Gaussian decay, so $k$ is continuous, even, and $1$-periodic, and $K(x,x')=k(x-x')$ is a continuous symmetric stationary kernel. In particular, the eigenvalue family of part (b) is the family of Fourier coefficients $\lambda_m=\hat k(m)$ of this $k$, which we now compute in closed form. The Gaussian translates have finite total integral, so Fubini's theorem allows the series to be integrated termwise. In each term the substitution $s=t+w$ leaves the exponential unchanged, because $e^{-2\pi i m^\top t}=e^{-2\pi i m^\top s}\,e^{2\pi i m^\top w}$ and $e^{2\pi i m^\top w}=1$ for $w\in\mathbb Z^{d_x}$. The shifted cubes $w+[0,1]^{d_x}$ partition $\mathbb R^{d_x}$, so the sum of the shifted integrals collapses into a single integral over $\mathbb R^{d_x}$,
\[
\begin{aligned}
    \hat k(m)
    &=\sum_{w\in\mathbb Z^{d_x}}\int_{[0,1]^{d_x}}e^{-\|t+w\|^2/(2\ell^2)}\,e^{-2\pi i m^\top t}\,dt\\
    &=\sum_{w\in\mathbb Z^{d_x}}\int_{w+[0,1]^{d_x}}e^{-\|s\|^2/(2\ell^2)}\,e^{-2\pi i m^\top s}\,ds\\
    &=\int_{\mathbb R^{d_x}}e^{-\|s\|^2/(2\ell^2)}\,e^{-2\pi i m^\top s}\,ds .
\end{aligned}
\]
This last integral is the Fourier transform of the Gaussian at frequency $m$. It factorizes over the coordinates, and in one dimension completing the square gives
\[
    \int_{\mathbb R}e^{-s^2/(2\ell^2)}\,e^{-2\pi i\xi s}\,ds
    =e^{-2\pi^2\ell^2\xi^2}\int_{\mathbb R}e^{-(s+2\pi i\ell^2\xi)^2/(2\ell^2)}\,ds
    =\sqrt{2\pi}\,\ell\,e^{-2\pi^2\ell^2\xi^2},
\]
where the shifted Gaussian integral equals $\sqrt{2\pi}\,\ell$ after moving the contour back to the real axis by Cauchy's theorem. Hence
\begin{equation}\label{eq:rbf_eigenvalues_exact}
    \lambda_m=\hat k(m)=(2\pi\ell^2)^{d_x/2}\exp\bigl(-2\pi^2\ell^2\|m\|^2\bigr),
    \qquad m\in\mathbb Z^{d_x}.
\end{equation}

We finally move to calculate the specific decay rate of $\eta_R$ for the RBF kernel. Write $c\coloneqq2\pi^2\ell^2$, so $\lambda_m=C_0\,e^{-c\|m\|^2}$ by~\eqref{eq:rbf_eigenvalues_exact}, with $C_0=(2\pi\ell^2)^{d_x/2}$. For $T\ge1$, the number of lattice points $m\in\mathbb Z^{d_x}$ with $\|m\|\le T$ is at most $(2\lfloor T\rfloor+1)^{d_x}\le(3T)^{d_x}$, since such points satisfy $\|m\|_\infty\le T$ and each coordinate then takes at most $2\lfloor T\rfloor+1$ integer values. Set $c_{d_x}\coloneqq3^{d_x}$ and $T_R\coloneqq(R/c_{d_x})^{1/d_x}$. Let first $R\ge c_{d_x}$, so that $T_R\ge1$ and the ball of radius $T_R$ contains at most $R$ lattice points. Since $\sum_{r\le R}\lambda_r$ is the sum of the $R$ largest values of the family $\{\lambda_m\}$, it dominates the sum of $\lambda_m$ over any $R$ lattice points, in particular over a set of $R$ points containing this ball, so
\[
    \sum_{r>R}\lambda_r
    =\sum_{m\in\mathbb Z^{d_x}}\lambda_m-\sum_{r\le R}\lambda_r
    \le\sum_{\|m\|>T_R}\lambda_m .
\]
For $\|m\|>T_R$, splitting $e^{-c\|m\|^2}\le e^{-cT_R^2/2}\,e^{-c\|m\|^2/2}$ yields
\[
    \sum_{\|m\|>T_R}\lambda_m
    \le C_1\,e^{-cT_R^2/2},
    \qquad
    C_1\coloneqq C_0\Bigl(\sum_{j\in\mathbb Z}e^{-cj^2/2}\Bigr)^{d_x}<\infty ,
\]
where the lattice sum factorizes over coordinates and each factor is finite by comparison with a geometric series. Substituting $T_R^2=(R/c_{d_x})^{2/d_x}$ and invoking part (a),
\[
    \eta_R
    \le2\sum_{r>R}\lambda_r
    \le2C_1\exp\bigl(-c'R^{2/d_x}\bigr),
    \qquad
    c'\coloneqq\frac{c}{2}\,c_{d_x}^{-2/d_x}=\frac{\pi^2\ell^2}{9},
\]
for all $R\ge c_{d_x}$. For $1\le R<c_{d_x}$, part (a) gives $\eta_R\le2\sum_m\lambda_m=2k(0)$, while $\exp(-c'R^{2/d_x})\ge e^{-\pi^2\ell^2}$, so the same bound holds for all $R\ge1$ after enlarging the prefactor. This proves $\eta_R\lesssim e^{-c'R^{2/d_x}}$ with $c'$ depending only on $\ell$ and $d_x$, completing the proof.

\subsection{Proof of Proposition~\ref{prop:kernel_truncation_error}}
\label{proof:kernel_error}
The entropy terms are the same in $S_\epsilon$ and $S_\epsilon^R$, so
\[
\begin{aligned}
    |S_\epsilon(\{\mu^i\}_{i=1}^M,\nu)-S_\epsilon^R(\{\mu^i\}_{i=1}^M,\nu)|
    &\le \sup_{\substack{\pi^i\in\Gamma(\mu^i,\nu)\\ i\in[M]}}
    \left|
    \sum_{i,j=1}^M
    \iint
    \mathcal E_R(x,x')\|z-z'\|^2
    \,d\pi^i(x,z)d\pi^j(x',z')
    \right| .
\end{aligned}
\]
By Assumption~\ref{assumption2.1}, the residual kernel $\mathcal E_R=K-K_R$ is positive semidefinite and satisfies
\[
    |\mathcal E_R(x,x')|
    \le
    \sqrt{\mathcal E_R(x,x)\mathcal E_R(x',x')}
    \le
    \eta_R .
\]
Therefore,
\[
\begin{aligned}
    |S_\epsilon(\{\mu^i\}_{i=1}^M,\nu)-S_\epsilon^R(\{\mu^i\}_{i=1}^M,\nu)|
    &\le
    \eta_R
    \sup_{\substack{\pi^i\in\Gamma(\mu^i,\nu)\\ i\in[M]}}
    \sum_{i,j=1}^M
    \iint
    \|z-z'\|^2
    \,d\pi^i(x,z)d\pi^j(x',z') .
\end{aligned}
\]
Since each $\pi^i$ has $z$-marginal $\nu$, each term in the double sum equals
\[
    \iint_{\mathcal Z\times\mathcal Z}
    \|z-z'\|^2\,d\nu(z)d\nu(z').
\]
so the quantity inside the supremum does not depend on the choice of the couplings. Thus
\[
    |S_\epsilon(\{\mu^i\}_{i=1}^M,\nu)-S_\epsilon^R(\{\mu^i\}_{i=1}^M,\nu)|
    \le
    M^2\eta_R
    \iint
    \|z-z'\|^2\,d\nu(z)d\nu(z') .
\]
Let $Z,Z'\sim\nu$ be independent, then
\[
    \iint \|z-z'\|^2\,d\nu(z)d\nu(z')
    =
    \mathbb E\|Z-Z'\|^2
    \le
    2\mathbb E\|Z\|^2 .
\]
The $4$-sub-Weibull assumption implies

\[
    \mathbb E\|Z\|^2 \le \sqrt{2\log 2}\,\sigma 
\]
and therefore
\[
    \iint \|z-z'\|^2\,d\nu(z)d\nu(z')
    \lesssim
    \sigma,
\]
which gives
\[
    \left|
    S_\epsilon(\{\mu^i\}_{i=1}^M,\nu)
    -
    S_\epsilon^R(\{\mu^i\}_{i=1}^M,\nu)
    \right|
    \lesssim
    M^2\sigma\,\eta_R .
\]

\subsection{Proof of Theorem~\ref{thrm:pop_duality}} \label{proof:pop_duality}
We start from substituting rank-$R$ approximation of $K(x,x^\prime)$, $\Phi_R(x)^\top\Phi_R(x^\prime)$, into~\Cref{eq:main_objective} and obtain:
\begin{equation*}
    \begin{aligned}
        S^R_{\epsilon}(\{\mu^i\},\nu) = &\inf_{\{\pi_i \}} \left\{ \right.2 \sum_{i} \int \Phi_R(x)^\top \left(\sum_{j}\int \Phi_R(x^\prime)d\mu^j(x^\prime) \right) \norm{z}^2 d\pi^i(x,z) - 2\sum_{i,j}\operatorname{tr}{V_iV_j^\top}
        \\& +\epsilon \sum_{i=1}^M D_{KL}(\pi^i \parallel \mu^i \otimes \nu) \left. \right\}  \quad ,
    \end{aligned}
\end{equation*}
where $V_i = \int \Phi_R(x) z^\top d\pi^i(x,z)$. 
Using the identity
\begin{equation*}
    -2\norm{\sum_i V_i}_F^2 = \inf_{A} \left(2 \norm{A}^2_F - 4\sum_i \left<A, V_i \right> \right) 
\end{equation*}
and note that the optimum is achieved at $A^*=\sum_i V_i$, we have
\begin{equation*}
    S^R_{\epsilon}(\{\mu^i\},\nu) = \inf_{A} \left\{2\norm{A}^2_F +\sum_i \OT^i_{A,m_\Phi,\epsilon}(\mu^i,\nu) \right\}\quad ,
\end{equation*}
where the cost function $c_{A,m_\Phi}^i(x,z)=2\left(\Phi_R(x)^\top \left(\sum_{j}\int \Phi_R(x^\prime)d\mu^j(x^\prime) \right) \right) \norm{z}^2 - 4\left< A,\Phi_R(x)z^\top \right>$. Since $A^*=\sum_i V_i$, 
\begin{equation*}
    \norm{A}_F \leq M\kappa\left( \int \norm{z}^4 d\nu \right)^{1/4} \leq 2^{1/4}M\kappa\sqrt{\sigma} .
\end{equation*}
Therefore the infimum is achieved in $A_M = \{A:\norm{A}_F \leq 2^{1/4}M\kappa\sqrt{\sigma} \}$.

\newpage
\section{Proofs for Section~\ref{sec:theory}}
\subsection{Proof of Proposition~\ref{thrm:sample_complexity}}\label{proof:thrm:sample_complexity}
We will prove the two-sample case while the one-sample case is a byproduct.
First, we rescale \(z\rightarrow \epsilon^{-1/2}z\) so that
\[
    S^R_\epsilon(\{\mu^i\},\nu)
    =
    \epsilon S^R_1(\{\mu^i\},\nu^\epsilon),
\]
where \(\nu^\epsilon\) is \(4\)-sub-Weibull with parameter
\(\sigma^2/\epsilon^2\). 
Similarly,
\[
    S^R_\epsilon(\{\hat\mu_n^i\},\hat\nu_n)
    =
    \epsilon
    S^R_1(\{\hat\mu_n^i\},\hat\nu_n^\epsilon).
\]
We set \(\epsilon=1\) and recover the general
case \(\epsilon>0\) at the end.
Let \(L>0\) be a radius to be specified later, and define
\[
    \mathcal A_L
    \coloneqq
    \left\{
        A\in\mathbb R^{R\times d_z}:\|A\|_F\le L
    \right\}.
\]
Recall that the population EOT cost is
\[
    c_{A,m_\Phi}^i(x,z)
    =
    2\,\Phi_R(x)^\top m_\Phi\,\|z\|^2
    -
    4\left\langle A,\Phi_R(x)z^\top\right\rangle_F,
\]
where
\[
    m_\Phi
    =
    \sum_{j=1}^M
    \int_{\mathcal X^j}\Phi_R(x')\,d\mu^j(x').
\]
For the empirical objective, define
\[
    \hat m_\Phi
    =
    \sum_{j=1}^M
    \int_{\mathcal X^j}\Phi_R(x')\,d\hat\mu_n^j(x')
    =
    \frac1n\sum_{j=1}^M\sum_{k=1}^n\Phi_R(x_k^j),
\]
and
\[
    \hat c_{A,\hat m_\Phi}^i(x,z)
    =
    2\,\Phi_R(x)^\top\hat m_\Phi\,\|z\|^2
    -
    4\left\langle A,\Phi_R(x)z^\top\right\rangle_F .
\]
So similarly the empirical EOT problem with cost $\hat c_A^i(x,z)$ is
\[
    \OT^i_{A,\hat m_\Phi,\epsilon}(\hat{\mu}^i_n,\hat{\nu}_n)
    \coloneqq
    \inf_{\pi^i\in\Gamma(\hat{\mu}^i_n,\hat{\nu}_n)}
    \left\{
        \int \hat c_A^i(x,z)\,d\pi^i(x,z)
        +
        \epsilon
        D_{\mathrm{KL}}\!\left(
            \pi^i\,\middle\|\,\hat{\mu}^i_n \otimes \hat{\nu}_n
        \right)
    \right\}.
\]

Combined with Theorem~\ref{thrm:pop_duality} and the same lifting argument
applied to the empirical objective, we have
\begin{align}
    &
    \left|
    S^R_1(\{\mu^i\},\nu)
    -
    S^R_1(\{\hat\mu_n^i\},\hat\nu_n)
    \right|
    \nonumber \\
    &\le
    \sum_{i=1}^M
    \sup_{A\in\mathcal A_L}
    \left|
    \OT^i_{A,m_\Phi,1}(\mu^i,\nu)
    -
    \OT^i_{A,\hat m_\Phi,1}(\hat\mu_n^i,\hat\nu_n)
    \right|
    \nonumber \\
    &\le
    \sum_{i=1}^M
    \sup_{A\in\mathcal A_L}
    \left|
    \OT^i_{A,m_\Phi,1}(\mu^i,\nu)
    -
    \OT^i_{A,m_\Phi,1}(\hat\mu_n^i,\hat\nu_n)
    \right|
    \label{eq:ot_sum}
    \\
    &\quad+
    \sum_{i=1}^M
    \sup_{A\in\mathcal A_L}
    \left|
    \OT^i_{A,m_\Phi,1}(\hat\mu_n^i,\hat\nu_n)
    -
    \OT^i_{A,\hat m_\Phi,1}(\hat\mu_n^i,\hat\nu_n)
    \right|.
    \label{eq:cost_diff}
\end{align}

We first control the cost-difference term in~\Cref{eq:cost_diff}. For
fixed \(A\) and fixed \(i\), the two EOT problems have the same marginals
\((\hat\mu_n^i,\hat\nu_n)\) and the same KL reference measure
\(\hat\mu_n^i\otimes\hat\nu_n\). They differ only through the cost. Hence
\begin{align*}
   &
   \left|
   \OT^i_{A,m_\Phi,1}(\hat\mu_n^i,\hat\nu_n)
   -
   \OT^i_{A,\hat m_\Phi,1}(\hat\mu_n^i,\hat\nu_n)
   \right|
   \\
   &\le
   \sup_{\pi^i\in\Gamma(\hat\mu_n^i,\hat\nu_n)}
   \int
   \left|
   c_{A,m_\Phi}^i(x,z)-\hat c_{A,\hat m_\Phi}^i(x,z)
   \right|d\pi^i(x,z)
   \\
   &=
   \sup_{\pi^i\in\Gamma(\hat\mu_n^i,\hat\nu_n)}
   2\int
   \left|
   \Phi_R(x)^\top(m_\Phi-\hat m_\Phi)
   \right|
   \|z\|^2
   d\pi^i(x,z)
   \\
   &\le
   2\kappa\|m_\Phi-\hat m_\Phi\|
   \int\|z\|^2d\hat\nu_n(z),
\end{align*}
where the last inequality uses \(\|\Phi_R(x)\|\le\kappa\). Summing over
\(i\) gives
\begin{equation}
\label{eq:bound_cost_intermediate}
    \sum_{i=1}^M
    \sup_{A\in\mathcal A_L}
    \left|
    \OT^i_{A,m_\Phi,1}(\hat\mu_n^i,\hat\nu_n)
    -
    \OT^i_{A,\hat m_\Phi,1}(\hat\mu_n^i,\hat\nu_n)
    \right|
    \le
    2M\kappa\|m_\Phi-\hat m_\Phi\|
    \int\|z\|^2d\hat\nu_n(z).
\end{equation}

Next,
\begin{align}
    \mathbb E\|m_\Phi-\hat m_\Phi\|
    &=
    \mathbb E
    \left\|
    \sum_{j=1}^M
    \left[
        \frac1n\sum_{k=1}^n\Phi_R(x_k^j)
        -
        \int\Phi_R(x)d\mu^j(x)
    \right]
    \right\|
    \nonumber \\
    &\le
    \sum_{j=1}^M
    \mathbb E
    \left\|
        \frac1n\sum_{k=1}^n
        \left(
            \Phi_R(x_k^j)
            -
            \int\Phi_R(x)d\mu^j(x)
        \right)
    \right\|
    \nonumber \\
    &\le
    \frac{M\kappa}{\sqrt n},
    \label{eq:cost_bound_expectation}
\end{align}
where the first inequality follows the triangle inequality. For the second, fix \(j\) and set \(Y_k\coloneqq\Phi_R(x_k^j)-\int\Phi_R(x)d\mu^j(x)\). The vectors \(Y_1,\dots,Y_n\) are i.i.d.\ and mean zero, so the cross terms \(\mathbb E[Y_k^\top Y_l]\), \(k\ne l\), vanish and 
\[
    \mathbb E\left\|\frac1n\sum_{k=1}^nY_k\right\|^2
    =
    \frac{\mathbb E\|Y_1\|^2}{n}
    =
    \frac1n
    \left(
        \mathbb E\|\Phi_R(x_1^j)\|^2
        -
        \left\|\mathbb E\,\Phi_R(x_1^j)\right\|^2
    \right)
    \le
    \frac{\kappa^2}{n},
\]
where the last inequality uses \(\|\Phi_R(x)\|\le\kappa\). Jensen's inequality then gives \(\mathbb E \norm{Y_1} \leq \left(\mathbb E \norm{Y_1}^2\right)^{1/2}\), and summing over \(j\in[M]\) gives the display.


We also need to specify the radius \(L\). Define
\[
    \widetilde\sigma
    \coloneqq
    \sigma
    \vee
    \inf\left\{
        s>0:
        \int
        \exp\left(\frac{\|z\|^4}{2s^2}\right)
        d\hat\nu_n(z)
        \le 2
    \right\},
\]
then both \(\nu\) and \(\hat\nu_n\) are \(4\)-sub-Weibull with parameter
\(\widetilde\sigma^2\). We now choose
\[
    L
    =
    M\kappa(2\widetilde\sigma^2)^{1/4}.
\]
With this choice, the population outer infimum may be restricted to \(\mathcal A_L\), since \(\widetilde\sigma\ge\sigma\). Conditionally on the samples, the same argument as in the proof of Theorem~\ref{thrm:pop_duality} shows that the empirical outer optimizer also belongs to \(\mathcal A_L\). Indeed, if the latent marginal is $\nu$ or $\hat\nu_n$, then the optimizing matrix has the form
\[
    A^*=\sum_{i=1}^M V_i,
    \qquad
    V_i=\int \Phi_R(x)z^\top\,d\pi^i(x,z),
\]
and therefore
\[
    \|A^*\|_F
    \le
    \sum_{i=1}^M
    \int \|\Phi_R(x)\|\,\|z\|\,d\pi^i(x,z)
    \le
    M\kappa
    \left(\int\|z\|^2d\nu(z)\right)^{1/2}
    \le
    M\kappa(2\widetilde\sigma^2)^{1/4}
    =
    L .
\]
Thus the choice of \(\mathcal A_L\) above is valid for both the population and empirical objectives.

For each modality \(i\), denote $(\theta^i,\psi^i)$ to be the EOT potentials for $\OT^i_{A,m_\Phi,1}(\mu^i,\nu)$. For more details of EOT potentials, please refer to \citet{nutz2022entropic}. We now control the regularity of the EOT potentials for the frozen cost \(c_{A,m_\Phi}^i\), uniformly over \(A\in\mathcal A_L\). For notation simplicity, we omit the modality superscript from here for $\mu$, $\theta$, and $\psi$. More specifically, for any $\sigma$, let $\mathcal{F}_{\sigma}$ be the class of $\mathcal{C}^\infty(B_R(\kappa))$ functions $\theta(u)$, satisfying:
\begin{equation*}
\begin{aligned}
    & \theta(u) \leq 9\sqrt{2}M\kappa^2\sigma, \\
    & -\theta(u) \leq \log2 + M\kappa^2\left(3\sqrt{2}\sigma+ 4\cdot2^{1/4}\sigma^{1/2} \right) + 8M^2\kappa^4\sigma^2(1+2^{1/4}\sigma^{1/2})^2  \\
    & \left|D_u^\alpha \theta(u) \right| \leq C_\alpha (M\kappa)^{|\alpha|}\left(1+\sigma^{|\alpha|} \right) \left(1+\left(M\kappa^2(\sigma^{5/2}+\sigma^3)+M^2\kappa^4(\sigma^4+\sigma^{9/2}+\sigma^5)\right)^{|\alpha|/2}\right),
\end{aligned}
\end{equation*}
for every multi-index $\alpha$ with $1\le\alpha$, and some constant $C_\alpha>0$ that depends on $\alpha$. Similarly, for any $\sigma$, let $\mathcal{G}_\sigma$ be the class of $\mathcal{C}^\infty(\BR^{d_z})$ functions $\psi(z)$ satisfying:
\begin{equation*}
\begin{aligned}
    &\psi(z) \leq 2M\kappa^2 \norm{z}^2 + 4\cdot2^{1/4}M\kappa^2\sigma^{1/2} \norm{z} + 3\sqrt{2}M\kappa^2\sigma, \\
    & -\psi(z) \leq 2M\kappa^2\norm{z}^2 +4\cdot2^{1/4}M\kappa^2\sigma^{1/2}\norm{z}+5\sqrt{2}M\kappa^2 \sigma, \\
    & \left| D^\alpha_z\psi(z) \right| \leq  K_\alpha (M\kappa^2)^{|\alpha|}(1+\norm{z}+2^{1/4}\sigma^{1/2})^{|\alpha|},
\end{aligned}
\end{equation*}
for every multi-index $\alpha$ and some constant $K_\alpha>0$ that depends on $\alpha$. For the following analysis, we also introduce the induced $x$-space class:
\[
    \widetilde{\mathcal F}_\sigma^i
    =
    \left\{
        f:\mathcal X^i\to\mathbb R:
        f(x)=\theta(\Phi_R(x)),\
        \theta\in\mathcal F_\sigma
    \right\}.
\]
We slightly abuse the notation and use $\theta$ both for the function on $\BR^R$ (an element of $\mathcal{F}_\sigma$) and for its lift $\theta\circ\Phi_R$ on $\BR^{d_x}$ (an element of $\Tilde{\mathcal{F}}_\sigma$) when the meaning is clear from the context. 

The EOT dual problem for each corresponding pair of marginals is initially optimized over $L^1(\mu)\times L^1(\nu)$. The following Lemma~\ref{lemma:1} shows that, for each fixed \(A\in\mathcal A_L\), every fixed-cost population or empirical EOT problem appearing in~\Cref{eq:ot_sum} admits an optimal dual pair
\[
    (\theta\circ\Phi_R,\psi)
    \in
    \widetilde{\mathcal F}_{\widetilde\sigma}
    \times
    \mathcal G_{\widetilde\sigma}.
\]
Consequently, restricting the corresponding dual optimization to this product class does not change the EOT value.


\begin{lemma}[Uniform regularity of EOT potentials]
\label{lemma:1}
Fix a modality \(i\) and fix \(m_\Phi\) as defined above. Let \(\mu\in\mathcal P(\mathcal X^i)\). Suppose that \(\nu\) is \(4\)-sub-Weibull with parameter \(\sigma^2\), and fix \(A\) such that
\[
    \|A\|_F
    \le
    M\kappa(2\sigma^2)^{1/4}.
\]
Then there exist
\[
    \theta\in\mathcal F_\sigma,
    \qquad
    \psi\in\mathcal G_\sigma,
\]
such that
\[
    (\theta\circ\Phi_R,\psi)
\]
is a pair of optimal EOT potentials for \(\OT^i_{A,m_\Phi,1}(\mu,\nu)\).
\end{lemma}
The proof of Lemma~\ref{lemma:1} is in Section~\ref{sec:lemma1}.

Therefore, we define Hölder classes $\mathcal{F}_s$ and $\mathcal{G}_s$
\begin{align*}
    &\mathcal{F}_s = \{\theta:B_R(\kappa)\rightarrow\BR: |\theta|\leq C_{M,\kappa}, |D^\alpha\theta|\leq C_{M,\kappa,\alpha} \quad \forall|\alpha|\leq s  \}, \\
    & \mathcal{G}_s = \{\psi:\BR^{d_z}\rightarrow\BR: |\psi|\leq K_{M,\kappa}(1+\norm{\cdot}^2), |D^\alpha\psi|\leq K_{M,\kappa,\alpha}(1+\norm{\cdot}^s) \quad \forall|\alpha|\leq s  \}.
\end{align*}

By Lemma~\ref{lemma:1}, applied with \(\sigma=\widetilde\sigma\), for every \(A\in\mathcal A_L\) and every population or empirical EOT problem appearing in \Cref{eq:ot_sum}, there exist optimal potentials \((\theta,\psi)\) satisfying
\[
    (1+\widetilde\sigma^{4s})^{-1}\theta
    \in
    \mathcal F_s,
    \qquad
    (1+\widetilde\sigma^s)^{-1}\psi
    \in
    \mathcal G_s.
\]

The next step will be decomposing the EOT problem into suprema of empirical processes, where the above derived regularization for the potentials will be used to derive the parametric convergence rate for the empirical $S^R_1$. We go back the the RHS of Eq~\Cref{eq:ot_sum} and again fix a single modality. Since we work in the feature space instead of the raw $\mathcal{X}$ space, we need to define the induced distributions in the feature space. Suppose $U=\Phi_R(\mathcal X)$, $U_k=\Phi_R(x_k)$, then let $\mu_\Phi=\text{dist}(U)$ and $\hat \mu_{\Phi,n}=\frac{1}{n}\sum_{i=1}^n\delta_{U_k}$. Then,
\[
\left(\mu_\Phi - \hat \mu_{\Phi,n} \right)\theta = \int_{\mathcal{X}} \theta(\Phi_R(x))d(\mu-\hat{\mu}_n)(x).
\]
By the same dual-comparison argument as in Corollary~2 of
\citet{mena2019statistical}, uniformly over \(A\in\mathcal A_L\),
\begin{align}
\label{eq:ot_bound_1}
&
\sup_{A\in\mathcal A_L}
\left|
\OT^i_{A,m_\Phi,1}(\mu,\nu)
-
\OT^i_{A,m_\Phi,1}(\hat\mu_n,\hat\nu_n)
\right|
\\
&\qquad\lesssim
(1+\widetilde\sigma^{4s})
\sup_{\theta\in\mathcal F_s}
\left|
(\mu_\Phi-\hat\mu_{\Phi,n})\theta
\right|
+
(1+\widetilde\sigma^s)
\sup_{\psi\in\mathcal G_s}
\left|
(\nu-\hat\nu_n)\psi
\right|.
\nonumber
\end{align}

We will bound $\sup_{\theta\in \mathcal{F}_s}|(\mu_\Phi - \hat \mu_{\Phi,n})\theta|$ and $\sup_{\psi \in \mathcal{G}_s}|(\nu-\hat{\nu}_n)\psi|$ in expectation and handle $\Tilde{\sigma}$. From the Eq~\Cref{eq:ot_bound_1} we have
\begin{equation*}
\begin{aligned}
    & \BE\left[(1+\Tilde{\sigma}^{4s})\sup_{\theta\in \mathcal{F}_s}|(\mu_\Phi - \hat \mu_{\Phi,n}\theta| \right] \leq \sqrt{\BE\left[(1+\Tilde{\sigma}^{4s})^2\right] \BE\left[(\sup_{\theta\in \mathcal{F}_s}(\mu_\Phi - \hat \mu_{\Phi,n})\theta)^2 \right] }, \\
    & \BE\left[(1+\Tilde{\sigma}^{s})\sup_{\psi\in \mathcal{G}_s}|(\nu-\hat{\nu}_n)\psi| \right] \leq \sqrt{\BE\left[(1+\Tilde{\sigma}^{s})^2\right] \BE\left[(\sup_{\psi\in \mathcal{G}_s}(\nu-\hat{\nu}_n)\psi)^2 \right] }.
\end{aligned}
\end{equation*}
By Theorem 3.5.1 of~\citet{gine2021mathematical}, we have 
\begin{equation}\label{eq:fisrt_bound}
    \BE\left[\left(\sup_{\theta\in \mathcal{F}_s}(\mu_\Phi-\hat\mu_{\Phi,n})\theta\right)^2 \right] \lesssim \frac{1}{n} \BE \left(\int_0^{\sqrt{\sup_{\theta\in \mathcal{F}_s}\norm{\theta}^2_{L^2(\hat \mu_{\Phi,n})}}} \sqrt{\log 2N(\mathcal{F}_s,L^2(\hat \mu_{\Phi,n}),\delta)}d\delta \right)^2,
\end{equation}
where $N$ is the $\delta$-covering number of the class $\mathcal{F}_s$ and observe that 
\begin{equation*}
\begin{aligned}
    \sup_{\theta\in \mathcal{F}_s}\norm{\theta}^2_{L^2(\hat{\mu}_n)} \leq C^2_{M,\kappa}. 
\end{aligned}
\end{equation*}

Next, since $\norm{\Phi_R(x)}\leq\kappa$, both $\mu_\Phi$ and $\hat\mu_{\Phi,n}$ are supported on $B_R(\kappa)$. Let $\mathcal{X}^1=B_R(\kappa+1)$, then, by Theorem 2.7.1 in \citet{van1996weak}, with $s_1=\lceil R/2 \rceil+1$,
\begin{align*}
    \log N(\mathcal{F}_{s_1},L^2(\hat{\mu}_n),\delta) \leq \log N(\mathcal{F}_{s_1},\norm{\cdot}_{\infty,B_R(\kappa)},\delta)\lesssim_{R} (\kappa+1)^R\delta^{-R/s_1}. 
\end{align*}
Substituting back into Eq~\Cref{eq:fisrt_bound}, we have
\begin{align*}
    \BE\left[\left(\sup_{\theta\in \mathcal{F}_{s_1}}(\mu-\hat{\mu}_n)f\right)^2 \right] &\lesssim_{R} \frac{1}{n}\,\left((\kappa+1)^{R/2}\int_0^{C_{M,\kappa}}\delta^{-R/2s_1} d\delta \right)^2\\
    & \leq \frac{1}{n}(\kappa+1)^R\frac{C_{M,\kappa}^{2-R/s_1}}{\left(1-\frac{R}{2s_1} \right)^2}. 
\end{align*}

For $\BE\left[(\sup_{\psi\in \mathcal{G}_s}(\nu-\hat{\nu}_n)\psi)^2 \right]$, similarly we have
\begin{equation}\label{eq:second_bound}
    \BE\left[\left(\sup_{\psi\in \mathcal{G}_s}(\nu-\hat{\nu}_n)\psi\right)^2 \right] \lesssim \frac{1}{n} \BE \left(\int_0^{\sqrt{\sup_{\psi\in \mathcal{G}_s}\norm{\psi}^2_{L^2(\hat{\nu}_n)}}} \sqrt{\log 2N(\mathcal{G}_s,L^2(\hat{\nu}_n),\delta)}d\delta \right)^2.
\end{equation}
This case is similar to Theorem 2 of \citet{zhang2024gromov}. First observe that
\begin{equation*}
    \sup_{\psi \in \mathcal{G}_s} \norm{\psi}^2_{L^2(\hat{\nu}_n)} \leq K_{R,d_z} \frac{1}{n}\sum_{l=1}^n (1+\norm{z_l}^4)\leq K_{R,d_z}(1+\Tilde{\sigma}^2L_z)
\end{equation*}
where $L_z\coloneqq\frac{1}{n}\sum_{l=1}^n\exp(\frac{\norm{z_l}^4}{2\Tilde{\sigma}^2})$ satisfies $L_z\leq 2$ almost surely. To control the covering number, we need to partition our unbounded domain. More specifically, define $Q_q \coloneqq [-2^q\sqrt{\Tilde{\sigma}},2^q\sqrt{\Tilde{\sigma}}]^{d_z}$ for $q\in \mathbb{N}_0$ and partition $\BR^{d_z}$ into $I_q=Q_q\backslash Q_{q-1}$. To verify the conditions of Corollary 2.7.4 of Theorem~\citet{van1996weak}, note that the Lebesgue measure of each $\{ z\in \BR^{d_z}:\norm{z-I_q}\leq 1 \}$ is bounded by $K_{d_z}(1+2^{qd_z}\Tilde{\sigma}^{d_z/2})$. Moreover, $\hat{\nu}_n(I_q)=\frac{1}{n}\sum_{l=1}^n\vecone_{\{z_l\in I_q\}}\leq L_z\exp(-2^{4q-5})$. Finally, for any $q\in \mathbb{N}_0$ and $\psi \in \mathcal{G}_s$, the restriction $\psi|_{I_q}$ has a $\mathcal{C}^s(I_q)$-Hölder norm bounded by $K_{s,d_z}(1+\Tilde{\sigma}^{s/2})2^{qs}$. With $s_2=\lceil d_z/2 \rceil+1$, $V=d_z/s_2$ and $r=2$, we have
\begin{equation*}
\begin{aligned}
    &\log N(\delta,\mathcal{G}_{s_2},L^2(\hat{\nu}_n)) \\ \leq& \log N_{[\,]}(2\delta,\mathcal{G}_{s_2},L^2(\hat{\nu}_n)) \\
    \leq& K_{R,d_z}\delta^{-V}L_z^{V/r}\left( \sum_{q=0}^\infty(1+2^{qd_z}\Tilde{\sigma}^{d_z/2})^{r/(V+r)} \exp(-2^{4q-5})^{V/(V+r)} \left((1+\Tilde{\sigma}^{s_2/2})2^{qs_2}\right)^{Vr/(V+r)}\right)^{(V+r)/r} \\
    \leq& K_{R,d_z}\delta^{-V}L_z^{V/r} \left(1+\Tilde{\sigma}^{\frac{d_z+s_2V}{2}} \right) \left(\sum_{l=0}^\infty \exp\left(\frac{(-2^{4q-5})V}{V+r}\right)2^{\frac{(s_2V+d_z)qr}{V+r}} \right)^{(V+r)/r} \\
    \leq & K_{R,d_z}\delta^{-\frac{d_z}{s_2}}L_z^{\frac{d_z}{2s_2}}(1+\Tilde{\sigma}^{d_z}). 
\end{aligned}
\end{equation*}
Substituting back to Eq~\Cref{eq:second_bound}, we have
\begin{equation*}
\begin{aligned}
    \BE\left[(\sup_{\psi\in \mathcal{G}_s}(\nu-\hat{\nu}_n)\psi)^2 \right] &\lesssim_{R,d_z} \frac{1}{n} \BE \left[ \left( \int_0^{\sqrt{1+\Tilde{\sigma}^2L_z}}\sqrt{\delta^{\frac{d_z}{s_2}}L_z^{\frac{d_z}{2s_2}}(1+\Tilde{\sigma}^{d_z})}d\delta \right)^2 \right] \\
    &\lesssim_{R,d_z} \frac{1+\Tilde{\sigma}^{d_z}}{n}\BE\left[ \left(L_z^{d_z/4s_2}\int_0^{\sqrt{1+\Tilde{\sigma}^2L_z}} \delta^{-d_z/2s_2} d\delta\right)^2 \right] \\
    &\lesssim_{R,d_z} \frac{1}{n}(1+\Tilde{\sigma}^{d_z+2}). 
\end{aligned}
\end{equation*}

Finally, we take $s=s_1\vee s_2$. It remains to bound moments of $\Tilde{\sigma}$. Fix $k \in \mathbb{N}$, define
\begin{equation*}
    \tau^2_k = \sigma^2 \vee \left( \frac{k\sigma^2}{n} \sum_{l=1}^n\exp(\norm{z_l}^4/(2k\sigma^2)) \right)
\end{equation*}
so that $\nu,\hat{\nu}_n$ are both 4-sub-Weibull with parameter $\tau_k^2$. By Lemma 4 of \citet{mena2019statistical} and similar to Theorem 2 of \citet{zhang2024gromov}, we have
\begin{equation*}
    \BE[\Tilde{\sigma}^{2k}] \lesssim_k\sigma^{2k}.
\end{equation*} 
Hence, set $q_{R,z}=\left\lceil\frac{R\vee d_z}{2}\right\rceil$, we have
\begin{equation*}
\begin{aligned}
    & \BE\left[(1+\Tilde{\sigma}^{4s})\sup_{\theta\in \mathcal{F}_s}|(\mu-\hat{\mu}_n)\theta| \right] \leq \sqrt{\BE[(1+\Tilde{\sigma}^{4s})^2]\frac{(\kappa+1)^R}{n\left( 1-\frac{R}{2s} \right)^2}} \lesssim_{M,R,d_z,\kappa} \frac{1+\sigma^{4q_{R,z}+4}}{\sqrt{n}}, \\
    & \BE\left[(1+\Tilde{\sigma}^{s})\sup_{\psi\in \mathcal{G}_s}|(\nu-\hat{\nu}_n)\psi| \right] \leq \sqrt{\BE[(1+\Tilde{\sigma}^{s})^2]\frac{(1+\sigma^{d_z+2})}{n} } \lesssim_{M,R,d_z,\kappa} \frac{1+\sigma^{2q_{R,z}+2}}{\sqrt{n}}. 
\end{aligned}
\end{equation*}
Combine with the bounds in~\Cref{eq:cost_diff}, we obtain the final results.

\newpage
\section{Proofs for Section~\ref{sec:optimization}}
\subsection{Proof of Theorem~\ref{thrm:sample_duality}}\label{proof:sample_duality}
The objective in~\Cref{eq:sum_obj} can be rewritten in the following matrix form. 
\begin{equation*}
\begin{aligned}
    \min_{\hat{\Pi}} & \left\langle\hat{\Pi}^{\top}  W\hat{\Pi}, D( \hat Z)\right\rangle = \min_{\hat{\Pi}} \operatorname{tr}{\hat{\Pi}^{\top}  W\hat{\Pi} D( \hat Z)} \\
    \text{subject to} &\quad \pi^i\vecone = \frac{1}{n}\vecone, \text{and}\quad \pi^i{}^{\top}\vecone = \frac{1}{n}\vecone \qquad \forall i\in [M]
\end{aligned} \quad,
\end{equation*}
where $D( \hat Z)=diag( \hat Z \hat Z^\top)\vecone^\top + \vecone diag( \hat Z \hat Z^\top)^\top - 2 \hat Z \hat Z^\top$. Therefore, we have:
\begin{equation*}
    \begin{aligned}
        \hat S^{\hat R}_{\epsilon}(\{\hat{\mu}^i_n\},\hat{\nu_n})= \min_{\hat{\Pi}} \operatorname{tr}{\hat{\Pi}^{\top} LL^\top\hat{\Pi} D( \hat Z)} +\epsilon \sum_{i=1}^M D_{KL}(\hat{\pi}^i \parallel \hat{\mu}^i \otimes \hat{\nu}).
    \end{aligned}
\end{equation*}

Using the identity
\begin{equation*}
    -\norm{M}^2 = \inf_A \left(\norm{A}^2 - 2\innerproduct{A,M} \right),
\end{equation*} 

we have

\begin{equation*}
    \begin{aligned}
        \hat S^{\hat R}_{\epsilon}(\{\hat{\mu}^i_n\},\hat{\nu_n})=&\min_{\hat{\Pi}} \quad \operatorname{tr}{\hat{\Pi}^\top LL^\top \hat{\Pi} diag( \hat Z \hat Z^\top)\vecone^\top} +\operatorname{tr}{\hat{\Pi}^\top  LL^\top \hat{\Pi} \vecone diag( \hat Z \hat Z^\top)^\top} - \\&2\norm{L^\top \hat{\Pi}  \hat Z}^2_F +\epsilon \sum_{i=1}^M D_{KL}(\hat{\pi}^i \parallel \hat{\mu}^i \otimes \hat{\nu}) \\
         =& \min_{\hat{\Pi}}\quad  2 \left<\frac{1}{n}L L^\top \vecone  \hat Z_{\mathrm{norm}}^{\top}, \hat{\Pi} \right> -2\norm{L^\top \hat{\Pi}  \hat Z}^2_F +\epsilon \sum_{i=1}^M D_{KL}(\hat{\pi}^i \parallel \hat{\mu}^i \otimes \hat{\nu}) \\
          =& \min_{\hat{A}}  2\norm{\hat{A}}^2_F + \sum_{i=1}^M \OT^i_{C_i,\epsilon}(\hat{\mu}^i,\hat{\nu}). 
    \end{aligned} 
\end{equation*}

\subsection{Proof of Proposition~\ref{prop:empirical_mercer_realization}}
\label{pf:prop:empirical_mercer_realization}
\begin{align*}
    \BE \left[ \left|
    S_\epsilon^R(\{\hat\mu_n^i\}_{i=1}^M,\hat\nu_n)-
    \hat S_\epsilon(\{\hat\mu_n^i\}_{i=1}^M,\hat\nu_n)
    \right|\right] =& \BE \left[\left| \sum_{i,j} \sum_{k,k',\ell,\ell'} \left( K_R(x_k^i,x^j_{k'})-K(x_k^i,x^j_{k'})\right)\hat{\pi}^i_{k,\ell}\,\hat{\pi}^j_{k',\ell'} \,\norm{z_{\ell}-z_{\ell'}}^2 \right|\right] \\ 
    \le& \eta_R \cdot \BE \left[\left| \sum_{i,j} \sum_{k,k',\ell,\ell'}\hat{\pi}^i_{k,\ell}\,\hat{\pi}^j_{k',\ell'} \,\norm{z_{\ell}-z_{\ell'}}^2 \right|\right] 
\end{align*}
by Assumption~\ref{assumption2.1}. Now, 
\begin{align*}
     \sum_{k,k'}\hat{\pi}^i_{k,\ell}\,\hat{\pi}^j_{k',\ell'} = \left( \sum_k \hat{\pi}^i_{k,\ell} \right)\left( \sum_{k'}\hat{\pi}^j_{k',\ell'} \right) = \frac{1}{n}\cdot\frac{1}{n} . 
\end{align*}
In addition, $\nu$ is $4$-sub-Weibull with parameter $\sigma^2$, so 
\begin{align*}
    \BE \left[ \left|
    S_\epsilon^R(\{\hat\mu_n^i\}_{i=1}^M,\hat\nu_n)-
    \hat S_\epsilon(\{\hat\mu_n^i\}_{i=1}^M,\hat\nu_n)
    \right|\right] \leq M^2\eta_R\sigma. 
\end{align*}

\subsection{Proof of Proposition~\ref{prop:pivoted_cholesky_empirical_error}}
\label{pf:pivoted_cholesky_empirical_error}
We follow the proof in Section~\ref{pf:prop:empirical_mercer_realization}, replacing $K_R(x_k^i,x^j_{k'})-K(x_k^i,x^j_{k'})$ by $W_{(i,k),(j,k')}-(LL^\top)_{(i,k),(j,k')}$, which is bounded in absolute value by $\rho_{\hat R}$ uniformly over the index pairs. Applied conditionally on the samples, the same argument gives
\[
    \left|
    \hat S_\epsilon(\{\hat{\mu}_n^i\}_{i=1}^M,\hat{\nu}_n)
    -
    \hat S_\epsilon^{\hat R}(\{\hat{\mu}_n^i\}_{i=1}^M,\hat{\nu}_n)
    \right|
    \le
    M^2\rho_{\hat R}
    \iint_{\mathcal Z\times\mathcal Z}
    \|z-z'\|^2\,d\hat{\nu}_n(z)d\hat{\nu}_n(z').
\]
The residual $\rho_{\hat R}$ is a function of the modality samples $\{x_k^i\}$, while $\hat{\nu}_n$ is a function of the latent points $\{z_\ell\}$, and the two collections are drawn independently. The expectation of the product therefore factorizes,
\[
    \mathbb E
    \left[
    \rho_{\hat R}
    \iint
    \|z-z'\|^2\,d\hat{\nu}_n(z)d\hat{\nu}_n(z')
    \right]
    =
    \mathbb E[\rho_{\hat R}]\cdot
    \mathbb E
    \left[
    \iint
    \|z-z'\|^2\,d\hat{\nu}_n(z)d\hat{\nu}_n(z')
    \right].
\]

\newpage
\section{Proof of the Auxiliary Lemma}
\subsection{Proof of Lemma~\ref{lemma:1}}\label{sec:lemma1}
Fix $A\in \mathcal A_M$, where $\mathcal{A}_M=\{A:\norm{A}_F \leq M\kappa(2\sigma^2)^{1/4} \}$, assume the existence of optimal potentials $(\theta, \psi) \in L^1(\mu)\times L^1(\nu)$. Let $(\theta^i_0, \psi_0) \in L^1(\mu^i)\times L^1(\nu)$ be optimal EOT potentials for cost $c_A^i$ and assume WLOG normalization, $\int \theta^i_0 d\mu^i = \int \psi_0 d\nu(z)= \frac{1}{2}\OT^i_{A,1}(\mu^i,\nu)$. For a more detailed discussion of EOT potentials, please refer to \citet{nutz2022entropic}. Define new functions $\theta^i$ and $\psi$ as 
\begin{align*}
    &\theta^i(x) \coloneqq -\log \int e^{\psi_0(z)-c_A^i(x,z)} d\nu(z), \\
    &\psi(z) \coloneqq -\log \int e^{\theta^i_0(x)-c_A^i(x,z)} d\mu^i(x). 
\end{align*}
By Jensen's inequality, we have
\begin{align*}
    \theta^i(x) \leq \int c_A^i(x,z) d\nu(z) - \frac{1}{2}\OT^i_{A,1}(\mu^i,\nu). 
\end{align*}
Next, by the definition, 
\begin{align*}
    \int c_A^i(x,z) d\nu(z) &= 2\left(\Phi_R(x)^\top \left(\sum_{j}\int \Phi_R(x^\prime)d\mu^j(x^\prime) \right) \right) \int \norm{z}^2 d\nu(z) -4\left< A,\Phi_R(x)\left(\int z d\nu(z) \right)^\top \right> \\
    & \leq 2\left(\Phi_R(x)^\top \left(\sum_{j}\int \Phi_R(x^\prime)d\mu^j(x^\prime) \right) \right) \left(\sqrt{2}\sigma \right) + 4M\kappa^2(2\sigma^2)^{1/4}\left(2^{1/4}\sqrt{\sigma} \right) \\
    & \leq 6\sqrt{2}M\kappa^2 \sigma, 
\end{align*}
with the assumption $\|\Phi_R(x)\| \leq \kappa$.  

We can upper bound $-\OT^i_{A,1}(\mu^i, \nu)$ very similarly and obtain $-\OT^i_{A,1}(\mu^i, \nu) \leq 6\sqrt{2}M\kappa^2\sigma$. Together we have, 
\begin{align*}
    \theta^i(x) \leq 9\sqrt{2}M\kappa^2\sigma. 
\end{align*}

Similarly, we derive the upper bound for $\psi(z)$:
\begin{align*}
    \int c_A^i(x,z) d\mu^i(x) &= 2 \left< \int\Phi_R(x) d\mu^i(x),\sum_{j}\int \Phi_R(x^\prime)d\mu^j(x^\prime) \right> \norm{z}^2 - 4\left<A, \int\Phi_R(x) d\mu^i(x)z^\top \right> \\
    & \leq 2 \norm{\int\Phi_R(x) d\mu^i(x)} \norm{\sum_{j}\int \Phi_R(x^\prime)d\mu^j(x^\prime) }\norm{z}^2 + 4 \norm{A}_F \norm{\int\Phi_R(x) d\mu^i(x)}  \norm{z} \\
    & \leq 2M\kappa^2 \norm{z}^2 + 4\cdot2^{1/4}M\kappa^2\sigma^{1/2} \norm{z}.
\end{align*}
Therefore, 
\begin{align}\label{eq:upper_psi}
    \psi(z) \leq 2M\kappa^2 \norm{z}^2 + 4\cdot2^{1/4}M\kappa^2\sigma^{1/2} \norm{z} + 3\sqrt{2}M\kappa^2\sigma. 
\end{align}

For the lower bound, consider
\begin{align*}
    -\theta^i(x) &\leq \log \int e^{2M\kappa^2 \norm{z}^2 + 4\cdot2^{1/4}M\kappa^2\sigma^{1/2} \norm{z} + 3\sqrt{2}M\kappa^2\sigma+4\left<A,\Phi_R(x)z^\top \right> - 2\left(\Phi_R(x)^\top \left(\sum_{j}\int \Phi_R(x^\prime)d\mu^j(x^\prime) \right) \right) \norm{z}^2} d\nu(z) \\
    & \leq 3\sqrt{2}M\kappa^2\sigma + \log \int e^{4M\kappa^2\norm{z}^2 + 8\cdot 2^{1/4}M\kappa^2\sigma^{1/2}\norm{z}} d\nu(z) \\
    & \leq M\kappa^2\left(3\sqrt{2}\sigma+ 4\cdot2^{1/4}\sigma^{1/2} \right) + \log\int e^{4M\kappa^2\norm{z}^2+4\cdot 2^{1/4}M\kappa^2\sigma^{1/2}\norm{z}^2} d\nu(z) \\
    & \leq \log2 + M\kappa^2\left(3\sqrt{2}\sigma+ 4\cdot2^{1/4}\sigma^{1/2} \right) + 8M^2\kappa^4\sigma^2(1+2^{1/4}\sigma^{1/2})^2, 
\end{align*}
where we have used $2\norm{z}\leq 1+\norm{z}^2$ inside the exponents and $\int e^{t\norm{z}^2}d\nu(z) \leq 2e^{t^2\sigma^2/2}$ for a $\sigma^2$-sub-Gaussian random variable $\norm{z}^2$.
Similarly, 
\begin{align*}
    -\psi(z) & \leq \log \int e^{5\sqrt{2}M\kappa^2 \sigma+4\left<A,\Phi_R(x)z^\top \right> - 2\left(\Phi_R(x)^\top \left(\sum_{j}\int \Phi_R(x^\prime)d\mu^j(x^\prime) \right) \right) \norm{z}^2} d\mu^i(x) \\
    & \leq 2M\kappa^2\norm{z}^2 +4\cdot2^{1/4}M\kappa^2\sigma^{1/2}\norm{z}+5\sqrt{2}M\kappa^2 \sigma.
\end{align*}

Next, we want to bound the derivatives of $\theta^i(x)$. For a multi-index $\alpha,\beta \in \mathbb{N}_0^R$ where $\mathbb{N}_0 = \mathbb{N}\cup \{0\}$, let $D^\alpha$ denote the operator $\frac{\partial^{|\alpha|}}{\partial^{\alpha_1} u_1,\dots,\partial^{\alpha_R} u_R}$, where $|\alpha|=\sum_{a=1}^R \alpha_a$ and $D^0f=f$. 
Since $c_A^i$ only depends $x$ through $\Phi_R(x)$, write $u=\Phi_R(x) \in \BR^R$ and write $\Tilde{c}_A^i(u,z)=2u^\top\left(\sum_j \int\Phi_R(x^\prime)d\mu^j(x^\prime) \right)\norm{z}^2-4\left<A,uz^\top \right>$. 

Applying the multivariate Faa di Bruno formula \citep{constantine1996multivariate} to $D_u^\alpha \log \int e^{\psi_0(z)-\Tilde{c}_A^i(u,z)}d\nu(z)$ yields a finite sum, indexed by partitions of $\alpha$, of products of terms $\frac{D_u^\beta\int e^{\psi_0(z)-\tilde c_A^i(u,z)}d\nu(z)}{\int e^{\psi_0(z)-\tilde c_A^i(u,z)}d\nu(z)}, \quad |\beta|\leq |\alpha|$, with the number of terms depending only on $|\alpha|$. Thus, for us to establish the magnitude of derivatives bounds, it suffices to bound $\left|D_u^\beta\int e^{\psi_0(z)-\Tilde{c}_A^i(u,z)}d\nu(z)\right|$. Note that $D^\eta_u \Tilde{c}_A^i(u,z)=0$ for $|\eta|\geq2$. Hence, $D_u^\beta e^{-\Tilde{c}_A^i(u,z)}=e^{-\Tilde{c}_A^i(u,z)}\prod_{k=1}^R\left( -\partial_{u_k}\Tilde{c}_A^i(u,z) \right)^{\beta_k}$. Observe that

\begin{equation*}
\begin{aligned}
    \norm{\nabla_u(-\Tilde{c}_A^i(u,z))}&\leq 2\norm{ \left(\sum_j \int \Phi_R(x^\prime)d\mu^j(x^\prime) \right)\norm{z}^2 - 4A z } \\ 
    & \leq 2 M \kappa (\norm{z}^2 + 2\cdot2^{1/4}\sigma^{1/2}\norm{z}).
\end{aligned}
\end{equation*}
Consequently, we have 
\begin{equation*}
\begin{aligned}
    \left|D_u^\beta\int e^{\psi_0(z)-\Tilde{c}_A^i(u,z)}d\nu(z)\right| & \leq \int e^{\psi_0(z)-\Tilde{c}_A^i(u,z)} \prod_{k=1}^R\left|\partial_{u_k}(\Tilde{c}_A^i(u,z))\right|^{\beta_k} d\nu(z)  \\
    & \leq \int e^{\psi_0(z)-\Tilde{c}_A^i(u,z)} \norm{\nabla_u(\Tilde{c}_A^i(u,z))}^{|\beta|} d\nu(z)  \\
    & \leq  (2M\kappa)^{|\beta|}  \int e^{\psi_0(z)-\Tilde{c}_A^i(u,z)}(\norm{z}^2 + 2\cdot2^{1/4}\sigma^{1/2}\norm{z})^{|\beta|} d\nu(z) . 
\end{aligned}
\end{equation*}
We proceed to bound $\left| \frac{\int e^{\psi_0(z)-\Tilde{c}_A^i(u,z)}(\norm{z}^2 + 2\cdot2^{1/4}\sigma^{1/2}\norm{z})^{|\beta|} d\nu(z)}{\int e^{\psi_0(z)-\Tilde{c}_A^i(u,z)} d\nu(z)} \right|$ by splitting the integral in the numerator into $\norm{z} <\tau$ and $\norm{z} \geq \tau$, for some constant $\tau$ to be specified later. For $\norm{z} <\tau$, we have 
\begin{equation*}
    \left| \frac{\int e^{\psi_0(z)-\Tilde{c}_A^i(u,z)}(\norm{z}^2 + 2\cdot2^{1/4}\sigma^{1/2}\norm{z})^{|\beta|} d\nu(z)}{\int e^{\psi_0(z)-\Tilde{c}_A^i(u,z)} d\nu(z)} \right| \leq (\tau^2 + 2\cdot2^{1/4}\sigma^{1/2}\tau)^{|\beta|}.
\end{equation*}
For the region $\norm{z} \geq \tau$, consider the numerator, we have
\begin{equation*}
\begin{aligned}
    & \int_{\norm{z} \geq \tau} e^{\psi_0(z)-\Tilde{c}_A^i(u,z)}\norm{z}^{|\beta|}(\norm{z} + 2^{5/4}\sigma^{1/2})^{|\beta|} d\nu(z) \\
    \leq & \, e^{3\sqrt{2}M\kappa^2\sigma} \int_{\norm{z} \geq \tau} e^{4M\kappa^2\norm{z}^2+8\cdot2^{1/4}M\kappa^2\sigma^{1/2}\norm{z}} \norm{z}^{|\beta|}(\norm{z} + 2^{5/4}\sigma^{1/2})^{|\beta|} d\nu(z) \\
    \leq &\, e^{3\sqrt{2}M\kappa^2\sigma+4\cdot2^{1/4}M\kappa^2\sigma^{1/2}} \left( \int_{\norm{z} \geq \tau} e^{8M\kappa^2(1+2^{1/4}\sigma^{1/2})\norm{z}^2} d\nu(z) \right)^{1/2} \\
    & \quad \times \left(\int_{\norm{z}\geq \tau}  \norm{z}^{2|\beta|}(\norm{z} + 2^{5/4}\sigma^{1/2})^{2|\beta|} d\nu(z) \right)^{1/2} \\
    \leq & \, \sqrt{2} e^{3\sqrt{2}M\kappa^2\sigma+4\cdot2^{1/4}M\kappa^2\sigma^{1/2}+16M^2\kappa^4(1+2^{1/4}\sigma^{1/2})^2\sigma^2} \\
    & \quad \times \left(\int_{\norm{z}\geq \tau}  \norm{z}^{2|\beta|}(\norm{z} + 2^{5/4}\sigma^{1/2})^{2|\beta|} d\nu(z) \right)^{1/2},
\end{aligned}
\end{equation*}
where we again have used $2\norm{z}\leq 1+\norm{z}^2$ and Cauchy–Schwarz inequality for the second inequality, and $\int e^{t\norm{z}^2}d\nu(z) \leq 2e^{t^2\sigma^2/2}$ for a $\sigma^2$-sub-Gaussian random variable $\norm{z}^2$ for the third inequality. For the remaining integral, apply Cauchy–Schwarz inequality again, we have 
\begin{equation*}
\begin{aligned}
     & \int_{\norm{z}\geq \tau}  \norm{z}^{2|\beta|}(\norm{z} + 2^{5/4}\sigma^{1/2})^{2|\beta|} d\nu(z) \\
     \leq & \left(\int_{\norm{z}\geq \tau}  \norm{z}^{4|\beta|} d\nu(z)\right)^{1/2} \left(\int_{\norm{z}\geq \tau} (\norm{z} + 2^{5/4}\sigma^{1/2})^{4|\beta|} d\nu(z)\right)^{1/2} \\
    \leq & \left(\int_{\norm{z}\geq \tau}  \norm{z}^{4|\beta|} d\nu(z)\right)^{1/2} \left( \sigma^{2|\beta|} 2^{9|\beta|-1} \int_{\norm{z}\geq \tau} 1+\left(\frac{\norm{z}}{2^{5/4}\sigma^{1/2}} \right)^{4|\beta|}d\nu(z) \right)^{1/2} ,
\end{aligned}
\end{equation*}
with 
\begin{equation*}
\begin{aligned}
    & \left|\int_{\norm{z}\geq \tau}  \norm{z}^{4|\beta|} d\nu(z)\right| \leq e^{-\frac{\tau^4}{4\sigma^2}}\int e^{\frac{\norm{z}^4}{4\sigma^2}}\norm{z}^{4|\beta|} d\nu(z) \leq \sqrt{2}(2\sigma^2)^{|\beta|}e^{-\frac{\tau^4}{4\sigma^2}} \sqrt{(2|\beta|)!} \\
    & \left|\int_{\norm{z}\geq \tau} 1 d\nu(z)\right| \leq e^{-\frac{\tau^4}{4\sigma^2}}\int e^{\frac{\norm{z}^4}{4\sigma^2}} d\nu(z)  \leq  \sqrt{2}e^{-\frac{\tau^4}{4\sigma^2}} .
\end{aligned}
\end{equation*}
Substitute back and have
\begin{equation*}
\begin{aligned}
    &\int_{\norm{z}\geq \tau}  \norm{z}^{2|\beta|}(\norm{z} + 2^{5/4}\sigma^{1/2})^{2|\beta|} d\nu(z) \\
     \leq& e^{-\frac{\tau^4}{8\sigma^2}} 2^{5|\beta|-\frac{1}{4}}\sigma^{2|\beta|}((2|\beta|)!)^{\frac{1}{4}} \left( \left|\int_{\norm{z}\geq \tau} 1 d\nu(z)\right|+(2^{5/4}\sigma^{1/2})^{-4|\beta|}\left|\int_{\norm{z}\geq \tau}  \norm{z}^{4|\beta|} d\nu(z)\right|  \right)^{1/2} \\
    \leq & e^{-\frac{\tau^4}{4\sigma^2}} 2^{5|\beta|-\frac{1}{4}}\sigma^{2|\beta|}((2|\beta|)!)^{\frac{1}{4}}\left(\sqrt{2}+2^{-4|\beta|+\frac{1}{2}}((2|\beta|)!)^{\frac{1}{2}} \right)^{1/2}. 
\end{aligned}
\end{equation*}
Finally for the numerator, 
\begin{equation*}
\begin{aligned}
    & \int_{\norm{z} \geq \tau} e^{\psi_0(z)-\Tilde{c}_A^i(u,z)}\norm{z}^{|\beta|}(\norm{z} + 2^{5/4}\sigma)^{|\beta|} d\nu(z) \\ 
    \leq & \, e^{3\sqrt{2}M\kappa^2\sigma+4\cdot2^{1/4}M\kappa^2\sigma^{1/2}+16M^2\kappa^4(1+2^{1/4}\sigma^{1/2})^2\sigma^2-\frac{\tau^4}{8\sigma^2}} 2^{\frac{5}{2}|\beta|+\frac{3}{8}}\sigma^{|\beta|}((2|\beta|)!)^{\frac{1}{8}}\left(\sqrt{2}+2^{-4|\beta|+\frac{1}{2}}((2|\beta|)!)^{\frac{1}{2}} \right)^{1/4}. 
\end{aligned}
\end{equation*}
For the denominator, from Eq~\Cref{eq:upper_psi} we have 
\begin{equation*}
    \left|\int e^{\psi_0(z)-\Tilde{c}_A^i(u,z)} d\nu(z)\right|^{-1} = e^{\theta^i(x)} \leq e^{5\sqrt{2}M\kappa^2\sigma}.
\end{equation*}
Therefore, 
\begin{equation*}
\begin{aligned}
     &\left| \frac{\int e^{\psi_0(z)-\Tilde{c}_A^i(u,z)}(\norm{z}^2 + 2^{5/4}\sigma\norm{z})^{|\beta|} d\nu(z)}{\int e^{\psi_0(z)-\Tilde{c}_A^i(u,z)} d\nu(z)} \right| \\ \leq & (\tau^2 + 2^{5/4}\sigma^{1/2}\tau)^{|\beta|} + 2^{\frac{5}{2}|\beta|+\frac{3}{8}}\sigma^{|\beta|}((2|\beta|)!)^{\frac{1}{8}}\left(\sqrt{2}+2^{-4|\beta|+\frac{1}{2}}((2|\beta|)!)^{\frac{1}{2}} \right)^{1/4} \\& \times e^{8\sqrt{2}M\kappa^2\sigma+4\cdot2^{1/4}M\kappa^2\sigma^{1/2}+16M^2\kappa^4(1+2^{1/4}\sigma^{1/2})^2\sigma^2-\frac{\tau^4}{8\sigma^2}} .
\end{aligned}
\end{equation*}
Choose $\tau>0$ such that $\frac{\tau^4}{8\sigma^2} \geq 8\sqrt{2}M\kappa^2\sigma+4\cdot2^{1/4}M\kappa^2\sigma^{1/2}+16M^2\kappa^4(1+2^{1/4}\sigma^{1/2})^2\sigma^2$, more specifically, $$\tau^4=C(M\kappa^2(\sigma^{5/2}+\sigma^3)+M^2\kappa^4(\sigma^4+\sigma^{9/2}+\sigma^5))$$ for some constant $C >0$. Substituting back and applying AM-GM inequality, we have
\begin{equation*}
\begin{aligned}
    &\left| \frac{\int e^{\psi_0(z)-\Tilde{c}_A^i(u,z)}(\norm{z}^2 + 2^{5/4}\sigma\norm{z})^{|\beta|} d\nu(z)}{\int e^{\psi_0(z)-\Tilde{c}_A^i(u,z)} d\nu(z)} \right| \\
    \leq & C_\beta\left(\left( M\kappa^2(\sigma^{5/2}+\sigma^3)+M^2\kappa^4(\sigma^4+\sigma^{9/2}+\sigma^5)  \right)^{|\beta|/2} +\sigma^{|\beta|}\right)
\end{aligned}
\end{equation*}
thus
\begin{equation*}
\begin{aligned}
    &\left|D_u^\alpha \Tilde{\theta}^i(u) \right| \\ 
    \leq& C_\alpha (M\kappa)^{|\alpha|}\left(1+\sigma^{|\alpha|} \right) \left(1+\left(M\kappa^2(\sigma^{5/2}+\sigma^3)+M^2\kappa^4(\sigma^4+\sigma^{9/2}+\sigma^5)\right)^{|\alpha|/2}\right)  .
\end{aligned}
\end{equation*}

Finally, we are left to bound the derivatives of $\psi(z)$. The set up is exact as the derivatives of $\theta^i(x)$ and we will slightly abuse the multi-index notations. Similarly as before, observe that
\begin{equation*}
\begin{aligned}
    \norm{\nabla_z(-c_A^i(x,z))} &\leq 4 \norm{\Phi_R(x)} \norm{\left(\sum_j \int \Phi_R(x^\prime)d\mu^j(x^\prime) \right)}\norm{z} + 4\norm{A}_F\norm{\Phi_R(x)} \\
    & \leq 4M\kappa^2(\norm{z}+2^{1/4}\sigma^{1/2})
\end{aligned}
\end{equation*}
and 
\begin{equation*}
    \norm{\nabla^2_z (-c_A^i(x,z))}_{op} \leq 4M\kappa^2.
\end{equation*}
Note that the derivatives vanish with orders larger than 2. Consequently, apply the Faa di Bruno formula for the numerator, 
\begin{equation*}
\begin{aligned}
    \left| D_z^\beta\int e^{\theta^i_0(x)-c_A^i(x,z)}d\mu^i(x) \right| & \leq \int e^{\theta^i_0(x)-c_A^i(x,z)} \sum_{\text{partitions} \{B_1,\dots B_r'\} \text{ of }\beta}\prod_{b=1}^{r'} \left|D_z^{B_b}(-c_A^i(x,z)\right|d\mu^i(x) \\
    & \leq (4M\kappa^2)^{|\beta|}(1+\norm{z}+2^{1/4}\sigma^{1/2})^{|\beta|}\int e^{\theta^i_0(x)-c_A^i(x,z)} d\mu^i(x)
\end{aligned}
\end{equation*}

Observing cancellations with the denominator, we conclude that
\begin{equation*}
\begin{aligned}
    \left| D^\alpha_z\psi(z) \right| \leq  C_\alpha (M\kappa^2)^{|\alpha|}(1+\norm{z}+2^{1/4}\sigma^{1/2})^{|\alpha|}. 
\end{aligned}
\end{equation*}

\newpage
\section{Additional Experiment Results}
\label{sec:additiona_exp}


We conduct an ablation study on the latent dimension $d$, which is fixed to $32$ in Section~\ref{sec:experiments}. We evaluate \name and all baselines at $d\in\{32,48,64,128\}$, keeping every other component of the experimental protocol unchanged. The purpose of the ablation is to test whether the advantage of \name over the baselines is robust to the choice of the latent dimension. Tables~\ref{tab:coco6k_full} and~\ref{tab:coco9k_full} report the results on the MS--COCO 6K and MS--COCO 9K subsets, where the rows with $d=32$ coincide with Table~\ref{tab:coco6k}.

On the smaller 6K subset (Table~\ref{tab:coco6k_full}), the conclusions drawn at $d=32$ persist across the entire range. \name achieves the highest purity at every dimension and is the only method with positive silhouette scores, while all baselines yield negative silhouette values. It also obtains the best R@1 and R@5 scores at every dimension, with retrieval performance nearly unchanged for $d \in \{48,64,128\}$ and lower at $d=32$. These results suggest that \name learns a shared representation that is effective for retrieval and geometrically coherent, and that both properties are only mildly sensitive to the choice of $d$. 

On the larger 9K subset (Table~\ref{tab:coco9k_full}), the retrieval gains persist. \name achieves the best R@1 and R@5 scores at every dimension and again yields consistently positive silhouette scores. The purity results are more mixed in this setting, with CSA-based baselines obtaining higher purity at several dimensions. Thus, the 9K results suggest that \name gives the strongest bidirectional retrieval performance and the most consistently separated shared geometry, while class-purity scores on the larger subset remain competitive across methods.

\begin{table*}[ht]
\centering
\caption{Alignment and retrieval performances on the MS--COCO 6k subset.}
\begin{tabular}{c|l|cc|cc}
\hline
    & & \multicolumn{2}{c|}{Alignment} & \multicolumn{2}{c}{Retrieval} \\
\hline
    $d$ & \multicolumn{1}{c|}{Model} & \textbf{Purity} & \textbf{Silhouette} & \textbf{R@1} & \textbf{R@5} \\
\hline
    \multirow{4}{*}{32} & CSA & 0.638 & -0.210 & 0.257 & 0.542 \\
    & CSA+$\textit{STRUCTURE}$ & 0.654 & -0.196 & 0.290 & 0.595 \\
    & MLP+$\textit{STRUCTURE}$ & 0.638 & -0.156 & 0.154 & 0.390 \\
    & \name(Ours) & 0.719 & 0.117 & \textbf{0.380} & \textbf{0.674} \\
\hline
    \multirow{4}{*}{48} & CSA & 0.634 & -0.166 & 0.328 & 0.632 \\
    & CSA+$\textit{STRUCTURE}$ & 0.638 & -0.156 & 0.355 & 0.676 \\
    & MLP+$\textit{STRUCTURE}$ & 0.633 & -0.131 & 0.182 & 0.437 \\
    & \name(Ours) & 0.722 & 0.106 & \textbf{0.443} & \textbf{0.727} \\
\hline
    \multirow{4}{*}{64} & CSA & 0.632 & -0.140 & 0.359 & 0.666 \\
    & CSA+$\textit{STRUCTURE}$ & 0.650 & -0.131 & 0.386 & 0.707 \\
    & MLP+$\textit{STRUCTURE}$ & 0.634 & -0.118 & 0.199 & 0.467 \\
    & \name(Ours) & 0.713 & 0.086 & \textbf{0.450} & \textbf{0.725} \\
\hline
    \multirow{4}{*}{128} & CSA & 0.639 & -0.091 & 0.393 & 0.693 \\
    & CSA+$\textit{STRUCTURE}$ & 0.649 & -0.087 & 0.414 & 0.726 \\
    & MLP+$\textit{STRUCTURE}$ & 0.625 & -0.108 & 0.234 & 0.505 \\
    & \name(Ours) & 0.702 & 0.078 & \textbf{0.452} & \textbf{0.726} \\
\hline
\end{tabular}
\label{tab:coco6k_full}
\end{table*}

\begin{table*}[ht]
\centering
\caption{Alignment and retrieval performances on the MS--COCO 9k subset.}
\begin{tabular}{c|l|cc|cc}
\hline
    & & \multicolumn{2}{c|}{Alignment} & \multicolumn{2}{c}{Retrieval} \\
\hline
    $d$ & \multicolumn{1}{c|}{Model} & \textbf{Purity} & \textbf{Silhouette} & \textbf{R@1} & \textbf{R@5} \\
\hline
    \multirow{4}{*}{32} & CSA & 0.673 & -0.150 & 0.263 & 0.549 \\
    & CSA+$\textit{STRUCTURE}$ & 0.689 & -0.143 & 0.304 & 0.613 \\
    & MLP+$\textit{STRUCTURE}$ & 0.638 & -0.133 & 0.084 & 0.239 \\
    & \name(Ours) & 0.633 & 0.104 & \textbf{0.413} & \textbf{0.706} \\
\hline
    \multirow{4}{*}{48} & CSA & 0.687 & -0.115 & 0.328 & 0.629 \\
    & CSA+$\textit{STRUCTURE}$ & 0.676 & -0.116 & 0.358 & 0.683 \\
    & MLP+$\textit{STRUCTURE}$ & 0.636 & -0.107 & 0.128 & 0.350 \\
    & \name(Ours) & 0.633 & 0.112 & \textbf{0.460} & \textbf{0.757} \\
\hline
    \multirow{4}{*}{64} & CSA & 0.669 & -0.100 & 0.374 & 0.680 \\
    & CSA+$\textit{STRUCTURE}$ & 0.678 & -0.102 & 0.408 & 0.730 \\
    & MLP+$\textit{STRUCTURE}$ & 0.654 & -0.096 & 0.180 & 0.455 \\
    & \name(Ours) & 0.629 & 0.092 & \textbf{0.507} & \textbf{0.785} \\
\hline
    \multirow{4}{*}{128} & CSA & 0.683 & -0.067 & 0.405 & 0.728 \\
    & CSA+$\textit{STRUCTURE}$ & 0.699 & -0.076 & 0.442 & 0.756 \\
    & MLP+$\textit{STRUCTURE}$ & 0.643 & -0.081 & 0.304 & 0.615 \\
    & \name(Ours) & 0.620 & 0.070 & \textbf{0.506} & \textbf{0.785} \\
\hline
\end{tabular}
\label{tab:coco9k_full}
\end{table*}

\vskip 0.2in

\end{document}

%% file: macros.tex
\usepackage{xspace}

\usepackage{xcolor}

\newcommand{\name}{JK-EGW\xspace}

\newcommand{\norm}[1]{\left\lVert#1\right\rVert}

\newcommand{\innerproduct}[1]{\left\langle#1\right\rangle}

\newcommand{\vecone}{\mathbf{1}}

\newcommand{\BR}{\mathbb{R}}
\newcommand{\BE}{\mathbb{E}}

\newcommand{\OT}{\texttt{OT}}
